\documentclass[a4paper,10pt]{amsart}
\usepackage[utf8]{inputenc}
\usepackage{lmodern}
\usepackage[T1]{fontenc}
\usepackage{microtype} 	
\usepackage[all]{xy}
\usepackage{comment}
\usepackage{enumerate}
\usepackage{scalerel, stackengine}
\usepackage{graphicx}
\usepackage{soul}
\usepackage{amsmath,amssymb,amsfonts,amsthm,latexsym,mathrsfs,tikzsymbols}
\usepackage{mathtools}
\usepackage{ stmaryrd } 
\usepackage{xcolor}
\usepackage{hyperref}
\usepackage{cleveref}
\hypersetup{colorlinks=true,linkcolor=gray,citecolor=gray}
\usepackage[hyperpageref]{backref} 
\usepackage{booktabs}
\usepackage{longtable} 
\usepackage{pdflscape} 
\usepackage{array}

\newtheorem{lemma}{Lemma}[section]
\newtheorem{proposition}[lemma]{Proposition}
\newtheorem{theorem}[lemma]{Theorem}

\newtheorem{corollary}[lemma]{Corollary}

\newtheorem{maintheorem}{Theorem}

\newtheorem{maincorollary}[maintheorem]{Corollary}

\theoremstyle{definition}

\newtheorem{definition}[lemma]{Definition}

\theoremstyle{remark}
\newtheorem{remark}[lemma]{Remark}
\newtheorem{example}[lemma]{Example}

\newtheorem*{remark*}{Remark}
\newtheorem*{problem*}{Problem}

\numberwithin{equation}{section}

\newcommand\Ker{{\rm Ker}}

\newcommand\diag{{\rm diag}}
\newcommand{\GL}{\operatorname{GL}}
\newcommand{\SL}{\operatorname{SL}}

\newcommand{\Ma}{\operatorname{M}}

\DeclareMathOperator{\Aut}{Aut}

\newcommand\Hom{{\rm Hom}}

\newcommand\Span{\operatorname{Span}}

\newcommand\Irr{\operatorname{Irr}}

\DeclareMathOperator{\Res}{Res}
\DeclareMathOperator{\Ind}{Ind} 

\newcommand\cd{\operatorname{cd}}

\newcommand{\ZZ}{\mathcal{Z}}

\newcommand{\U}{{\mathcal U}}
\newcommand{\qa}[3]{\left(\frac{#1, #2}{#3}\right)}
\newcommand\PCI{\operatorname{PCI}}
\newcommand\Mexc{\operatorname{(M_{exc})}}
\newcommand\Cen{\operatorname{Cen}}

\newcommand\ot{\otimes}
\newcommand\op{\oplus}
\newcommand\mc{\mathcal}

\newcommand\wh{\widehat}
\stackMath
\newcommand\reallywidehat[1]{%
\savestack{\tmpbox}{\stretchto{%
  \scaleto{%
    \scalerel*[\widthof{\ensuremath{#1}}]{\kern-.6pt\bigwedge\kern-.6pt}%
    {\rule[-\textheight/2]{1ex}{\textheight}}
  }{\textheight}%
}{0.5ex}}%
\stackon[1pt]{#1}{\vstretch{1.5}{\tmpbox}}%
}
\DeclarePairedDelimiter\abs{\lvert}{\rvert}

\newcommand{\N}{{\mathbb N}}
\newcommand{\Z}{{\mathbb Z}}
\newcommand{\F}{{\mathbb F}}
\newcommand{\Q}{{\mathbb Q}}
\newcommand{\R}{{\mathbb R}}
\newcommand{\C}{{\mathbb C}}

 \newcommand\restr[2]{{
   \left.\kern-\nulldelimiterspace 
   #1 
   \right|_{#2} 
   }}

\theoremstyle{definition}

\renewcommand{\Mexc}{\ensuremath{\operatorname{(M_{exc})}}} 

\newcommand{\Fsq}{F}                      
\newcommand{\Fit}{\mathrm{Fit}}             
\newcommand{\Nor}{\operatorname{N}}       
\newcommand{\core}{\operatorname{core}}
\newcommand{\exc}{\operatorname{exp}}     
\newcommand{\cl}{\operatorname{cl}}       
\DeclareMathOperator{\ord}{o}             
\DeclareMathOperator{\Dic}{Dic}
\newcommand{\HH}{\mathbb{H}}            
\newcommand{\HHt}{\mathbb{H}_3}           
\newcommand{\HHrt}{\mathbb{H}_{\sqrt3}}   
\newcommand{\HHrtt}{\mathbb{H}_{\sqrt2}}  
\newcommand{\Lam}{\Lambda}                
\newcommand{\cp}[3]{(#1/#2,\,#3)}         
\newcommand{\Comp}{\mathcal{C}}           

\newcommand{\sg}[1]{\mathrm{SmallGroup}(#1)}
\newcommand{\sgid}[2]{[#1,#2]}            

\let\oldtocsubsection=\tocsubsection

\renewcommand{\tocsubsection}[2]{\hspace{2em}\oldtocsubsection{#1}{#2}}

\definecolor{Vino}{rgb}{0.256,0,0}

\usepackage[margin=1.4in]{geometry}

\thanks{The first author is grateful for financial support from the FWO and the F.R.S.–FNRS under the Excellence of Science (EOS) program (project ID 40007542). The second author is grateful to Fonds Wetenschappelijk Onderzoek Vlaanderen - FWO (grant 88258), and le Fonds de la Recherche Scientifique - FNRS (grant 1.B.239.22) for financial support.}

\begin{document}
\title{Representing in Low rank II: Virtual structure and Isomorphism problem}
\author{Robynn Corveleyn}
\author{Geoffrey Janssens}
\author{Doryan Temmerman}

\address{(Robynn Corveleyn) \newline Institut de Recherche en Math\'ematique et Physique, UCLouvain, Chemin du Cyclotron 2, 1348 Louvain-la-Neuve, Belgium \newline Email address: {\tt robynn.corveleyn@uclouvain.be}}

\address{(Geoffrey Janssens) \newline Département de mathématique, Université Libre de Bruxelles and \newline Department of Mathematics and Data Science, Vrije Universiteit Brussel,
Pleinlaan $2$, 1050 Elsene, Belgium
\newline Email address: {\tt geoffrey.janssens@ulb.be}}

\address{(Doryan Temmerman) \newline AI Lab, Vrije Universiteit Brussel, 
Pleinlaan 2, 1050 Elsene, Belgium\newline Email address: {\tt doryan.temmerman@vub.be}}

\begin{abstract}
We group theoretically characterise the finite groups $G$ and number fields $F$ such that all $F$-representations of reduced degree at least $2$ are so-called exceptional. If $G$ is not a $2$-group, we obtain a full classification. Secondly, we determine whether such groups are uniquely determined by their group algebra $FG$. In particular, we obtain an instance of the virtual structure problem of Kleinert and answer the rational isomorphism problem for such groups. 
\end{abstract}

\maketitle

\newcommand\blfootnote[1]{%
  \begingroup
  \renewcommand\thefootnote{}\footnote{#1}%
  \addtocounter{footnote}{-1}%
  \endgroup
}

\tableofcontents

\section{Introduction}\label{sec:intro}
\begingroup
\setlength{\parskip}{3pt}
Let $G$ be a finite group and $F$ a number field. The irreducible $F$-representations of $G$ are encoded in the simple components of the Wedderburn--Artin decomposition
\begin{equation}\label{eq:WA-intro}
FG\cong\Ma_{n_1}(D_1)\times\cdots\times\Ma_{n_\ell}(D_\ell),
\end{equation}
where $D_i$ are finite-dimensional division $F$-algebras. A recurrent theme in representation theory is to understand the interplay between the group $G$ and the isomorphism types of the simple components in \eqref{eq:WA-intro}. This can be approached from two sides. On the one hand, one can prescribe the simple components and ask which finite groups $G$ have this decomposition. A classical instance is the description by Sehgal \cite{SehNilpotent} of the groups $G$ such that $\Q G$ is a direct product of division algebras. On the other hand, one can ask how much of $G$ is remembered by the decomposition \eqref{eq:WA-intro}. For example, the rational group algebra determines an abelian group by a theorem of Perlis--Walker \cite{PW}. Nevertheless, Dade \cite{Dade} constructed non-isomorphic groups $G$ and $H$ with isomorphic group algebras over any field.

An actively studied point of view is that both perspectives meet in the unit group $\U(RG)$, where $R$ is the ring of integers of $F$. Indeed, $\U(RG)$ is a finitely presented group which has finite index in the direct product of the unit groups of the projections of $RG$ onto the simple components. Thus the first question asks for which finite groups $G$ the group $\U(RG)$ is, up to finite index, built from pieces of a prescribed kind. For the second, note that the commensurability class of $\U(RG)$ only depends on $FG$. Hence groups which are not distinguishable by their group algebra over $F$ will also not be distinguishable by the virtual structure of their unit groups.

This paper is the sequel to \cite{CJT1}, where the above philosophy was initiated for the representations of low rank. More precisely, \cite{CJT1} is concerned with the following type of simple algebras, whose presence often breaks the standard methods used to study $\U(RG)$.

\begin{definition}\label{def_exc_comps}
A finite-dimensional simple $F$-algebra $B$ is called \emph{exceptional} if it is isomorphic to one of the following:
\begin{enumerate}[(1)]
\item a non-commutative division algebra which is not a totally definite quaternion algebra over a number field,
\item a matrix algebra $\Ma_2(D)$ with $D\in\{\Q,\ \Q(\sqrt{-d}),\ \qa{-a}{-b}{\Q}\mid a,b,d\in\N_0\}$.
\end{enumerate}
If $B$ is of the latter form we call it an \emph{exceptional matrix algebra}, and in the former case an \emph{exceptional division algebra}. If $B\cong FGe$ for some $e\in\PCI(FG)$, then $B$ is called an \emph{exceptional component} of $FG$. Here $\PCI(FG)$ denotes the set of primitive central idempotents of $FG$.
\end{definition}

 We refer to the introduction of \cite{CJT1} for a detailed account of why these algebras are exceptional in practice. Let us only mention that the exceptional matrix algebras are exactly the matrix algebras having a maximal order $\Ma_n(\mc{O})$ such that $\SL_n(\mc{O})$ does not satisfy the congruence subgroup problem.
 
 The aim of this paper is to understand the finite groups $G$ and number fields $F$ such that $F G$ has the following property:
\[
\Mexc\quad\text{For all $e \in \PCI(FG)$ such that $F Ge\cong\Ma_n(D)$ with $n\ge2$, $\Ma_n(D)$ is exceptional.}
\]
So with respect to the standard methods mentioned above, they represent the “most degenerate” groups. In \cite{CJT1}, we characterised property $\Mexc$ in terms of topological and homological properties of $\U(R G)$, where $R$ is the ring of integers of a number field $F$.

In this paper, we solve the two questions mentioned in the first paragraph for this property. Namely, we characterise the finite groups $G$ such that their rational group algebra $\Q G$ has $\Mexc$, and we investigate how much of such a group $G$ is remembered by $\Q G$.

Firstly, we describe the structure of such groups and give a precise classification of the non-nilpotent ones. Together with the sequel \cite{CJT3}, which addresses the nilpotent case, this yields a new case of Kleinert's virtual structure problem. Secondly, we determine which of the non-nilpotent groups with $\Mexc$ are determined by their rational group algebra, i.e.\ satisfy the rational isomorphism problem. Somewhat to our surprise, not all of them do.

We now elaborate in more detail on the main results of this paper.

\subsection*{The virtual structure problem for exceptional components}

In \cite{KleinertSurvey}, Kleinert formulated the problem to obtain ‘unit theorems' for the unit group of orders $\mc{O}$ in finite-dimensional semisimple algebras $A$. More precisely, the problem is to determine the class of groups $\mc{G}(A)$ such that most of the generic unit groups of $A$ are contained in $\mathcal G(A)$. A generic unit group of $A$ is a subgroup of finite index in the group of reduced norm $1$ elements of an order in $A$. Formulated a bit differently, one could consider a class of groups $\mc{G}$ and classify all semisimple algebras $A$ such that $A$ has a unit theorem for $\mc{G}$.

If we restrict ourselves to group algebras, then the simple factors are related to each other via the representation theory of the group basis. As such the following problem was distilled in \cite{JesRioCrelle}.

\begin{problem*}[Virtual structure problem]
Let $\mc{G}$ be a class of groups. Classify the finite groups $G$ such that $\U(\Z G)$ has a subgroup of finite index lying in $\mc{G}$.
\end{problem*}

For a historical overview on the virtual structure problem, we refer the reader to the introduction of \cite{J1}. In all the cases solved so far, the answer is formulated in terms of the simple components of $\Q G$. In light of this, we formulated the following variant on the level of the simple components in \cite[Section 2]{CJT1}.

\begin{problem*}[Block virtual structure problem]
Let $\mc{P}$ be a property of a simple $F$-algebra. Classify the group algebras $FG$ such that, for each primitive central idempotent $e$ of $FG$ in some subset $\Omega$ \textup{(}depending on $\mc{P}$\textup{)} of $\PCI(FG)$, the simple algebra $FGe$ has property $\mc{P}$.
\end{problem*}

If $\mc{P}$ is the property that the unit group of an order in $FGe$ virtually belongs to a class of groups $\mc{G}$ which behaves well with respect to direct products, then this problem is equivalent to the virtual structure problem for $\mc{G}$, by \cite[Lemma 2.1]{CJT1}. However, for an arbitrary property $\mc{P}$, the virtual structure problem and its blockwise variant probably do not coincide.

In this paper, we consider the block virtual structure problem for $\mc{P}$ the property of being exceptional, and the subset $\Omega$ consisting of those $e\in\PCI(F G)$ such that $F Ge$ is not a division algebra. In other words, we consider the classification of the finite groups $G$ such that $F G$ has $\Mexc$.

Up to the presence of exceptional division components, this instance does correspond to a virtual structure problem. Indeed, in \cite[Theorem E]{CJT1} we proved that $F G$ has $\Mexc$ and no exceptional division components if and only if $\U(R G)$ has a subgroup of finite index which is a direct product of finitely generated abelian groups and of virtually indecomposable groups virtually mapping onto a non-abelian free group.

The classification obtained in this paper and its sequel \cite{CJT3} completes a line of research started 30 years ago with among others the papers \cite{JesLealRio,JesRioCrelle,PitaRioRuiz,FreebyFree}. The most general result until now concerned the finite groups $G$ such that $\Q G$ has the stronger property that
\begin{enumerate}[(i)]
\item all matrix components are of the form $\Ma_2(\Q(\sqrt{-d}))$ with $d\in\N$, and
\item $\Q G$ has no exceptional division components,
\end{enumerate}
which were classified in \cite{JesRioCrelle,FreebyFree} as the finite groups $G$ such that $\U(\Z G)$ is virtually a direct product of free-by-free groups. Considering property $\Mexc$ generalises the considerations of \emph{loc.\ cit.\ }in two ways. Namely, it also allows matrix algebras over quaternion algebras, as well as arbitrary division components. Furthermore, we consider group algebras over an arbitrary number field. Anecdotally, the MathSciNet review of \cite{JesRioCrelle} states that they obtain ``the most far-reaching result possible'', which indeed had seemed to be the case for by now almost 20 years.

It was proven in \cite[Section 3]{CJT1} (recalled in \Cref{Iff for Mexc over larger field} below), that classifying the groups $G$ and number fields $F$ such that $FG$ has $\Mexc$ reduces to (i) classifying the groups for which $\Q G$ has $\Mexc$ and (ii) obtaining the Wedderburn--Artin decompositions of those $\Q G$. Therefore, the main bulk of the paper is devoted to the case that $F= \Q$.

By \cite[Theorem 3.5]{SehNilpotent} the groups $G$ such that $\Q G$ has only division components are classified (see \eqref{condition only division}). Therefore we may assume that $\Q G$ has at least one matrix component. Our first main result describes the shape of $G$. Part of it was already obtained in \cite{CJT1}. We denote here and throughout the rest of the paper, the cyclic group of order $m$ by $C_m$.

\begin{maintheorem}[\Cref{broad form finite grps}, \Cref{necessary props 1,prop:fixed-inverted}, \Cref{necessary props 2} {\it \&} \Cref{thm:abelian-B}]\label{Main th A}
Let $G$ be a finite group such that $\Q G$ has a matrix component and satisfies $\Mexc$. Then 
\begin{enumerate}
    \item $G\cong A\rtimes Q$ where $A \cong C_3^n$, for some $n$, and $Q$ is a $2$-group,
    \item $Q$ has nilpotency class at most $3$,
    \item $[Q: \Cen_Q(A)] \leq 2$ and $\Fit(G)=A\times \Cen_Q(A)$,
\item $A=[A,Q]\times(A\cap\ZZ(G))$ and every element of $Q\setminus \Cen_Q(A)$ acts on $[A,Q]$ by inversion,
\item if $G$ is not nilpotent, then $\Fit(G)$ has an abelian subgroup of index at most $2$ which is normal in $G$.
\end{enumerate}
\end{maintheorem}

In particular $G$ is metabelian.  Furthermore, $Q$ is a Sylow $2$-subgroup of $G$ and $\Cen_Q(A)$ is the Sylow $2$-subgroup of $\Fit(G)$. We call $(Q,\Cen_Q(A))$ the \emph{Sylow pair} of $G$. A crucial ingredient for the remaining part of the paper is that \Cref{Main th A} implies that the isomorphism type of $G$ is determined by three invariants: the orders of $[A,Q]$ and of $A\cap\ZZ(G)$, and the Sylow pair $(Q,\Cen_Q(A))$ up to isomorphism of pairs (see \Cref{thm:abstract-iso}). Hence classifying the groups with $\Mexc$ boils down to determining which values of these three invariants occur.

Note that if $G$ is nilpotent, then $G=A\times Q$, and the problem reduces to a question on the $2$-group $Q$. In \Cref{thm:nilp-local-global} we prove that, for $Q$ non-abelian, $\Q Q$ has $\Mexc$ if and only if all the groups $Qe$, with $e\in\PCI(\Q Q)$, belong to an explicit list consisting of the cyclic $2$-groups and eleven groups of order at most $64$. The classification of the nilpotent groups with $\Mexc$ builds on this result and is the subject of the sequel \cite{CJT3}.

In the present paper we classify the non-nilpotent groups $G$ as above. This is done case-by-case for each possible value of the nilpotency class of $Q$ in \Cref{secrtion full characterisation solv}. Altogether, we obtain the following.

\begin{maintheorem}[\Cref{thm:mixed,classif Q abelian,thm:class2-abelian-kernel,thm:class2-nonabelian-list,thm:class3}]\label{Main th B}
Let $G$ be a finite non-nilpotent group. Then $\Q G$ has $\Mexc$ if and only if 
\begin{enumerate}
    \item $G=A\rtimes Q$ with $A \cong C_3^n$, for some $n$, and $Q$ is a $2$-group,
    \item $[Q:\Cen_Q(A)]=2$,
\end{enumerate}
and either $\Cen_Q(A)$ is elementary abelian, or $3\nmid|\ZZ(G)|$ and the pair $(Q,\Cen_Q(A))$ is isomorphic to one of the following pairs, for some $m\ge0$:
\begin{enumerate}[(i)]
\item $(\langle y\rangle\times E,\ \langle y^2\rangle\times E)$ with $E$ an abelian group such that $\exp(E) \mid 4$,
\item $(Q_8\times C_2^m,\ C_4\times C_2^m)$ or $(R\times C_2^m,\ C_R\times C_2^m)$, where 
$$R=\langle r,v\mid r^4=v^4=1,\ r^v=r^{-1}\rangle \text{ and } C_R=\langle rv,r^2\rangle,$$
\item $(N(i)\times C_2^m,\ C(i)\times C_2^m)$ for one of the seven pairs $(N(i),C(i))$ listed in \Cref{thm:class2-nonabelian-list},
\item $(SD_{16}\times C_2^m,\ D_8\times C_2^m)$ or $(Q_{16}\times C_2^m,\ Q_8\times C_2^m)$,
\end{enumerate}
\end{maintheorem}

Thus, maybe surprisingly, besides the two infinite families given by $\Cen_Q(A)$ elementary abelian and by $Q$ abelian, there are up to elementary abelian direct factors of $Q$ only eleven Sylow pairs.

For each of the groups in \Cref{Main th B} we also describe the Wedderburn--Artin decomposition of $\Q G$ in terms of group theoretical data of $G$, see \Cref{subsec:decompositions}. Finally, in \Cref{Classification other fields} we determine the fields $L\ne\Q$ such that $LG$ has $\Mexc$. Such fields only exist if $\Cen_Q(A)$ is elementary abelian or $Q$ is abelian, and they are imaginary quadratic.

Obtaining \Cref{Main th B} turned out to be a lot more involved than we expected. It required an especially good understanding of the strong Shoda pairs \cite{OlRiSi} of $G$. This has been possible thanks to the fact that  $G$ has an abelian normal subgroup $B$ of index at most $4$, see \Cref{Main th A}. In \Cref{thm:exhaustion} we show that, when $3\nmid|\ZZ(G)|$, the simple components of $\Q G$ on which $A$ acts non-trivially are parametrised by the $Q$-orbits of two explicit families of subgroups of $\Cen_Q(A)$. Each such component is a crossed product whose isomorphism type depends on relations between the latter subgroups and the elements of $Q\setminus \Cen_Q(A)$. In this way property $\Mexc$ is translated into a list of relations which have to hold in $Q$, see \Cref{thm:SSP-criterion}. Solving these relations is the main bulk of \Cref{secrtion full characterisation solv}.

Finally, we would like to discuss the relation to the milestone paper \cite{FreebyFree}. By \Cref{cor:index2}, if $G$ is not nilpotent, then $\U(\Z[\Fit(G)])$ is virtually a direct product of free-by-free groups. Thus the non-nilpotent groups with $\Mexc$ contain a subgroup of index $2$ classified in \cite{FreebyFree}. By \Cref{cor:which-Q} their Sylow $2$-subgroups are also among the groups classified in \emph{loc.\ cit}. It would therefore have been natural to start from the classification in \cite{FreebyFree} and to solve the associated extension problem. Nevertheless we decided to not depend on \cite{FreebyFree} and its predecessors, for the following reasons.
\begin{itemize}
\item In \cite{FreebyFree} the groups are described as the epimorphic images of certain families of groups given by a presentation. This is a full description, since every epimorphic image of such a group again satisfies the property. However, it is not a closed description up to isomorphism. The latter is what \Cref{Main th B} provides, through the three invariants above, and what is required for the rational isomorphism problem.
\item For the rational isomorphism problem we need the Wedderburn--Artin decomposition of $\Q G$ in terms of group theoretical data of $G$. This required the closed description and a precise understanding of the strong Shoda pairs (only partially obtained in \cite{FreebyFree}).
\item $\Mexc$ did not exclude the existence of exceptional division components, whereas being of Kleinian type does. Hence even when $\cd(G)=\{ 1,2 \}$ we still had to understand which extra groups arise compared to \cite{FreebyFree}. For instance, the groups in part (i) with $y$ of order at least $16$ have exceptional division components. Hence new types of simple components would have to be controlled along the extension.
\end{itemize}
In consequence, our methods also yield a new proof of the classification in \cite{FreebyFree} for non-nilpotent groups. We refer to Appendix \ref{appendix section B} for a detailed comparison between \Cref{Main th B} and \emph{loc.\ cit}.

\subsection*{The rational isomorphism problem}

Finally, we study the rational isomorphism problem, i.e.\ the problem of distinguishing a finite group among all finite groups through its rational group algebra.

\begin{problem*}[Rational isomorphism problem]
Let $G$ and $H$ be finite groups. If $\Q G\cong\Q H$ as $\Q$-algebras, is then $G\cong H$?
\end{problem*}

In general the answer is negative, as shown for the first time by the counterexample of Dade \cite{Dade}. Positive instances of the problem are scarce, as the authors are only aware of the case of abelian groups \cite{PW} and metacyclic groups \cite{GBdR1,GBdR2}. Note that Dade's examples are metabelian, indicating the difficulty to find new classes. Although finite groups satisfying $\Mexc$ are metabelian by \Cref{Main th A}, we show that the non-nilpotent ones are distinguishable up to a single family of exceptions. For this, we use the Wedderburn--Artin decompositions obtained along the proof of \Cref{Main th B}. Let $N(3)$ and $N(4)$ be two of the seven groups $N(i)$ (with \textsc{SmallGroupID}'s respectively $\sgid{32}{31}$ and $\sgid{32}{30}$) of \Cref{Main th B}.

\begin{maintheorem}[\Cref{rational Iso for Mexc}]\label{main th C}
Let $G=A\rtimes Q$ be a finite non-nilpotent group such that $\Q G$ has $\Mexc$, and let $H$ be a group such that $\Q H\cong\Q G$. Then the following holds: 
\begin{enumerate}
\item If $Q$ is not isomorphic to $N(3)\times C_2^m$ or $N(4)\times C_2^m$ for any $m\ge0$, then $H\cong G$.
\item If $Q\cong N(i)\times C_2^m$ with $i\in\{3,4\}$, then $H\cong G$ or $H\cong A\rtimes(N(7-i)\times C_2^m)$. 
\end{enumerate}
Moreover, $A\rtimes(N(3)\times C_2^m) \ncong A\rtimes(N(4)\times C_2^m)$  but their rational group algebras are isomorphic.
\end{maintheorem}

The smallest groups in the second case are $C_3\rtimes N(4)$ and $C_3\rtimes N(3)$ (with \textsc{SmallGroupID}'s respectively $\sgid{96}{142}$ and $\sgid{96}{143}$). Their group algebras are isomorphic over every field of characteristic $0$, but their integral group rings are not, see \Cref{rem:RIP-nature}. Interestingly, by the discussion at the start of the introduction, the unit groups of their integral group rings are nevertheless commensurable. Thus the virtual structure problem is blind for the difference between these two groups. Together with the case that $\Q G$ has only division components, \Cref{main th C} yields the following.

\begin{maincorollary}[\Cref{cor:RIP-summary}]
Let $G$ be a finite group such that $\Q G$ has $\Mexc$ and let $H$ be a finite group with $\Q H\cong\Q G$. Then $G\cong H$ in each of the following cases:
\begin{enumerate}[(a)]
\item $\Q G$ has only division components,
\item $G$ is non-nilpotent and the Sylow $2$-subgroup of $\Fit(G)$ is abelian,
\item $G$ is non-nilpotent and its Sylow $2$-subgroups have nilpotency class three,
\item $G$ is non-nilpotent with Sylow $2$-subgroups of nilpotency class two, not isomorphic to $N(3)\times C_2^m$ or $N(4)\times C_2^m$.
\end{enumerate}
\end{maincorollary}

Furthermore, for every non-nilpotent group $G$ with $\Mexc$ the algebra $\Q G$ determines the isomorphism types of $\Fit(G)$, of $\ZZ(G)$ and of $G/G'$, as well as the nilpotency class and the exponent of the Sylow $2$-subgroups, see \Cref{cor:RIP-invariants}. 

\endgroup
\vspace{5pt}
\noindent \textbf{Acknowledgements.} This research was supported through the program “Research in Pairs” by the Mathematisches Forschungsinstitut Oberwolfach in $2021$, we warmly thank the institute for its ideal working conditions. GJ would also like to thank the university of Ljubljana for its hospitality from March until July $2025$, where the major part of this article has been written.

\vspace{5pt}
\noindent \textbf{Declaration on the usage of AI.}
This work, from the initial research phase until the present written article, was completed without the use of artificial intelligence at any stage.

\vspace{0,3cm}
\noindent {\it Conventions and notations.}
Throughout the full article all groups denoted by a Latin letter are finite groups. All orders will be understood to be $\Z$-orders. The letter $F$ denotes a field of characteristic $0$ in \Cref{section recollections}. From \Cref{subsec:Lambda-components} onwards, it denotes the field $\Q(\sqrt{-3})$. We also use the following notations:
\begin{itemize}
\item $\PCI(FG)$ for the set of primitive central idempotents of $FG$ and $\Comp(FG)=\{FGe\mid e\in\PCI(FG)\}$ for the set of isomorphism types of the simple components of $FG$.
\item If $FGe\cong\Ma_n(D)$ with $n\ge2$, then we call $FGe$ a matrix component and otherwise a division component.
\item The degree of a central simple algebra $A$ is $\sqrt{\dim_{\ZZ(A)}A}$.
\item $Ge$ for the image of $G$ in $FGe$, for $e\in\PCI(FG)$ (see \eqref{map associated to central id}).
\item $\widehat H=\frac{1}{|H|}\sum_{h\in H}h$ for a subgroup $H$ of $G$.
\item $\cd(G)$ for the set of degrees of the irreducible complex characters of $G$ and $\pi(G)$ for the set of prime divisors of $|G|$.
\item $\ZZ(G)$, $\ZZ_2(G)$, $G'$, $\gamma_3(G)$, $\Phi(G)$ and $\Fit(G)$ respectively for the centre, the second centre, the commutator subgroup, the third term of the lower central series, the Frattini subgroup and the Fitting subgroup of $G$. Furthermore $G^{\mathrm{ab}}=G/G'$, $\cl(G)$ is the nilpotency class and $\exc(G)$ the exponent of $G$.
\item $\Cen_G(X)$ and $\Nor_G(X)$ for the centraliser and the normaliser of $X$ in $G$.
\item The notation $g^t$ will denote the conjugation $t^{-1}gt$ and $[g,t]=g^{-1}g^t$.
\item For a group element $g\in G$, we denote its order by $\ord(g)$. For an integer $n\ge1$, we denote its multiplicative order modulo $m$ by $\ord_m(n)$.
\item We denote the non-negative integers by $\N$, and the positive integers by $\N_0$.
\item $\zeta_n$ for a complex primitive $n$-th root of unity and $\phi(n)$ for Euler's phi function.
\item $\HH:=\qa{-1}{-1}{\Q}$, $\HHt:=\qa{-1}{-3}{\Q}$ and $\mathbb H_{\sqrt d}:=\qa{-1}{-1}{\Q(\sqrt d)}$ for $d\in\{2,3\}$. These are totally definite quaternion algebras over their (totally real) centres.
\item $\cp{E}{E_0}{\alpha}$ for the crossed product $E\oplus Eu$ with $ue=\sigma(e)u$ and $u^2=\alpha$, associated to a quadratic extension of fields $E/E_0$ with non-trivial automorphism $\sigma$ and $\alpha\in E_0^*$. It is isomorphic to $\Ma_2(E_0)$ if and only if $\alpha$ is a norm of $E/E_0$, and a quaternion division algebra otherwise \cite[Corollaries 4.7.5 and 4.7.7]{GS}.
\item $C_n$ for the cyclic group of order $n$. Furthermore $D_{2n}$, $Q_{2^n}$ and $SD_{2^n}$ denote the dihedral, generalised quaternion and semidihedral group of the indicated order, $\Dic_{4n}$ the dicyclic group of order $4n$ and $X\circ Y$ a central product of $X$ and $Y$.
\item $\sgid{n}{i}$ for $\sg{n,i}$ in the \textsc{SmallGroups} library of GAP \cite{SmallGroups}.
\item A subgroup $N$ of a group $X$ is called cocyclic if $X/N$ is cyclic.
\item For a finite abelian group $X$: $\Omega_1(X)=\{x\in X\mid x^2=1\}$.
\item For an elementary abelian $3$-group $V$: $s(V):=(|V|-1)/2$, the number of subgroups of index $3$.
\end{itemize}

\section{Finite groups with only exceptional higher simple components}

This paper will classify the finite groups $G$ and number fields $F$ such that $FG$ has the following property: 
\begin{align}
\tag{$\mathrm{M}_{\mathrm{exc}}$} \text{ For all } e \in \PCI(FG),&    \text{ such that } F Ge \cong \Ma_n(D)  \text{ with } n \geq 2, \\ \notag FGe &\: \cong \Ma_n(D) \text{ is exceptional}. 
\end{align}

In this section we recall and introduce the main tools for describing the simple components in case $F = \Q$. Furthermore, we will show that $\Mexc$ is preserved by subgroups.

\subsection{Recollections on part I}\addtocontents{toc}{\protect\setcounter{tocdepth}{2}}
\label{section recollections} 
In this section, we recall some results from the companion paper \cite{CJT1} which are relevant throughout the rest of this article.

In \cite[Theorem 3.5]{SehNilpotent} it was proven that 
\begin{equation}\label{condition only division}
\begin{array}{c}
\Q G \text{ has only division algebra components} \\
\Updownarrow \\
G \cong A \times C_2^{n} \times Q_8 \text{ with } \left\lbrace \begin{array}{l}
    A \text{ is abelian}   \\
    |A| \text{ and } o_{|A|}(2) \text{ are odd} 
\end{array} \right.
\end{array}
\end{equation}
In particular, $G$ is nilpotent of class at most $2$. Hence from now we may assume that $\Q G$ has at least one matrix component. In \cite{CJT1}, see also \cite[Proposition 4.1]{CdRZ}, we showed that property $\Mexc$ heavily restricts the possible simple components appearing in the Wedderburn--Artin decomposition of $\Q G$. 

Let $\cd(G)$ denote the set of degrees of the irreducible characters of $G$. Furthermore, for a field $K$ of characteristic zero and $a,b\in K^*$, we denote the associated \emph{quaternion algebra} by $\qa{a}{b}{K}$. Recall that this is the $K$-algebra with basis $\{ 1,i,j,ij\}$ and $i^2=a$, $j^2=b$, $ji=-ij$.

\begin{theorem}[{\cite[Theorem 3.4]{CJT1}}]
\label{thm componenten QG mexc}
Let $G$ be a finite group and suppose that $\Q G$ has a matrix component. If $\Q G$ has $\Mexc$, then the following hold:
\begin{enumerate}
\item The $1 \times 1$ components of $\Q G$ are either fields or quaternion algebras. More precisely, the non-commutative possibilities are: 
\[\{\qa{\zeta_{2^t}}{-3}{\Q(\zeta_{2^t})}, \qa{-1}{-1}{\Q},  \qa{-1}{-3}{\Q}, \qa{-1}{-1}{\Q(\sqrt{2})}, \qa{-1}{-1}{\Q(\sqrt{3})} \mid t\in \N_{\geq 3}\}.\]
\item The only exceptional matrix components are
\[\{\Ma_2(\Q(\sqrt{-d})), \Ma_2(\qa{-1}{-1}{\Q}), \Ma_2(\qa{-1}{-3}{\Q}) \mid d = 0, 1, 2, 3\}.\]
\end{enumerate}
In particular, $\cd (G) \subseteq \{ 1,2,4\}$ and $\deg (\Q G e) \mid 4$ for every $e \in \PCI(\Q G)$.
\end{theorem}

We call a simple $\Q$-algebra \emph{admissible} if it is a field or one of the algebras listed in \Cref{thm componenten QG mexc}, i.e.\ if it can occur as a simple component of a group algebra $\Q G$ with $\Mexc$ and having a matrix component. \smallskip

The proof of \Cref{broad form finite grps} below requires a finer understanding of the quotient groups $Ge$, with $e\in \PCI(\Q G)$. In \cite{CJT1}, we proved the following.

\begin{lemma}[{\cite[Lemma 3.9]{CJT1}}]\label{Prop components preserved for quotient}
    For any finite group $G$, normal subgroup $N \leqslant G$ and $e \in \PCI(\Q G)$, the following holds:
    \[\mathcal{C}(\Q[G/N]) \subseteq \mathcal{C}(\Q G)\,  \text{ and } \, \mathcal{C}(\Q[Ge]) \subseteq \mathcal{C}(\Q G).\]
    Consequently, 
    \[\mathcal{C}(\Q G) =  \displaystyle\bigcup\limits_{e \in \PCI(\Q G)} \mathcal{C}(\Q[Ge]).\]
    In particular, property $\Mexc$ is inherited by quotients.
\end{lemma}

Considering the role of the ground field, the following was proven in \cite{CJT1}. 

\begin{theorem}[{\cite[Theorem 3.1]{CJT1}}]\label{remark on larger field}
        Let $G$ be a finite group and $F\subseteq L$ fields of characteristic $0$. If $LG$ has $\Mexc$, then $F G$ also has $\Mexc$. Moreover, if $LG$ has a non-division component, then $L = \Q(\sqrt{-d})$ for some $d \in \N$.
\end{theorem}

In particular, if $LG$ has $\Mexc$, then $\Q G$ has $\Mexc$. More precisely, \cite[Theorem 3.1]{CJT1} obtained a necessary and sufficient condition for when $LG$ has $\Mexc$ in terms of $\Q G$. This criterion will be recalled in \Cref{section : other fields}. This shows the importance of first understanding the case of rational group algebras. Most of our work will be dedicated to this. 

\subsection{Describing simple components over the rationals}\label{sectie background SSP}
We require some methods to construct primitive central idempotents of $\Q G$. These methods were first introduced by Olivieri-del R\'io-Sim\'on \cite{OlRiSi}, see also \cite[Chapter 3]{EricAngel1} for a good introduction. 

\subsubsection{Background on Strong Shoda pair theory}

Recall that if $ H\unlhd G$, then 
\[\wh{H}:= \frac{1}{|H|} \sum_{h \in H} h\]
is a central idempotent in $\Q G$. Now, we set $\epsilon(H,H) := \wh{H}$, and for a strict normal subgroup $K$ of $H$ define
\begin{equation}\label{Def epsilon idempotents}
\epsilon (H,K) =  \prod\limits_{M/K \in \mc{M}(H/K)} (\wh{K} - \wh{M}),
\end{equation}
where $ \mc{M}(H/K)$ denotes the set of the non-trivial minimal normal subgroups of $H/K$. The construction results in a central idempotent in $\Q H$. To obtain a central idempotent in $\Q G$, one then associates to $K \unlhd H$ the element 
\begin{equation}\label{Def e(G,H,K)}
e(G,H,K) = \sum_{t \in \mc{T}} \epsilon(H,K)^{t},
\end{equation}
where $\mc{T}$ is a right transversal of $\Cen_G(\epsilon(H,K))$ in $G$. The element $e(G,H,K)$ is central in $\Q G$ and is a primitive idempotent when $(H,K)$ is a {\it Strong Shoda pair} of $G$, ({\it SSP }in short). A tuple $(H,K)$ is called a strong Shoda pair when $K \leq H \unlhd N_G(K)$, $H/K$ is cyclic and a maximal abelian subgroup of $N_G(K)/K$, and the $G$-conjugates of $\epsilon(H,K)$ are orthogonal.

To a central idempotent $e$, there is the associated epimorphism 
\begin{equation}\label{map associated to central id}
\varphi_e: G \rightarrow Ge, \ \ g \mapsto ge.
\end{equation}
Note that $Ge$ is a finite subgroup of $\U(\Q G)e$.

The following is a combination of \cite[Proposition 3.4.1, Theorems 3.4.2 \& 3.5.5 and Problem 3.5.1]{EricAngel1} and \cite[Lemma 3.4]{JM}.

\begin{theorem}[\cite{OlRiSi}]\label{form of SSP idempotent result}
With notations as above, $e :=e(G,H,K)$ is a primitive central idempotent of $\Q G$ if $(H,K)$ is a strong Shoda pair. Moreover, in that case $\Cen_G(\epsilon(H,K)) \cong N_G(K)$, $\ker(\varphi_{e}) = \text{core}_G(K) = \bigcap_{g \in G} K^g$ and 
$$\dim_{\Q} \Q G e =  [G:H] [G:N_G(K)] \phi([H:K]). $$
Furthermore,
$\Q G e \cong \Ma_{[G:N_G(K)]}\left(\Q(\zeta_{[H:K]}) * N_G(K)/H \right)$ for some explicit crossing (see \eqref{def crossing}). In particular, $\deg (\Q Ge) = [G:H]$.
\end{theorem}

\begin{remark}\label{rem:H=G}
For a SSP $(H,K)$ it might happen that $H=G$. This occurs precisely when  $G' \leq K$, which in turn is equivalent to $\Q G e(G,H,K)$ being commutative, see \cite[Lemma 2.4]{JespersSun}. We will also use the following two facts:
\begin{itemize}
\item $\ker(\varphi_e) = \core_G(K)$ is the kernel of the character induced from a linear character of $H$ with kernel $K$ \cite[Lemma 5.11]{Isaacs94}. So a normal subgroup $N$ of $G$ has trivial image in $Ge(G,H,K)$ if and only if $N\leq K$. Otherwise we say that $e(G,H,K)$, or the component $\Q Ge(G,H,K)$, is \emph{non-trivial on $N$}. 
\item If f $(H,K_1)$ and $(H,K_2)$ are strong Shoda pairs with $H\trianglelefteq G$, then $e(G,H,K_1)=e(G,H,K_2)$ if and only if $K_1$ and $K_2$ are conjugate in $G$ \cite[Exercise 3.4.3]{EricAngel1}.
\end{itemize}
\end{remark}

We specify the structure of the crossed product of $\Q(\zeta_{[H:K]}) * N_G(K)/H$ mentioned in \Cref{form of SSP idempotent result}, following \cite[Remark 3.5.6]{EricAngel1}. For ease of notation, let $N := N_G(K)$, $m := [H:K]$ and $H/K := \langle y \rangle$. The \emph{twisting} is given by:
\begin{equation}\label{def twisting}
\alpha \colon  N/H \rightarrow \Aut (\zeta_m), \quad \overline{x} \mapsto \alpha_{\overline{x}}(\zeta_m)= \zeta_m^{j},
\end{equation}
where $j$ is determined by \[(yK)^x =(yK)^j, \quad \text{ with  } xH = \overline{x}.\]
Moreover, the \emph{crossing} is given by:
\begin{equation}\label{def crossing}
f\colon  N/H \times N/H \rightarrow \Q(\zeta_m)^*, \quad  (\overline{n_1}, \overline{n_2}) \mapsto \zeta_m^{\ell},
\end{equation}
with $\ell$ uniquely determined by \[t_{\overline{n_1}} \cdot t_{\overline{n_2}} = y^{\ell}\cdot k \cdot t_{\overline{n_1n_2}},\] for $k \in K$ and  $\overline{n_i} := t_{\overline{n_i}}H$, for a fixed left transversal $\mc{T} = \{t_{\overline{n_i}} \mid t_{\overline{n_i}}H \in N/H\}$ of $H$ in $N$.\smallskip

Throughout, for a strong Shoda pair $(H,K)$ of $G$ we use the notation
\begin{equation}\label{eq:SSP-notation}
e:=e(G,H,K),\qquad N:=\Nor_G(K),\qquad m:=[H:K],\qquad H/K=\langle hK\rangle .
\end{equation}
We record the centre of $\Q Ge$, which follows from \Cref{form of SSP idempotent result}, and we make the crossed product $\Q(\zeta_m)\ast N/H$ explicit when $N/H$ has order two.

For this we recall some facts concerning quaternion algebras. If $L=F(\sqrt a)$ is a quadratic extension of $F$ and $\alpha\in F^*$, then the associated cyclic algebra is
\begin{equation}\label{eq:cp}
\cp{L}{F}{\alpha}:=L\oplus Lu \text{ with } u^2=\alpha,\quad u\beta u^{-1}=\bar\beta,
\end{equation}
where $\beta\mapsto\bar\beta$ is the non-trivial automorphism of $L/F$. It can also be considered as the crossed product of $L$ by $\mathrm{Gal}(L/F)$ with twisting $\alpha$, and $\cp{L}{F}{\alpha}\cong\qa{a}{\alpha}{F}$. An element $\alpha\in F^*$ is called a \emph{norm} of $L/F$ if $\alpha=N_{L/F}(\beta)=\beta\bar\beta$ for some $\beta\in L^*$. The norms form a subgroup of $F^*$. By \cite[Proposition 2.6.7]{EricAngel1}, $\cp{L}{F}{\alpha}\cong\Ma_2(F)$ if and only if $\alpha$ is a norm of $L/F$, otherwise it is a division algebra, and $\cp{L}{F}{\alpha}\cong\cp{L}{F}{\beta}$ if and only if $\alpha\beta^{-1}$ is a norm of $L/F$. \smallskip

\begin{proposition}\label{prop:centre}
With the notation \eqref{eq:SSP-notation}, the twisting \eqref{def twisting} is faithful\footnote{In other words, $\ker(\alpha)=1$.} and
\begin{equation}\label{eq:centre}
\ZZ (\Q Ge)=\Q(\zeta_m)^{N/H},\qquad \dim_{\Q} \ZZ(\Q Ge)=\frac{\phi(m)}{[N:H]}.
\end{equation}
If $[N:H]=2$ with $N/H=\langle yH\rangle$ and $y^2 \equiv h^j \pmod K$, then
\[
\Q(\zeta_m)\ast N/H\cong\cp{\Q(\zeta_m)}{\Q(\zeta_m)^{\langle y\rangle}}{\zeta_m^{\,j}}.
\]  
\end{proposition}
\begin{proof}
If $nH$ acts trivially on $H/K$, then $nK \in \Cen_{N/K}(H/K)$. Since $H/K$ is a maximal abelian subgroup of $N/K$, we get $n\in H$. Thus the twisting is faithful, and by \cite[Lemma 2.6.1]{EricAngel1} the crossed product $\Q(\zeta_m)\ast N/H$ is a central simple algebra over the fixed field $\Q(\zeta_m)^{N/H}$, which has degree $\phi(m)/[N:H]$ over $\Q$, yielding \eqref{eq:centre}.

Next, let $[N:H]=2$ and consider for \eqref{def crossing} the transversal $T=\{1,y\}$. Then $\Q(\zeta_m)\ast N/H=\Q(\zeta_m)\oplus \Q(\zeta_m) u_{\bar y}$ where the element $u_{\bar y}$ acts on $\Q(\zeta_m)$ as $\alpha_{\bar y}$, the non-trivial automorphism of $\Q(\zeta_m)/\Q(\zeta_m)^{\langle y\rangle}$. Furthermore, $u_{\bar y}^2=f(\bar y,\bar y)=\zeta_m^{\,j}$ because $y\cdot y=h^jk$ with $k\in K$. In summary, we have \eqref{eq:cp} with $\alpha=\zeta_m^{\,j}$.
\end{proof}

Finally, we will mainly use the theory of SSP for $G$ metabelian, in which case the SSPs have a convenient description.

\begin{theorem}[{\cite[Theorem 3.5.12]{EricAngel1}}]
\label{SSP in metabelian}
Let $G$ be a finite metabelian group and let $B$ be a maximal abelian subgroup of $G$ containing $G'$. Let $K \leq G$ be such that $C' \leq K \leq C$ for some subgroup $C$ of $G$ with $B \leq C \leq G$. Then $(H,K)$ is a strong Shoda pair if and only if the following hold:
\begin{enumerate}
    \item $H/K$ is cyclic,
    \item $H$ is maximal in the set $\{ C \leq G \mid B \leq C \text{ and } C' \leq K \leq C\}$.
\end{enumerate}
\end{theorem}

\subsubsection{Criterion for a single component to be exceptional}

We now consider a single strong Shoda pair $(H,K)$ of a finite group $G$ and determine when the associated simple component is an admissible  matrix algebra.  Recall from \Cref{thm componenten QG mexc} that this means that it is in the following set:
\begin{equation}\label{eq:exceptional-list}
\{ \Ma_2(\Q),\ \Ma_2(\Q(i)),\ \Ma_2(\Q(\sqrt{-2})),\ \Ma_2(\Q(\sqrt{-3})),\Ma_2(\HH),\ \Ma_2(\HHt)\}.
\end{equation}
In particular, since $\deg\Q Ge=[G:H]$, the component $\Q Ge$ can only be in \eqref{eq:exceptional-list} if $[G:H]\in\{2,4\}$. We now characterise which pairs give which algebra, according to the position of $N_G(K)$ between $H$ and $G$. 

\begin{proposition}[Criterion for a single component]\label{conditions comp exc}
Let $G$ be a finite group with $\pi(G)\subseteq\{2,3\}$ and let $(H,K)$ be a SSP of $G$ such that $H\ne G$ and $\Q Ge$ is not a division algebra. Then, using the notation \eqref{eq:SSP-notation}, $\Q Ge$ is in \eqref{eq:exceptional-list} if and only if one of the following occur:
\begin{enumerate}[(i)]
\item $N=H$, $[G:H]=2$ and $m\in\{1,2,3,4,6\}$. Then 
$$\Q Ge\cong\Ma_2(\Q(\zeta_m)).$$
\item $N=G$, $[G:H]=2$, $m\in\{3,4,6\}$ and, for $y\in G\setminus H$, $y$ inverts $H/K$ and $y^2\in K$. Then 
$$\Q Ge\cong\Ma_2(\Q).$$
\item  $N=G$, $[G:H]=2$, $m\in\{8,12\}$, and, for $y\in G\setminus H$, writing $h^y\equiv h^q$ and $y^2\equiv h^j\pmod K$, one has $(m,q)\in\{(8,3),(8,5),(12,5),(12,7)\}$ and $\zeta_m^{j}$ is a norm of $\Q(\zeta_m)/\Q(\zeta_m)^{\langle y\rangle}$. Then 
$$\Q Ge\cong\Ma_2(\Q(\zeta_m)^{\langle y\rangle})$$ 
with $\Q(\zeta_m)^{\langle y\rangle}=\Q(\sqrt{-2})$, $\Q(i)$, $\Q(i)$, $\Q(\sqrt{-3})$ respectively.
\item $[G:N]=[N:H]=2$, $m\in\{4,6\}$, and, for $y\in N\setminus H$, $y$ inverts $H/K$ and $y^2\equiv h^{m/2}\pmod K$. Then 
$$\Q Ge\cong\Ma_2(\HH) \text{ if } m=4 \text{ and } \Q Ge\cong\Ma_2(\HHt) \text{ if }  m=6.$$
\item  $N=G$, $[G:H]=4$, $m\in\{8,12\}$, and $z^2\notin K$ for $z\in G$ inverting $H/K$. Then, with $P$ a Sylow $2$-subgroup of $G/K$,
\[
\Q Ge\cong\begin{cases}
\Ma_2(\HH) & \text{if } m=8,\\
\Ma_2(\HH) & \text{if } m=12 \text{ and } P\cong Q_8\times C_2 \text{ or } SD_{16},\\
\Ma_2(\HHt) & \text{if } m=12 \text{ and } P\cong C_4\circ D_8 \text{ or } Q_{16}.
\end{cases}
\]
\end{enumerate}
\end{proposition}

\begin{proof}[Proof of \Cref{conditions comp exc}]
Put $r=[G:N]$ and $u=[N:H]$, so that $ru=[G:H]\in\{2,4\}$.

First let $[G:H]=2$. Since $\Q Ge$ is not a division algebra, it is exceptional if and only if it is $\Ma_2(\Q)$ or $\Ma_2(\Q(\sqrt{-d})$. If $u=1$, then $r=2$ and $\Q Ge\cong\Ma_2(\Q(\zeta_m))$. Hence $\phi(m)\le2$, which gives (i). 

If $u=2$, then $r=1$, $N=G$ and \eqref{eq:centre} gives $[\Q(\zeta_m)^{\langle y\rangle}:\Q]=\phi(m)/2\in\{1,2\}$, so that $m\in\{3,4,5,6,8,10,12\}$. Since $\pi(G)\subseteq\{2,3\}$, in fact $m\in\{3,4,6,8,12\}$. By \Cref{prop:centre}, $\Q Ge\cong\cp{\Q(\zeta_m)}{\Q(\zeta_m)^{\langle y\rangle}}{\zeta_m^j}$, which is a matrix algebra exactly when $\zeta_m^j$ is a norm of $\Q(\zeta_m)/\Q(\zeta_m)^{\langle y\rangle}$.  Here $G/H=\langle yH\rangle$ acts faithfully on $H/K=\langle hK\rangle\cong C_m$ (\Cref{prop:centre}), i.e.\ $h^y\equiv h^q\pmod K$ for some $q\in(\Z/m\Z)^*$ of order two, and $y^2K$ is fixed by this action, as $y$ commutes with $y^2$.

For $m\in\{3,4,6\}$ one has $(\Z/m\Z)^*=\{\pm1\}$, so $y$ inverts $H/K$ and $\Q(\zeta_m)^{\langle y\rangle}=\Q$. The elements of $H/K$ fixed by inversion are $K$ and, for $m$ even, $h^{m/2}K$. Hence $\zeta_m^j=\pm1$, with $\zeta_m^j=1$ if $m=3$. The parameter $1$ gives $\Ma_2(\Q)$, whereas $-1$ gives the division algebra $\HH$ (if $m=4$) or $\HHt$ (if $m=6$). Thus we obtained (ii).\smallskip

For $m\in\{8,12\}$, $(\Z/m\Z)^*\cong C_2^2$ has three elements of order two, and $\Q(\zeta_m)^{\langle y\rangle}$ is the fixed field of $\zeta_m\mapsto\zeta_m^q$. For $m=8$ these are $q=3,5,7$ with $\Q(\zeta_m)^{\langle y\rangle}=\Q(\sqrt{-2}),\Q(i),\Q(\sqrt2)$ and for $m=12$ they are $q=5,7,11$ with $\Q(\zeta_m)^{\langle y\rangle}=\Q(i),\Q(\sqrt{-3}),\Q(\sqrt3)$. An exceptional algebra $\Ma_2(\Q(\zeta_m)^{\langle y\rangle})$ requires $\Q(\zeta_m)^{\langle y\rangle}$ to be imaginary quadratic, which excludes $q=-1$, i.e.\ $q=7$ for $m=8$ and $q=11$ for $m=12$. In the remaining cases $\Q Ge\cong\cp{\Q(\zeta_m)}{\Q(\zeta_m)^{\langle y\rangle}}{\zeta_m^j}$ is $\Ma_2(\Q(\zeta_m)^{\langle y\rangle})$ exactly when $\zeta_m^j$ is a norm of $\Q(\zeta_m)/\Q(\zeta_m)^{\langle y\rangle}$. In summary, we showed (iii).\smallskip

Now let $[G:H]=4$. Then the only option in \eqref{eq:exceptional-list} is $\Ma_2(D)$ with $D\in\{\HH,\HHt\}$.  So by \eqref{eq:centre} we need $\phi(m)=[N:H]=u$. There are three cases, $(r,u)\in\{(4,1),(2,2),(1,4)\}$. If $(r,u)=(4,1)$, then $\Q Ge\cong\Ma_4(\Q(\zeta_m))$, which is not exceptional. \smallskip 

If $(r,u)=(2,2)$, then $\phi(m)=2$, so $m\in\{3,4,6\}$, and $N/H=\langle yH\rangle$ acts faithfully on $\Q(\zeta_m)$, i.e.\ $y$ inverts $H/K$. Then $\Q Ge\cong\Ma_2(\cp{\Q(\zeta_m)}{\Q}{\zeta_m^j})$ with $y^2\equiv h^j\pmod K$ (\Cref{prop:centre}), and $y^2K$, being fixed by the inversion, is $K$ or, for $m$ even, $h^{m/2}K$. Thus $\zeta_m^j=\pm1$, with $\zeta_m^j=1$ if $m=3$. The parameter $1$ gives $\Ma_4(\Q)$, which is not exceptional; the parameter $-1$ gives $\Ma_2(\HH)$ for $m=4$ and $\Ma_2(\HHt)$ for $m=6$. This is (iv). \smallskip

If $(r,u)=(1,4)$, i.e.\ $N=G$, then $\phi(m)=4$, so $m\in\{5,8,10,12\}$, and $m\in\{8,12\}$ as $\pi(G)\subseteq\{2,3\}$. Since $\Q Ge\cong\Q[G/K]e(G/K,H/K,1)$, we may replace $G$ by $G/K$ and assume that $K=1$. Then $H=\langle h\rangle$ is cyclic of order $m$ and $\Q Ge$ is the simple component attached to the faithful irreducible character of $G$ of degree four. By \Cref{prop:centre} the group $G/H$  embeds into $\mathrm{Gal}(\Q(\zeta_m)/\Q)\cong(\Z/m\Z)^*\cong C_2^2$. Hence $\Cen_G(H)=H$ and $G/H\cong C_2^2$ acts on $H$ as its full automorphism group. In particular there are elements $z\in G$ inverting $h$, and they form a coset of $H$. For such $z$ the element $z^2\in H$ is fixed by $z$, so $z^2\in\{1,h^{m/2}\}$, and $(zh^j)^2=z^2$ for all $j$. Thus the condition $z^2\ne1$ does not depend on the choice of $z$.

Since $|G|=[G:H]\,|H|=4m$, the group $G$ has order $32$ or $48$. A search in the SmallGroups library \cite{SmallGroups} with GAP shows that there are exactly ten groups of these orders with a cyclic normal subgroup $H$ of order $m$ such that $\Cen_G(H)=H$ and $[G:H]=4$. Five of them have $z^2=1$, namely $\sgid{32}{43}$, $\sgid{48}{15}$, $\sgid{48}{17}$, $\sgid{48}{38}$ and $\sgid{48}{41}$. The five others have $z^2\ne1$. They are $\sgid{32}{44}$ and the groups $\sgid{48}{40}$, $\sgid{48}{16}$, $\sgid{48}{39}$ and $\sgid{48}{18}$, whose Sylow $2$-subgroups are $Q_8\times C_2$, $SD_{16}$, $C_4\circ D_8$ and $Q_{16}$ respectively.

Suppose that $\Q Ge$ is an exceptional matrix algebra. Then $Ge$ is one of the groups listed in the table of Appendix \ref{appendix section} (see also \cite{EKVG}). The five groups with $z^2=1$ do not occur in that list, so $z^2\neq 1$. Conversely, the five groups with $z^2\neq 1$ occur in the list, and by the table of \cite{CJT1} (see also \cite[Table 2]{EKVG}) their faithful component of degree four is $\Ma_2(\HH)$ for $\sgid{32}{44}$, $\sgid{48}{40}$ and $\sgid{48}{16}$, and $\Ma_2(\HHt)$ for $\sgid{48}{39}$ and $\sgid{48}{18}$. This finishes the proof of (v).
\end{proof}

\subsection{\texorpdfstring{Necessary conditions for $\Mexc$ groups}{Necessary conditions for Mexc groups}}
In this section we first give several necessary group theoretical properties for the finite groups $G$ such that $\Q G$ has $\Mexc$. 
Moreover we will show that there exist subgroups of index two such that the matrix components are as in \cite{FreebyFree}.

The proof of \Cref{thm componenten QG mexc} to describe the possible simple algebras which can appear as a component of $\Q G$ with $\Mexc$, was achieved by describing which finite groups can be isomorphic to $Ge$ for $e \in \PCI( \Q G)$. This has the advantage of moreover yielding a first set of interesting properties for $G$. 

\begin{theorem}\label{broad form finite grps}
Let $G$ be a finite group and suppose that $\Q G$ has a matrix component. If $\Q G$ has $\Mexc$, then 
\begin{enumerate}
\item $G$ has\footnote{The proof will furthermore show that $[Ge:Be] \mid 4$ for every $e \in \PCI(\Q G)$.} an abelian normal subgroup $B$ with $\exp (G/B) \mid 2$.
\item $\pi(G) \subseteq \{2,3 \}$.
\item $\exp( \mc{Z}(G) \cap G') \mid 2$.
\item if $\Q G$ has no exceptional division components, then 
$$\exp (\mc{Z}(G)) \mid 4 \text{ or }6 \text{ and } \exp(G) \mid 24.$$
\end{enumerate}
In particular, $G$ is metabelian with $\cd (G) \subseteq \{ 1,2,4\}$ and $\deg (\Q G e) \mid 4$ for every $e \in \PCI(\Q G)$.
\end{theorem}
In \Cref{SSP and prop sectie} we will show that several of the above properties in fact hold in a stronger form, e.g. $G/B$ is an elementary abelian $2$-group of rank at most $2$.

\begin{proof}[Proof of \Cref{broad form finite grps}.]
By \Cref{Prop components preserved for quotient}, the proof reduces to a study of the groups $Ge$, $e\in\PCI(\Q G)$. The set $\PCI(\Q G)$ decomposes into $\PCI_1:=\{e\mid \Q Ge\cong D\}$ and $\PCI_2:=\{e\mid \Q Ge\cong\Ma_2(D)\}$, $D$ a rational division algebra. For $e\in\PCI_2$ the group $Ge$ is a finite subgroup of $\mathrm{GL}_2(\Q(\sqrt{-d}))$ or $\mathrm{GL}_2\big(\qa{-a}{-b}{\Q}\big)$ spanning the full matrix algebra. These have been classified in \cite[Theorem 3.7]{EKVG} (with one group added in \cite[Appendix A]{BJJKT}) and are listed in \Cref{appendix section} together with the relevant invariants. For $e\in\PCI_1$ the group $Ge$ is a finite subgroup of a division algebra, classified by Amitsur \cite{Amitsur}.

The above groups $Ge$ having $\Mexc$ were obtained in the proof of \cite[Theorem 3.4]{CJT1}, whose conclusion we now summarize.

Firstly, the groups $Ge$ with $e\in\PCI_2$ are written in table of the Appendix. Next, the only non-abelian groups $Ge$ with $e\in\PCI_1$ and such that $\Q[Ge]$ has $\Mexc$ are :
$$Q_8, \,  Q_{16},\, C_3 \rtimes Q_8, \, C_m\times Q_8 \text{ and } C_3\rtimes C_{2^n}$$
where $C_{2^n}$ acts by inversion on $C_3$ and the numbers $m$ and $\ord_m(2)$ are odd. 

Recall that for $n\geq 4$ one has
\begin{multline}\label{decomp grp with exp comp}
\Comp(\Q[C_3\rtimes C_{2^n}])=\Big\{\Q(\zeta_{2^\ell}),\ \HHt,\ \Ma_2(\Q),\ \Ma_2(\Q(i)),\\ \qa{\zeta_{2^t}}{-3}{\Q(\zeta_{2^t})}\ \Big|\ 1\le\ell\le n,\ 3\le t\le n-1\Big\}.
\end{multline}
Furthermore, 
\begin{align*}
\mc{C}(\Q[C_3 \rtimes C_{2}]) & = \{\Q, \Ma_2(\Q) \}, \qquad \mc{C}(\Q[C_3 \rtimes C_{4}]) = \{\Q, \Q(\sqrt{-1}), \qa{-1}{-3}{\Q}, \: \Ma_2(\Q)\},\\
\mc{C}(\Q[C_3 & \rtimes C_{8}]) = \{\Q, \Q(\sqrt{-1}), \Q(\zeta_8), \qa{-1}{-3}{\Q}, \: \Ma_2(\Q), \: \Ma_2(\Q(\sqrt{-1}))\}.
\end{align*}

Now, part (1) and (2) follows by using that $G$ embeds in the group $\prod_{e \in \PCI(\Q G)} Ge$. Indeed, note that all the aformentioned groups have an abelian normal subgroup $B_e$ with $Ge/B_e$ isomorphic to $C_2$ or $C_2 \times C_2$. Hence \[B := G \cap \hspace{-10pt} \prod_{e \in \PCI(\Q G)} \hspace{-10pt}B_e\] is an abelian normal subgroup of $G$ with $\exp(G/B) \mid 2$. Consequently, $G$ is metabelian. Additionally, for the simple algebras $\Ma_n(D)$ allowed by $\Mexc$, we see that $\Ma_2(D) \ot_{\mc{Z}(D)} \C$ is isomorphic to $\Ma_2(\C)$ or $\Ma_4(\C)$. By the first part of the proof, if $\Q Ge \cong D$, then $D \ot_{\Z(D)} \C$ is either $\C$ or $\Ma_2(\C)$. Hence indeed $\cd(G) \subseteq \{ 1,2,4 \}$ and $\deg (\Q Ge) \mid 4$ for all $e \in \PCI(\Q G)$.\vspace{0,2cm}

Since $(\chi(1),q) = 1$ for any odd prime $q$ and any $\chi \in \Irr_{\C}(G)$ it follows from the Ito--Michler theorem \cite[Corollary 12.34]{Isaacs94} that $G \cong A \rtimes Q$ with $A$ an abelian subgroup of odd order and $Q$ a Sylow $2$-subgroup of $G$. \vspace{0,2cm}

Note that it also follows from the description of the groups $Ge$ that $\pi(Ge) \subseteq \{2,3 \}$ for all $e \in \PCI(\Q G)$ such that $Ge$ is non-abelian and not isomorphic to $C_m \times Q_8$. In the latter case $\Q Ge \cong \qa{-1}{-1}{\Q(\zeta_m)}$. Therefore, if we take $g \in G$ with $o(g)$ prime different from $2$ or $3$, then $ge \in \Ker(\varphi_e)$ except for primitive central idempotents $e$ such that $\Q Ge$ is a field or a quaternion algebra with $ge \in \mc{Z}(\Q Ge)$. It follows that $g \in \mc{Z}(G)$. Now, as $G \cong A \rtimes Q$ with $Q$ a $2$-group, we can decompose $G$ into $\langle g \rangle \times H$ for some $H \leq G$. However by the assumption that $\Q G$ has a matrix component, which is necessarily exceptional, we obtain hence that $o(g) \mid 4 \text{ or } 6$, a contradiction.

Notice that the preceding argument proved that:
\begin{equation}\label{One degenerate case doesn't happen}
C_m \times Q_8 \text{ can only occur as a quotient of } G \text{ if } m= 3.
\end{equation}
However we also had the assumption that $o_m(2)$ is odd, which excludes the option $m=3$. \vspace{0,2cm}

Next we consider statement (3) and (4). First we note that the statements can be reduced to the same statement for the projections $Ge$ such that $\Q Ge$ is not a field. 

Indeed, take $1 \neq g \in G$ and let $\lambda: \langle g \rangle \rightarrow \Q(\zeta_{o(g)})$ be a faithful linear $\Q$-representation of $\langle g \rangle$. Subsequently, consider a decomposition into simple $\Q$-representations of the induction $\Ind_{\langle g \rangle}^{G}(\lambda) = \sigma_1 \oplus \cdots \oplus \sigma_{\ell}$.  By Frobenius' reciprocity, 
$$
\Hom_{\Q \langle g\rangle}(\lambda, \Res^G_{\langle g \rangle}(\sigma) ) \cong \Hom_{\Q G}(\Ind_{\langle g \rangle}^G(\lambda), \sigma)
$$
for any $1 \leq  i \leq \ell$. Hence $\Res^G_{\langle g \rangle}(\sigma)$ contains $\lambda$ as simple sub-representation and thus $\Res^G_{\langle g \rangle}(\sigma)$ is also injective. Furthermore, recall that $\ker(\Ind_{\langle g \rangle}^G(\lambda)) = \text{core}_G(\ker(\lambda))$ and therefore $\Ind_{\langle g \rangle}^G(\lambda)$ is a faithful representation of $G$. Now, one of these $\sigma_i$ must be satisfy that $\sigma_i(G)$ is non-abelian as otherwise $G' \subset \ker(\Ind_{\langle g \rangle}^G(\lambda)) = 1$ and hence $G$ abelian, contradicting the assumption that $\Q G$ has a matrix component. In conclusion, we have shown that there exists an irreducible representation $\sigma_j$ such that $ker(\sigma_j) \cap \langle g\rangle=1$ and $\Span_{\Q}\{ \sigma_j(G)\}$ is not a field. Hence 
\begin{equation}\label{order in terms of components}
o(g) \mid \max \{\exp(Ge) \mid \Q Ge \text{ is non-commutative}\}.
\end{equation}
Also note that if $g \in \ZZ(G)$ (resp. in $G'$), then $ge \in \ZZ (Ge)$ (resp. in $(Ge)'$). Thus in that case one can lower the upper-bound in \eqref{order in terms of components} with the exponent of the centres (resp. commutator subgroups).

Next let $\mc{E}_{nc} := \{ e \in \PCI(\Q G) \mid \Q G e \text{ non-commutative} \}$. Observe that the explanation above also shows that 
$$i: G \hookrightarrow \prod_{e\in \mc{E}_{nc}} Ge : g \mapsto (ge)_e$$
is an embedding. On the other hand $G$ maps surjectively on $Ge$. Combined these facts yield that $i(G) \cap \prod_{e \in \mc{E}_{nc}} \ZZ(Ge) = i(\ZZ(G)) \cong \ZZ(G)$. Therefore $\exp(G/\ZZ(G))$ divides $\text{lcm}\{ \exp(Ge/\ZZ(Ge)) \mid e \in \mc{E}_{nc} \}$.

Now consider the set 
$$\mc{E} := \{ e\in \PCI(\Q G) \mid \Q Ge \text{ is not a division component }\}.$$
As shown in the table of \Cref{appendix section}: $\exp(\mc{Z}(Gfe) \mid 4 \text{ or } 6$ for $e\in \mc{E}$. The same holds for $\exp( \mc{Z}Ge)$ when $\Q Ge$ is a non-commutative division component different from $\qa{\zeta_{2^t}}{-3}{\Q(\zeta_{2^t})}$. Indeed, verifying all $Ge$ which are $\Mexc$ one notices that such components either occur from a $Ge$ in the table or when $Ge \cong Q_{16}$ or $C_3 \rtimes Q_8$. In all those cases $\exp(\mc{Z}(Ge)) \mid 4 \text{ or } 6$. In all the cases mentioned one can also manually verify the facts that $\exp(\ZZ(Ge)\cap (Ge)') \mid 2, \, \exp(Ge)\mid 24$. and $\exp(Ge/\ZZ(Ge)) \mid 12$. Finally, the division component $\qa{\zeta_{2^t}}{-3}{\Q(\zeta_{2^t})}$ only arises for $Ge \cong C_3 \rtimes C_{2^n}$ with the action by inversion. For that group $\exp(\ZZ(Ge) \cap (Ge)') = 1$ and $\exp(G/\ZZ(G)) = 6$.

In conclusion, these numerical facts combined with the observations made earlier yield both statements (3) and (4).
\end{proof}

\subsection{\texorpdfstring{Behaviour of $\Mexc$ under subgroups}{Behaviour of Mexc under subgroups}}

In \Cref{Prop components preserved for quotient} we have seen that property $\Mexc$ is preserved under quotients. It also closed under taking subgroups as the following show.

\begin{proposition}\label{Mexc under subgroups}
Let $G$ be a finite group and $F$ a field of characteristic $0$ such that $FG$ has $\Mexc$. If $H \leq G$, then $FH$ also has $\Mexc$.
\end{proposition}
\begin{proof}
For $e\in\PCI(FG)$ let $S_e:=\Span_F\{ He\}$ be the $F$-span of $He$ in $FGe$. The projection $FH\to S_e: \alpha\mapsto\alpha e$, is a surjective homomorphism of algebras, so $S_e$ is semisimple and its simple components are simple components of $FH$. Moreover, since $FH\to\bigoplus_eS_e$ is injective, as $FG=\bigoplus_eFGe$, every simple component of $FH$ occurs in some $S_e$. Hence it suffices to show that every non-division simple component $T$ of $S_e$ is an exceptional matrix algebra.\smallskip

If $FGe$ is a division algebra, then $S_e$ has no zero divisors and is therefore a division algebra. So let $FGe\cong\Ma_2(D)$ with $D\in\{\Q,\Q(\sqrt{-d}),\qa{-a}{-b}{\Q}\}$. Let $f$ be the central idempotent of $S_e$ with $T=S_ef$. Then $T\subseteq f\Ma_2(D)f$, which is the endomorphism ring of the right $D$-module $fD^2$ and hence isomorphic to $D$ or to $\Ma_2(D)$ according as $fD^2$ has $D$-dimension $1$ or $2$. In the former case $T$ is a division algebra. In the latter case $f=1$, so $T=S_e$ is a simple subalgebra of $\Ma_2(D)$. Write $T\cong\Ma_n(D')$ with $D'$ a division algebra. As $T$ contains $n$ pairwise orthogonal non-zero idempotents we have that $n\leq 2$. If $n=1$ then $T$ is a division algebra. If $n=2$, the two diagonal matrix units $f_1,f_2$ of $T$ are orthogonal idempotents of $\Ma_2(D)$ with $f_1+f_2=1$, hence of rank one, and $D'\cong f_1Tf_1\subseteq f_1\Ma_2(D)f_1\cong D$. Thus $D'$ is a division subalgebra of $D$ containing $\Q$. Therefore, if $D=\Q$ or $D=\Q(\sqrt{-d})$, then $D'=\Q$ or $D'=D$. If $D$ is a totally definite quaternion algebra over $\Q$, then $D'$ equals $\Q$, $D$ or a maximal subfield of $D$. The latter is an imaginary quadratic field since $D\ot_\Q\R$ is a division algebra. In all cases we obtained that $T\cong\Ma_2(D')$ is an exceptional matrix algebra, as desired.
\end{proof}

\section{\texorpdfstring{A semidirect product splittings in case of $\Mexc$}{A semidirect product splittings in case of (Mexc)}}

This section has two aims. The first being to obtain a first rough form of $\Mexc$ groups. Among others, \Cref{necessary props 1} and \Cref{prop:fixed-inverted} will show the following.

\begin{theorem}\label{summary th structure Mexc}
Let $G$ be a finite group such that $\Q G$ has $\Mexc$. Then the following holds:
\begin{enumerate}
    \item $G \cong A \rtimes Q$ with $Q$ a $2$-group and $A \cong C_3^k$ for some $k$,
    \item $A=A_{-1}\times (A\cap\ZZ(G))$ with $A_{-1} := \{ a \in A \mid a^{y}= a^{-1} \text{ for all } y \in G \setminus \Fit(G) \}$,
    \item $A_{-1}= [A,Q]$,
    \item $[G:\Fit(G)] = [Q: \Cen_Q(A)] \leq 2$,
\end{enumerate}    
\end{theorem}

From the form of $\Mexc$ groups obtained in \Cref{summary th structure Mexc} we can provide a neat list of invariants that determine a $\Mexc$ group. For this we introduce  following terminology.
    \begin{definition}\label{def: Sylow pair}
Let $G$ be a finite group such that $\Q G$ has $\Mexc$. Then
\begin{itemize}
    \item The tuple $(Q,\Cen_Q(A))$ is called the \emph{Sylow pair} of $G$.
    \item Two pairs $(Q_1,C_1)$ and $(Q_2,C_2)$ of $2$-groups will be called \emph{isomorphic} if there is an isomorphism $Q_1\to Q_2$ mapping $C_1$ onto $C_2$.
\end{itemize}
\end{definition}

Note that $Q$ and $\Cen_Q(A)$ are Sylow $2$-subgroups of $G$ and of $\Fit(G)$ respectively, explaning the terminology.  From \Cref{summary th structure Mexc} we deduce that the isomorphism type of $G$ is determined by its Sylow pair and the orders of $A_{-1}= [A,Q]$ and $A\cap \ZZ(G)$.

\begin{corollary}\label{thm:abstract-iso}
Let $G=A\rtimes Q$ and $H=A_H\rtimes Q_H$ be finite groups having $\Mexc$. Then,
\begin{equation}\label{eq:iso-data}
G\cong H \Longleftrightarrow \left\{
\begin{array}{l}
|[A,Q]|=|[A_H,Q_H]|, \\
|A\cap\ZZ(G)|=|A_H\cap\ZZ(H)|,\\
(Q,\Cen_Q(A))\cong(Q_H,\Cen_{Q_H}(A_H)).
\end{array}
\right.
\end{equation}
\end{corollary}
Due to this corollary the classification results of \Cref{secrtion full characterisation solv} will be formulated in terms of the Sylow pair.
\begin{proof}
An isomorphism $G\to H$ maps $A=O_3(G)$ onto $A_H=O_3(H)$. Consequently, it also sends $[A,G]=[A,Q]$ onto $[A_H,Q_H]$, $A\cap\ZZ(G)$ onto $A_H\cap\ZZ(H)$ and $\Fit(G)=A\times C$ onto $\Fit(H)=A_H\times C_H$. Therefore it induces an isomorphism from  $G/A\cong Q$ to $H/A_H\cong Q_H$, which maps $C$ onto $C_H$.

Conversely, suppose that \eqref{eq:iso-data} holds and let $\alpha\colon Q\to Q_H$ be an isomorphism with $\alpha(C)=C_H$. Since $A=[A,Q]\times(A\cap\ZZ(G))$ and $A_H=[A_H,Q_H]\times(A_H\cap\ZZ(H))$ are elementary abelian $3$-groups whose factors have the same orders, there is an isomorphism $\beta\colon A\to A_H$ with $\beta([A,Q])=[A_H,Q_H]$ and $\beta(A\cap\ZZ(G))=A_H\cap\ZZ(H)$. By \Cref{prop:fixed-inverted}, an element $q\in Q\setminus \Cen_Q(A)$ inverts $[A,Q]$ and centralises $A\cap\ZZ(G)$. The same holds in $H$. As $\alpha(\Cen_Q(A))=\Cen_{Q_H}(A_H)$, this gives $\beta(a^q)=\beta(a)^{\alpha(q)}$ for all $a\in A$ and $q\in Q$. Hence $aq\mapsto\beta(a)\alpha(q)$ is an isomorphism $G\to H$, as desired.
\end{proof}

The second aim of this section is to obtain properties of $2$-groups having $\Mexc$. Subsequentely, the focus will be on the Sylow $2$-subgroup of a non-nilpotent $\Mexc$ group. 

\subsection{Groups with low character degrees}
This section and the next aim to give a precise description of the groups satisfying $\Mexc$. In light of \Cref{broad form finite grps}, we first consider the more general class of groups having character degrees of the irreducible complex characters only divisors of $4$. We denote by $\cd(G)$ the set of character degrees of irreducible complex representations of $G$. In this section we focus on the groups $G$ with $\cd(G) = \{1,4\}$.

More generally, whenever $\cd(G) \subseteq \{1,p^{j},p^{k}\}$ then $G$ is solvable of derived length at most $3$ by \cite[Theorem 12.15]{Isaacs94}. Furthermore, $(\chi(1),q) = 1$ for any prime $q$ different from $p$ and any $\chi \in \Irr_{\C}(G)$. This allows to apply the Ito-Michler theorem \cite[Corollary 12.34]{Isaacs94} which in this case yields that 
    \begin{equation}
        \label{G is A semidir Q}
        G \cong A \rtimes Q, 
    \end{equation}
where $A$ is an abelian $p'$-subgroup and $Q$ a Sylow $p$-subgroup of $G$. More precisely, $A$ is the direct product of all Sylow $q$-subgroups of $G$, for all $q\neq p$, which are normal and abelian under the above assumption on $\cd(G)$.

Using character theory and \cite{IsaacsPassman}, the following results will give further restrictions on the decomposition (\ref{G is A semidir Q}) in the case that $\abs{\cd(G)}=2$. Note that for a non-nilpotent group, one always has that $\abs{\cd(G)} > 1$.

\begin{proposition}
\label{lower character degree classif}
Let $G$ be a non-nilpotent finite group such that $\cd(G) \subseteq \{ 1, p^{i}, p^k\}$ with $0 \leq i < k$ and $p^k = \max(\cd(G))$. Let $G \cong A \rtimes Q$ be the decomposition from \cref{G is A semidir Q}. If $Q$ is abelian, then the following hold:
\begin{enumerate}
\item The Fitting subgroup $F(G)$ of $G$ is the unique maximal abelian subgroup, and has index $p^k$,
\item $C_G(A) = F(G) = \mc{Z}(G) \times G'$,
\item $G$ decomposes as a semidirect product \[G \cong N \rtimes \langle C,Z \rangle\] with $N$ an abelian subgroup of $F(G)$ containing $A$ and with $C$ and $Z$ $p$-groups that are at most $2$-generated, and such that $Z$ is central in $G$. Moreover, ${G}/\langle N, Z \rangle$ is isomorphic to $C_{p^k}$ or $C_{p^{i}} \times C_{p^{k-i}}$.
\end{enumerate}
Furthermore, if $\cd(G) = \{1, p^k \}$ with $k \geq 2$, then $Q$ is abelian. Conversely, for every group $G$ with $\abs{\cd(G)} = 2$ which satisfies the properties (1) to (3), one has $\cd(G) = \{1, p^k\}$.
\end{proposition}

In \Cref{necessary props 1} we will prove that if $\Q G$ has $\Mexc$, then $[G:F(G)]=2$. Hence, under the $\Mexc$ assumption, by \Cref{lower character degree classif} the Sylow $2$-subgroup of $G$ can only be abelian if $\cd(G) = \{ 1,2 \}$.

\begin{proof}
    Suppose $G$ is as in the statement. We prove (1) and (2) simultaneously. Since $Q$ is assumed abelian, in particular $G$ is solvable of derived length at most $2$. In particular $G'$ is abelian. It follows from the work of Taunt, \cite[VI, Satz 14.7 b)]{Huppert} that $F(G) = \mc{Z}(G) \times G'$. 
    By\footnote{The fact that the Sylow $p$-subgroup $Q$ is abelian yields that $G$ is a $\{p\}$-character group as in \cite[Definition 2.1.]{IsaacsPassman}} \cite[Theorem 2.2 (ii)]{IsaacsPassman}, $C_G(A)$ is a normal abelian subgroup of $G$ of index $p^k$.
    Additionally, $C_G(A)$ and $G/C_G(A)$ are abelian and hence it follows from \cite[VI, Satz 14.7 a)]{Huppert} that 
    \begin{equation}\label{Fitting is centraliser} 
    C_G(A) = \mc{Z}(G) \times (C_G(A) \cap G') = \mc{Z}(G) \times G' = F(G).
    \end{equation}
    It now also immediately follows that $F(G)$ is the unique maximal \emph{abelian} subgroup and that $[G:F(G)] = p^k$.  Moreover, by \cite[Theorem 2.2]{IsaacsPassman}, ${G}/{C_G(A)}$ is either isomorphic to $C_{p^{k}}$ or to $C_{p^{i}} \times C_{p^{k-i}}$.\vspace{0,1cm}
    
    We denote \[F := C_G(A) = F(G).\] Let $x,y \in Q$ be such that $\langle xF, yF \rangle = G/F$. More precisely, we choose $x,y \in Q$ as follows. If $G/F \cong C_{p^{i}} \times C_{p^{k-i}}$, we choose $x$ such that $xF$ has order $p^{i}$ and $yF$ has order $p^{k-i}$. If $G/F \cong C_{p^k}$, we choose $x$ such that $xF$ has order $p^k$, and we let $y$ be an arbitrary power of $x$. Let then $p^{\ell}$ be the order of $xF$. Then $F \cap \langle x \rangle = \langle x^{p^{\ell}} \rangle$, and $F\cap \langle y \rangle = \langle y^{p^{k-i}}\rangle$.

    Let $\langle z_1 \rangle$ be the cyclic subgroup of $F$ containing $x^{p^{\ell}}$ and which is maximal amongst the cyclic subgroups of $Q \cap F$ for this property. Similarly let $\langle z_2 \rangle$ be the maximal cyclic $Q\cap F$ above $y^{p^{k-i}}$. Denote $Z = \langle z_1,z_2 \rangle$.
    
    As $Q$ is abelian and due to the maximality of $\langle z_1 \rangle$ and $\langle z_2 \rangle$, there exists a $M \leqslant Q$ such that $Q \cap F(G) \cong M \times Z$. Moreover, since $F(G)$ is abelian, \cref{G is A semidir Q} implies that 
   \[
   F(G) = A \times \left(Q \cap F(G)\right),
   \]
   and it follows that 
   \[F(G) = N \times Z, \text{ with } N = A \times M.\]
   Denoting $C = \langle x,y \rangle$, we now show that \[G \cong  N \rtimes \langle C, Z \rangle. \] 
   Firstly,  $N \cap \langle C,Z \rangle = \{1 \}$. Indeed, $\langle C, Z \rangle$ is abelian, and since $F(G) = N \times Z$, we have that $N \cap \langle C, Z\rangle$ is at most $C$. So suppose $t \in N \cap \langle x,y \rangle \leqslant F$. Then $t = x^n y^m$ for some $n,m \in \Z$. But then $x^n = ty^{-m} \in F$ and since $F\cap \langle x \rangle = \langle x^{p^{\ell}} \rangle$, we obtain that $x^n = x^{kp^{\ell}} \in Z$ for some $k \in \Z$. But then $x^n \in Z \cap N$, and hence $x^n = 1$. Thus $t = y^m \in F$. Since $F \cap \langle y \rangle = \langle y^{p^{k-i}} \rangle \leqslant Z$, we obtain again that $t = y^m \in Z\cap N$ and hence $t = 1$.

    Additionally, $N \cdot \langle C,Z \rangle = G$. Indeed,
   \begin{align*}
       \abs{N} \abs{\langle C, Z\rangle} & = \abs{N} \left[\langle C,Z \rangle : \langle x,Z \rangle \right]\left[\langle x,Z \rangle :  Z \right]\\
       & = \begin{cases} \abs{N} \frac{\abs{\langle y \rangle}}{\abs{\langle y^{p^{k-i}} \rangle }} \frac{\abs{\langle x \rangle}}{\abs{\langle x^{p^{i}}\rangle}}  \abs{Z}  & \text{ if } G/F \cong C_{p^{k-i}} \times C_{p^{i}} \\
       \abs{N} \frac{\abs{\langle x \rangle}}{\abs{\langle x^{p^k}\rangle}} \abs{Z} & \text{ if }  G/F \cong C_{p^k}\end{cases}\\
       & = \abs{N} \abs{Z} \cdot p^k = \abs{F(G)} [G : F(G)] = \abs{G}.
   \end{align*} 
   Finally, we show that $N$ is normal in $G$. Since $N = A \times M$, and $A$ is normal in $G$ by construction, it suffices to show that $M$ is normal in $G$. We show that in fact $Z \times M$ is central in $G$. Indeed, if $t \in Z \times M \leqslant F$, it may be written as $t = t_1 t_2$ for unique $t_1 \in \mc{Z}(G)$ and $t_2 \in G'$, by \cref{Fitting is centraliser}. Furthermore, since $Q$ is an abelian Sylow $p$-subgroup (and hence $G' \leqslant A$), one has that the order of $t_2$ is coprime to $p$. But by definition $t \in Q$ has order a power of $p$, and it follows that $\langle t \rangle = \langle t^{o(t_2)} \rangle \subseteq \mc{Z}(G)$, as claimed.

    Suppose $\cd(G) = \{1,p^k\}$. If $Q$ is non-abelian, then \cite[Exercise 12.6]{Isaacs94} implies that $G$ is nilpotent, which is in contradiction with the hypothesis. Hence $Q$ is abelian. 
    For the converse, as $[G: F(G)]= p^k$ and $F(G)$ is a normal abelian subgroup, Ito's theorem \cite[Theorem 6.15]{Isaacs94} yields that $\chi(1)$ divides $p^k$ for all irreducible characters $\chi$. Note also that $G' \leq N \leq F(G)$, hence by \cite[Lemma 1]{HM} there exists some $\chi \in \Irr_{\C}(G)$ with $\chi(1)= [G: F(G)]= p^k$. The assumption $|\cd(G)|=2$ concludes the proof.
\end{proof}

\begin{remark}
The proof of \Cref{lower character degree classif} yields extra restrictions on the structure of $G$. For example $F(G) = C_G(A) \cong Z(G) \times G'$ by (\ref{Fitting is centraliser}). Additionally, denoting $C = \langle x,y \rangle$ and using \cite[Theorem 3.1.(iii)]{IsaacsPassman}, one can prove that the action of $x^{i}$ on $G'$ is free, for $1 \leq i \leq p^k$.
\end{remark}

\begin{example}
Consider the group $G := C_5 \rtimes C_8 = \langle a,b \mid a^5,b^8, a^b= a^3 \rangle$. It is easily shown that $\cd(G)= \{1,4 \}$ and $F(G) = C_G(\langle a \rangle)= \langle a, b^4 \rangle$. Thus, in the notation of \Cref{lower character degree classif}, $N = \langle a \rangle$ and $C = \langle b \rangle$. In this example, a complement for $F(G)$ itself does not exist. Hence the splitting in \Cref{lower character degree classif} is the finest possible in general. 
\end{example}

Next we consider \emph{nilpotent} groups $G$ having $\Mexc$ and $\cd(G) = \{1,4\}$. Then the decomposition in (\ref{G is A semidir Q}) is a direct product $G \cong A \times Q$. A precise classification in the nilpotent case remains elusive. Nevertheless note that one has $\{ 1,4\}=\cd(G) = \cd(Q)$. Now applying \cite[Theorem 3.10 \& Lemma 5.4]{IsaacsPassman} yields the following.

\begin{lemma}\label{props for nilpot with cd 1,4}
    Let $G$ be a nilpotent group with $\cd(G) = \{ 1, 4\}$. Then $G \cong A \times Q$ with $A$ an odd abelian group and $Q$ a $2$-group satisfying:
    \begin{enumerate}
        \item $Q$ has nilpotency class $2$,
        \item $[Q,Q]$ and $Q/\mc{Z}(Q)$ are elementary abelian $2$-groups.
    \end{enumerate}
\end{lemma}

\begin{remark}
    If $\cd(G) = \{1,2,4 \}$, the statement of Lemma \ref{props for nilpot with cd 1,4} does not hold in general. For example the group 
    $$C_8 \rtimes (C_2 \times C_2) := \langle a,b,c \mid a^8 = b^2= c^2 =1, bab= a^3, cac = a^5, bc=cb \rangle $$
    has\footnote{This group has Small Group Id [32,43] and the stated information can be obtained via GAP or the GroupNames database.} $\cd(G) = \{1,2,4 \}$ but is nilpotent of class $3$ and $G' \cong C_4$.
\end{remark}

\subsection{Restriction on Sylow subgroups and action}

Recall from (\ref{G is A semidir Q}) that $G \cong A \rtimes Q$ with $A$ an odd abelian group and $Q$ a $2$-group. Denote by $\varphi:  Q \rightarrow \Aut (A)$ the action of $Q$ on $A$ and 
\begin{equation}\label{eq:C}
C:=\Ker(\varphi)=\Cen_Q(A).
\end{equation}
We will now give two of the main necessary properties on $Q, A$ and the action. 

\begin{proposition}\label{necessary props 1}
Let $G$ be a finite group such that $\Q G$ has $\Mexc$, with at least one matrix component. Then the following hold:
\begin{enumerate}
    \item $G \cong A \rtimes Q$ with $A$ an elementary abelian $3$-group and $Q$ a $2$-group,
    \item $[G: F(G)] = [Q: \Ker (\varphi)]$ divides $2$,
\end{enumerate}
\end{proposition}

The fact that $[G : F(G)] \mid 2$, together with the form of $G$ given by (1), can be interpreted very concretely, see \Cref{prop:fixed-inverted} below. Also, in \Cref{thm:abelian-B} we will show that, for non-nilpotent $G$, the subgroup $\Ker (\varphi)$ has an abelian subgroup of index at most $2$ which is normal in $Q$. In other words, $G$ has an abelian normal subgroup $B$ with $[G:B]\leq 4$.

\begin{proof}[Proof \Cref{necessary props 1}]
 As noted in \eqref{G is A semidir Q} one has a decomposition $G \cong A \rtimes_{\varphi} Q$ with $Q$ a $2$-group and $A$ an abelian $2'$-group. Note that $Q \cap F(G) = \Ker(\varphi)$ and thus in particular $F(G) = A \times \Ker(\varphi)$. By \Cref{broad form finite grps}~(3), it then follows that ${Q}/{\Ker(\varphi)}$ is an elementary abelian $2$-group, since ${Q}/{\Ker(\varphi)} \cong G/F(G)$, and $F(G)$ contains the group $B$ given by \Cref{broad form finite grps}~(3) since it is the maximal normal nilpotent subgroup of $G$. In particular $\exp (G/F(G)) \mid 2$ and $\exp (Q/(F(G)\cap Q))=2$.
 
 From \Cref{broad form finite grps}, it follows that the prime divisors of $\abs{G}$ are given by $\{2,3\}$. Hence $A$ is a $3$-group. To see (1), it remains to prove that $\exp(A) \mid 3$. We treat the cases $G$ nilpotent and $G$ non-nilpotent separately.

Suppose first that $G$ is nilpotent. Then since $G \cong A \times Q$, we have $\Q G \cong \Q A \otimes_{\Q} \Q Q$. Moreover, since $A$ is abelian, \[\Q A \cong \bigoplus_{d \: \mid \: \abs{A}}n_d \Q(\zeta_d),\] with $n_d$ the number of cyclic subgroups of $A$ of order $d$. Since $\Q G$ has a matrix component $M_2(D)$ by assumption, it follows that this component appears in the decomposition of $\Q Q$. Hence $\Q G$ has a component $\Ma_2(D \ot_{\Q} \Q(\zeta_d))$, for some division ring $D$ and some $d \: \mid \: \abs{A}$. Note that $\Q(\zeta_d)$ is contained in the centre of $D \ot_{\Q} \Q(\zeta_d)$. However since $d$ is odd, we see that the latter can only be an exceptional matrix component if $d = 3$ as $[\Q(\zeta_d): \Q] = \varphi(d)> 2$ otherwise. Hence $A$ is an elementary abelian $3$-group. 

Suppose now that $G$ is non-nilpotent. 
We consider $\overline{G} := A \rtimes (Q/\Ker(\varphi))$, which has property $\Mexc$ as well by \Cref{Prop components preserved for quotient}. Moreover, the action of ${Q}/{\Ker(\varphi)}$ is faithful by construction. As ${Q}/{\Ker(\varphi)}$ is elementary abelian it follows from \Cref{lower character degree classif} that $[\overline{G}: F(\overline{G})] \mid 4$. Additionally, by construction $F(\overline{G})=A$. Hence the statement will follow if we prove it for groups $A \rtimes C_2$ and $A \rtimes (C_2 \times C_2)$ with faithful action.

Let $C$ be $C_2$ or $C_2 \times C_2$. In both cases it follows from \cite[Lemma 1.2]{IsaacsPassman}  that $A \cong \left(A \cap \mc{Z}(A \rtimes C)\right) \times P$ for some characteristic subgroup $P\leqslant A\rtimes C$. Additionally, $P$ is non-trivial, since otherwise $A \leqslant \mc{Z}(A \rtimes C)$, a contradiction with the non-nilpotency assumption. Now, 
\begin{equation*}
P \rtimes C \cong G/\left(\mc{Z}(A \rtimes C) \cap A \right),
\end{equation*}
and in particular $P \rtimes C$ has $\Mexc$ by \Cref{Prop components preserved for quotient}. Let $x \in P$ have maximal order, say $o(x) = 3^m$. Then there is a subgroup $K \leqslant P$ such that $P \cong \langle x \rangle \times K$ (since $P \leqslant A$ is abelian).

There exists a SSP $(H,K)$ for $K$ as above. Indeed, the set 
\[S := \{E \leqslant P \rtimes C_2 \mid P \leqslant E \text{ and } E' \leqslant K \leqslant E\}\]
is non-empty since it contains $P$, and it is finite by construction. Then $(H,K)$ is given by a maximal element $H \in S$ by \Cref{SSP in metabelian}. From \cite[Lemma 3.4]{JM}, we have 
\begin{equation*}
    \dim_{\Q} \left( \Q[P \rtimes C]e(P \rtimes C, H, K) \right)  =  [P \rtimes C : H]\cdot  [H : N_H(K)] \cdot \phi([H:K]).
\end{equation*}
Moreover, since $P \rtimes C$ has $\Mexc$ by construction, the above $\Q$-dimension is an element of $\{4,8,16\}$ by \Cref{broad form finite grps}. Since $[P : K] = \abs{\langle x \rangle} = 3^m$, we have that 
\[\phi([P : K]) = 3^{m-1}(3-1).\]
Additionally, $\phi([P : K])$ divides $\phi([H : K])$ and hence we also have \[ 2\cdot 3^{m-1} \mid \dim_{\Q}\left( \Q[P \rtimes C]e(P \rtimes C, H, K) \right) \in \{4,8,16\}.\]
It follows that $m = 1$. In particular, by definition of $m$, $A$ is an elementary abelian $3$-group. Statement (1) follows.

\vspace{0,2cm}
We now show $(2)$. First note that $[Q : \Ker(\varphi)] = [G:F(G)]$ follows immediately since $A \leq F(G)$. Note also that $\Ker(\varphi)$ is a normal subgroup of $G$.
We now have that 
\[G/\Ker(\varphi) \cong A \rtimes Q/\Ker(\varphi) \cong C_3^n \rtimes C_2^k,\]
for some $n,k\in \N$ and where the action is faithful.

Since $A = C_3^n$ its automorphism group $\Aut(C_3^n)$ is isomorphic to $\GL_n(\mathbb{F}_3) = \GL(n,3)$. Furthermore, a faithful action of $C_2^k$ on $A$ is equivalent to an injection $ C_2^k  \hookrightarrow \GL_n(\mathbb{F}_3)$.

Every involution in $\GL(n,3)$ is diagonalisable with eigenvalues $\pm 1$, and any commuting family of such involutions is simultaneously diagonalisable. Thus any elementary abelian $2$-subgroup of $\GL(n,3)$ is conjugate to a subgroup of the diagonal matrices
\[
\{\diag(\varepsilon_1,\dots,\varepsilon_n) \mid \varepsilon_i= \pm 1 \}
\;\cong\; C_2^n.
\]
Hence the largest possible value of $k$ is $n$. Furthermore, for each possible value of $k \in \{0,\ldots,n\}$, one can choose the generators of $A$ and $Q/\Ker(\varphi)$ in such a way that the action is given by:
\[
\varphi\colon C_2^k  \longrightarrow \GL_n(\mathbb{F}_3),
\quad
y_i \mapsto \diag(1,\dots,1,\underset{i\text{-th}}{-1},1,\dots,1),
\]
for $ i\in\{1,\dots,k\}$, i.e. $y_i$ acts by inversion on the $i$-th factor of $A$ and trivially on the others. 

Now as $G$ has $\Mexc$, one has that $\cd (G) \subseteq \{1,2,4 \}$ by \Cref{broad form finite grps}. In particular $\cd(G/\Ker(\varphi)) \subseteq \{ 1,2,4 \}$ as well. Moreover, from the above description it is clear that $F(G/\Ker(\varphi)) \cong A$. Thus \Cref{lower character degree classif} yields that 
\[2^k= [G/\Ker(\varphi) : F(G/\Ker(\varphi))] = \max(\cd(G)).\]
Therefore $k \leq 2$. Since $2^k = [G:F(G)]$ as well, it suffices to discard the possibility $k=2$ to obtain statement (2).
By the above equality, that possibility only occurs if $4 \in \cd(G/\Ker(\varphi))$. This implies that $\Ma_2(\qa{-a}{-b}{\Q}) \in \mathcal{C}(G/\Ker(\varphi))$, since $G$ has $\Mexc$. On the other hand, if $n \gneq k=2$, then it follows from the description of the action $\varphi$ that there are some elements in $A$ which are centralised by $Q/\Ker(\varphi)$. In particular, $A \cap \ZZ(G/\Ker(\varphi)$ is non-trivial. However, this would imply that there exists a simple component $\Q(\zeta_3) \otimes \Ma_2(\qa{-a}{-b}{\Q})$ of $G$, which contradicts the fact that $G$ has $\Mexc$. Hence all that remains to understand is the case that $k = n =2. $ Note that in this case, \[G/\Ker(\varphi) \cong C_3^2 \rtimes C_2^2 \cong D_6^2.\] As $\Ma_2(\Q) \in \mc{C}(\Q [D_6])$, we have that $\Ma_2(\Q) \otimes_{\Q} \Ma_2(\Q) \cong \Ma_4(\Q) \in\mc{C}(\Q[D_6 \times D_6)]$. Thus $G/\Ker(\varphi)$ does not have $\Mexc$, a contradiction. It follows that $k = 1$, and in particular, $[G : F(G)] = 2$, as desired.
\end{proof}

\Cref{necessary props 1} yiels a decomposition of $A$ and $G'$ which will be instrumental in the study of non-nilpotent groups.

\begin{proposition}\label{prop:fixed-inverted}
Let $G$ be a finite group such that $\Q G$ has $\Mexc$ and at least one matrix component. Then:
\begin{enumerate}
\item\label{decomposition A} $A=A_{-1}\times (A\cap\ZZ(G))$ with
\[A_{-1} := \{ a \in A \mid a^{y}= a^{-1} \text{ for all } y \in G \setminus F(G) \}.\]
\item $A\cap\ZZ(G) = \Cen_A(Q)$ and $A_{-1} = [A,Q],$
\item\label{Comm sbgrp G} $G'=A_{-1}\times Q'$,
\item If $G$ is non-nilpotent, then
$$\Fit(G) \text{ is abelian } \Leftrightarrow \cd(G)\subseteq\{1,2\}.$$
\end{enumerate}
\end{proposition}
Recall that by \Cref{broad form finite grps}, one has that $\cd(G) \subseteq \{1,2,4 \}$ if $\Q G$ has $\Mexc$. Hence \Cref{prop:fixed-inverted} also implies that $\Fit(G)$ is abelian if and only $4 \notin \cd(G)$. Also note that, since $[G:\Fit(G)]=2$, one has that $a\in A_{-1}$ if and only if $a^y=a^{-1}$ for a single element $y\in G\setminus\Fit(G)$.
\begin{proof}
The first part follows from the fact $[G:F(G)]=2$ obtained in \Cref{necessary props 1}.  Indeed, remark that for a given $a \in A$ and $y \notin \Cen_G(\langle a\rangle)$, one has that $a^y = a^{-1}$ or $a^y \notin \langle a\rangle$, since $o(a)=3$ and $y^2 \in F(G)$. Now the elements in $A$ for which $\langle a \rangle$ is not normalised by $y$ can be obtained from an element in $A_{-1}$ and $A \cap \mc{Z}(G)$. Indeed, if $a^y \notin \langle a\rangle$ then $a^{-1}a^{y} \in A_{-1}$ and $aa^{y}\in A \cap\mc{Z}(G)$. Furthermore, $a = (a^{-1}a^{y}) \cdot (aa^{y})^{-1},$ which proves \eqref{decomposition A}. \smallskip

Next, the decomposition \eqref{Comm sbgrp G} follows from \eqref{decomposition A}.  Indeed, since $G \cong A \rtimes Q$, it follows that $G' = \langle [A,Q], Q' \rangle$. Hence to see \eqref{Comm sbgrp G}, it suffices to show that $[A,Q] = A_{-1}$, which is the statement. If $a \in A_{-1}$ and $y \in Q \setminus (Q \cap F(G))$, then $[a,y] = a^2 = a^{-1}$, and hence $a \in [A,Q]$. Conversely, the only non-trivial generators of $[A,Q]$ are those of the form $[a,y]$ with $a \in A \setminus (A\cap \mc{Z}(G))$ and $y \in Q \setminus (Q \cap F(G))$. But for any such $a$ and $y$, one has that $[a,y]^y = [a,y]^{-1}$, i.e. $[a,y] \in A_1$. Hence $[A,Q] = A_{-1}$.\medskip

Finally we consider the criterions for $F(G)$ to be abelian. Since $G$ is metabelian by \Cref{broad form finite grps} we can use \cite[Lemma 12.25]{Isaacs94} which says that 
\begin{equation}\label{max deg fitting}
\max\{ \chi(1) \mid \chi \in \Irr_{\C}(F(G)) \} \leq \frac{\max\{ \psi(1) \mid \psi \in \Irr_{\C}(G) \}}{[G:F(G)]}.
\end{equation}
If $\cd (G)= \{1,2 \}$ the right hand side of \eqref{max deg fitting} equals $1$ by \Cref{necessary props 1} and hence $F(G)$ is abelian. Conversely if $F(G)$ is abelian, then $\cd (G)= \{1,2 \}$ as $[G:F(G)]=2$ by \Cref{necessary props 1}.
\end{proof}

\subsection{On $2$-groups having \texorpdfstring{$\Mexc$}{(Mexc)} group}\label{subsec:sylow2}
Suppose that $\Q G$ has $\Mexc$ and hence $G=A\rtimes Q$ by \Cref{necessary props 1}. Then $\Q Q$ is a quotient of $\Q G$ and therefore $\Q Q$ also has $\Mexc$ by \Cref{Prop components preserved for quotient}. In this subsection we will collect properties of $2$-groups with $\Mexc$ which are needed in the sequel: (i) a bound on the nilpotency class and the exponent of $Q$ (\Cref{necessary props 2}), (ii) a description in terms of irreducible images (\Cref{thm:nilp-local-global}) and (iii) properties of the second centre when $Q$ has class three (\Cref{thm:nilp-class3-structure}).

\subsubsection{Nilpotency class and exponent}

The following result will motivate later on to classify the finite groups having $\Mexc$ depending on nilpotency class of a Sylow $2$-subgroup. 

\begin{lemma}\label{necessary props 2}
Let $G$ be a finite group and $Q$ a Sylow $2$-subgroup. Suppose that $\Q G$ has $\Mexc$ and has a matrix component. Then the following holds:
\begin{enumerate}
    \item The nilpotency class of $Q$ is at most $3$.
    \item \label{no exc if Q nonab} If $Q$ is non-abelian, then $\Q G$ has no exceptional division components. In particular 
    $$\exc(G)\mid 24 \text{ and } \exc(Q)\mid 8.$$
    \item If $Q$ is class $2$, then $\exp( Q/\mc{Z}(Q)) =2$.
\end{enumerate}
\end{lemma}
\begin{proof}
Note that $G$ maps onto $Ge$. If we choose a Sylow $2$-subgroup $Q$ of $G$, then $Qe$ will also be a Sylow $2$-subgroup of $Qe$. Therefore $cl(Qe) \leq cl(Q)$. Furthermore $Q \leq \prod_e Qe$. This yields that 
\begin{equation}\label{class Q}
cl(Q) = \max\{ cl(Qe) \mid e \in \PCI(\Q G) \}.
\end{equation}
Now verifying the possible $\Mexc$ groups classified in  the proof of \Cref{broad form finite grps} we get (1). \vspace{0,2cm}

For (2): The proof of \Cref{broad form finite grps} shows that the only way that $\Q G$ can have exceptional division components is if $G$ maps onto on either (I) a group $C_3 \rtimes C_{2^m}$ with $m \geq 4$ and the action by inversion or (II) a group of the form $C_m \times Q_8$ in the family b.iii) with $m$ odd. In \Cref{One degenerate case doesn't happen} it was obtained that case (II) can't occur as we assumed that $\Q G$ has a matrix component. Suppose that it would map on a group of the form $C_3 \rtimes C_{2^m}$, then there is an element $x \in Q$ with $v_2(o(x)) \geq m$.

Next, if $\Q G$ has $\Mexc$, then also the quotient $\Q Q. $ Since $Q$ is a non-abelian nilpotent $2$-group, it has no exceptional division components. Indeed, it can not map onto one of the non-nilpotent groups of type (I). Further it can only map on $Q_8$ in type (II) whose rational group algebra however has no exceptional division components.  Now \Cref{broad form finite grps} yields that $\exp(Q) \mid 24$ and hence a divisor of $8$, being a $2$-group. This shows that $m \leq v_2(o(x)) \leq 3$ and hence $G$ does not map onto a group of type (I) either, yielding (2).\smallskip

For (3) note that $\ZZ(Q) = Q \cap \prod_{e\in \PCI(\Q G)} \ZZ(Qe)$ and hence $\exp(Q/\ZZ(Q)) = \exp(Q /Q\cap \prod_e \ZZ(Qe)) \mid \exp (\prod_e Qe/\ZZ(Qe))$. Furthermore by \eqref{class Q} we have that $cl(Qe) \leq 2$ for all $e$. Also $\Q G$ has no exceptional components by (2), thus again looking at the groups obtained in the proof of \Cref{broad form finite grps}, combined with this extra restriction, yields that $\exp( Qe / \ZZ (Qe)) \mid 2$ for all $e$, yielding (3). 
\end{proof}

\subsubsection{Irreducible images}\label{subsec:local-global-crit}
Let $Q$ be a $2$-group and $e \in \PCI(\Q G)$ such that $\Q Qe$ is an exceptional matrix component. Then the quotient $Qe \cong Q/\ker(\varphi_e)$ of $Q$ is one of the groups in the table of Appendix \ref{appendix section}. The latter moreover records which quotients have $\Mexc$. We will now characterise $2$-groups having $\Mexc$ in terms of the quotients $Q/\ker(\varphi_e$ for $e \in \PCI( \Q G)$. For this consider the following list of $2$-groups:
\begin{multline}\label{eq:listL}
\mc{L}:=\{C_{2^n}\ (n\ge0),\ D_8,\ Q_8,\ Q_{16},\ M_{16}=\sgid{16}{6},\ SD_{16}=\sgid{16}{8},\ C_4\circ D_8=\sgid{16}{13},\\ \sgid{32}{8},\ \sgid{32}{11},\ \sgid{32}{44},\ \sgid{32}{50},\ \sgid{64}{137}\}.
\end{multline}
\begin{remark*}
Here $M_{16}=C_8\rtimes C_2$ with $C_2$ acting by $a\mapsto a^5$ and $\sgid{32}{50}\cong D_8\circ Q_8=2^{1+4}_-$ is the extraspecial group of minus type. The members of class at most two are the cyclic groups, $D_8$, $Q_8$, $M_{16}$, $P$ and $2^{1+4}_-$. The others have class three. The members with a character of degree four are $\sgid{32}{8}$, $\sgid{32}{44}$, $2^{1+4}_-$ and $\sgid{64}{137}$.
\end{remark*}

\begin{proposition}\label{thm:nilp-local-global}
Let $Q$ be a non-abelian $2$-group. Then $\Q Q$ has $\Mexc$ if and only if $Qe\in\mc{L}$ for every $e\in\PCI(\Q Q)$.
\end{proposition}
\begin{proof}
Suppose that $\Q Q$ has $\Mexc$ and let $e\in\PCI(\Q Q)$. The groups $Ge$ for $G$ with $\Mexc$ were described in the proof of \Cref{broad form finite grps}, and we only have to select the $2$-groups among them. If $\Q Qe$ is a division algebra, then $Qe$ is abelian, hence cyclic, or $Qe$ is $Q_8$ or $Q_{16}$. If $\Q Qe$ is not a division algebra, then $Qe$ is one of the groups in the table of \Cref{appendix section}, and $\Q[Qe]$ has $\Mexc$ by \Cref{Prop components preserved for quotient}. The $2$-groups in that table for which $\Mexc$ holds are exactly the non-cyclic members of \eqref{eq:listL} different from $Q_8$ and $Q_{16}$.

Conversely, suppose that $Qe\in\mc{L}$ for every $e\in\PCI(\Q Q)$. Each member of $\mc{L}$ has $\Mexc$, and $\Q Qe$ is a simple component of $\Q[Qe]$ by \Cref{Prop components preserved for quotient}. Hence every simple component of $\Q Q$ is a division algebra or an exceptional matrix algebra.
\end{proof}

\subsubsection{The case of nilpotency class three}\label{subsec:nilp-class3-structure}
When $\cl(Q)=3$, then \Cref{thm:nilp-local-global} can be used to obtain more information on the layers of the upper central series of $Q$. Since $\bigcap_{e\in\PCI(\Q Q)}\ker(\varphi_e)=1$, the group $Q$ embeds into $\prod_eQe$ as a subdirect product, i.e. with surjective projections. 

Denote by $\ZZ_2$ the second centre and $\gamma_3$ the third term of the lower central series.

\begin{theorem}\label{thm:nilp-class3-structure}
Let $Q$ be a $2$-group of class $3$ such that $\Q Q$ has $\Mexc$. Then:
\begin{enumerate}
\item $\gamma_3(Q)$, $\ZZ_2(Q)/\ZZ(Q)$ and $Q/\ZZ_2(Q)$ are elementary abelian. 
\item $\ZZ_2(Q)$ is abelian of exponent dividing $4$.
\end{enumerate}
In particular $Q^2Q'\le\ZZ_2(Q)$.
\end{theorem}

The proof will use that for a subdirect product $S\le G_1\times G_2$ one has
\begin{equation}\label{eq:subdirect-centres}
\ZZ(S)=S\cap(\ZZ(G_1)\times\ZZ(G_2)),\qquad \ZZ_2(S)\le\ZZ_2(G_1)\times\ZZ_2(G_2),\qquad \gamma_3(S)\le\gamma_3(G_1)\times\gamma_3(G_2),
\end{equation}
Indeed, as the projections $\pi_i$ are surjective, $\pi_i(\ZZ(S))\le\ZZ(G_i)$ and $\pi_i([y,S])=[\pi_i(y),G_i]$ for $y\in S$.

\begin{remark}\label{rem:nilp-class3-maxab}
Let $Q$ be as in \Cref{thm:nilp-class3-structure}. Then $\ZZ_2(Q)$ is an abelian normal subgroup of $Q$ containing $Q^2Q'$. Hence every subgroup $A$ of $Q$ containing $\ZZ_2(Q)$ is normal with $Q/A$ elementary abelian. If $A$ is a maximal abelian subgroup of $Q$ containing $\ZZ_2(Q)$, then $\Cen_Q(A)=A$, and $[Q:A]\ge4$ if $4\in\cd(Q)$ by Ito's theorem \cite[Theorem 6.15]{Isaacs94}.
\end{remark}
\begin{proof}[Proof of \Cref{thm:nilp-class3-structure}]
We start by proving the first part. By \Cref{thm:nilp-local-global} every $R_{\chi}$ belongs to $\mc{L}$. For $R\in\mc{L}$ the group $\gamma_3(R)$ has order at most $2$ (it is contained in $\ZZ(R)\cap R'$, which has order $2$ by \Cref{appendix section}), $R/\ZZ(R)$ is elementary abelian if $\cl(R)\le2$ (\Cref{thm:nilp-local-global}), and for the members of class three $R/\ZZ(R)$ is one of $D_8$, $\sgid{16}{3}$, $C_2\times D_8$, $\sgid{32}{27}$ (see Appendix \ref{appendix section}), whose centres are elementary abelian. So $\ZZ_2(R)/\ZZ(R)=\ZZ(R/\ZZ(R))$ is elementary abelian for every $R\in\mc{L}$. By \eqref{eq:subdirect-centres}, $\gamma_3(Q)$ embeds into $\prod_\chi\gamma_3(R_\chi)$ and $\ZZ_2(Q)/\ZZ(Q)$ embeds into $\prod_\chi\ZZ_2(R_\chi)/\ZZ(R_\chi)$, which gives the first two assertions. Finally, by \Cref{Prop components preserved for quotient}, $Q/\ZZ(Q)$ is a $2$-group of class two with $\Mexc$, so $Q/\ZZ_2(Q)=(Q/\ZZ(Q))/\ZZ(Q/\ZZ(Q))$ is elementary abelian by \Cref{thm:nilp-local-global}. Thus $Q^2\le\ZZ_2(Q)$, and $Q'\le\ZZ_2(Q)$ because $\cl(Q)=3$.\medskip

Now we will prove the second assertion by induction on $|Q|$. If $\ZZ(Q)$ is cyclic, then $Q$ has a faithful irreducible character, so $Q$ itself is a member of $\mc{L}$ of class three, i.e.\ $Q_{16}$, $SD_{16}$, $\sgid{32}{11}$, $\sgid{32}{8}$, $\sgid{32}{44}$ or $\sgid{64}{137}$, and $\ZZ_2(Q)$ is abelian of exponent $4$: it is the cyclic group $\langle r^2\rangle$ for $Q_{16}=\langle r,s\rangle$ and $SD_{16}=\langle r,u\rangle$, and for the four other groups it is isomorphic to $C_2\times C_4$, as one checks on the matrix models of these groups given in Appendix \ref{appendix section}. So let $\ZZ(Q)$ be non-cyclic. Then $\Omega_1(\ZZ(Q))$ has at least three subgroups of order two, and $\gamma_3(Q)$ is one of its subgroups by (1). Hence we can choose two distinct subgroups $N_1,N_2\le\ZZ(Q)$ of order two with $\gamma_3(Q)\not\le N_i$. The quotients $G_i:=Q/N_i$ have $\Mexc$ (\Cref{Prop components preserved for quotient}), class three (as $\gamma_3(G_i)=\gamma_3(Q)N_i/N_i\ne1$) and smaller order, so $\ZZ_2(G_i)$ is abelian of exponent dividing $4$ by induction. As $N_1\cap N_2=1$, the group $Q$ embeds into $G_1\times G_2$ as a subdirect product, and \eqref{eq:subdirect-centres} gives $\ZZ_2(Q)\le\ZZ_2(G_1)\times\ZZ_2(G_2)$, so that $\ZZ_2(Q)$ is abelian of exponent dividing $4$.
\end{proof}

\subsection{Non-nilpotent groups: an abelian normal subgroup of index at most four}\label{subsec:fourterm}
In this subsection $G=A\rtimes Q$ is a finite \emph{non-nilpotent} group such that $\Q G$ has $\Mexc$ and the main result is the following.

\begin{theorem}\label{thm:abelian-B}
Let $G$ be a finite non-nilpotent group such that $\Q G$ has $\Mexc$. Then $\Fit(G)_2= \Cen_Q(A)$ has an abelian subgroup $L$ of index at most $2$ which is normal in $Q$. For every such $L$, the group $B:=A\times L$ is a maximal abelian subgroup of $G$ containing $G'$, and $[G:B]=[Q:L]\le4$.
\end{theorem}

\begin{remark}\label{rem:abelian-B-nilpotent}
For a $2$-group $Q$ with $\Mexc$ the analogue of \Cref{thm:abelian-B} (an abelian normal subgroup with quotient $C_2$ or $C_2^2$) holds when $\cd(Q)\subseteq\{1,2\}$, but not in general. Indeed, the group $\sgid{256}{28688}$ has class two, $\Mexc$ and a character of degree four, and its maximal abelian subgroups have index $8$ (see \cite{CJT3}).
\end{remark}

 The proof of \Cref{thm:abelian-B} requires \Cref{prop:nonlinear-C} below. The latter determines the non-commutative simple components of $\Q[\Fit(G)_2]$, using that $\Fit(G)$ inherits $\Mexc$ from $G$. For simplicity we use following notation in the sequel of the section:
 $$C := \Fit(G)_2 = \Cen_Q(A).$$

\begin{proposition}\label{prop:nonlinear-C}
Let $G$ be a finite non-nilpotent group such that $\Q G$ has $\Mexc$, and let $e\in\PCI(\Q C)$ be such that $\Q Ce$ is non-commutative. Then $Ce$ is isomorphic to $D_8$ or $Q_8$, and $e$ is central in $\Q G$. Consequently the following holds:
\begin{enumerate}[(i)]
\item $\cd(C)\subseteq\{1,2\}$ and $\cl(C)\leq 2$, 
\item $C'$ is elementary abelian and $C^2\le\ZZ(C)$, 
\item every non-commutative simple component of $\Q[\Fit(G)]$ is isomorphic to $\Ma_2(\Q)$, $\Ma_2(\Q(\sqrt{-3}))$ or $\HH$.
\end{enumerate}
\end{proposition}
\begin{proof}
Put $R:=Ce$ and $S:=\Q Ce$, so that $R$ is a non-abelian $2$-group spanning $S$. The centre $Z$ of $S$ is isomorphic to a subfield of the cyclotomic field $\Q(\zeta_{|C|})$, as it is generated by the values of an irreducible character of $C$. In particular $\sqrt{-3}\notin Z$, and $\Q(\sqrt{-3})\otimes_\Q S\cong Z(\sqrt{-3})\otimes_ZS$ is a simple algebra with centre $Z(\sqrt{-3})$, of the same degree as $S$.

\emph{Step 1: $\Q(\sqrt{-3})\otimes_\Q S$ is a simple component of $\Q[\Fit(G)]$.} Since $G$ is not nilpotent, $A\ne1$. By the theorem of Perlis--Walker theorem, $\Q A\cong\Q\oplus\Q(\sqrt{-3})^{s}$ with $s=(|A|-1)/2\ge1$. Hence
\[
\Q[\Fit(G)]=\Q[A\times C]\cong\Q A\otimes_\Q\Q C\cong\Q C\oplus(\Q(\sqrt{-3})\otimes_\Q\Q C)^{s},
\]
and $\Q(\sqrt{-3})\otimes_\Q S$ is a direct summand of $\Q[\Fit(G)]$. By \Cref{Mexc under subgroups}, $\Q[\Fit(G)]$ has $\Mexc$.

\emph{Step 2: $Ce$ is isomorphic to $D_8$ or $Q_8$.} Suppose first that $S$ is not a division algebra. Then neither is $\Q(\sqrt{-3})\otimes_\Q S$, so by Step 1 it is an exceptional matrix algebra. Its centre contains $\Q(\sqrt{-3})$, so by \Cref{def_exc_comps} it is $\Ma_2(\Q(\sqrt{-3}))$. Thus $Z=\Q$ and $S$ is a non-division algebra of degree $2$ over $\Q$, i.e.\ $S\cong\Ma_2(\Q)$. Now $D_8$ is the only $2$-group spanning $\Ma_2(\Q)$ (\Cref{appendix section}), so $R\cong D_8$. Suppose next that $S$ is a division algebra. By Amitsur's classification \cite{Amitsur}, the non-abelian $2$-group $R$ is a generalised quaternion group $\langle a,b\mid a^{2^{n-1}}=1,\ b^2=a^{2^{n-2}},\ a^b=a^{-1}\rangle$ with $n\ge3$. Then $S=\Q(\zeta)\oplus\Q(\zeta)b$ with $\zeta=\zeta_{2^{n-1}}$, $b\zeta=\zeta^{-1}b$ and $b^2=-1$, and $Z=\Q(\zeta+\zeta^{-1})$. The elements $i:=\zeta^{2^{n-3}}$ and $b$ satisfy $i^2=b^2=-1$ and $bi=-ib$, so $S\cong\qa{-1}{-1}{Z}$. In $\Q(\sqrt{-3})$ one has $-1=\omega^2+(\omega^2)^2$ for a primitive third root of unity $\omega$, so $\qa{-1}{-1}{K}\cong\Ma_2(K)$ for every field $K$ containing $\Q(\sqrt{-3})$. Hence $\Q(\sqrt{-3})\otimes_\Q S\cong\Ma_2(Z(\sqrt{-3}))$, which by Step 1 is an exceptional matrix algebra. As before this forces $Z=\Q$, i.e.\ $n=3$ and $R\cong Q_8$, with $S\cong\HH$.

In both cases we obtained that $\Q(\sqrt{-3})\otimes_\Q S\cong\Ma_2(\Q(\sqrt{-3}))$.

\emph{Step 3: $e$ is central in $\Q G$.} Since $e$ is central in $\Q C$ and $A$ centralises $C$, it suffices to show that $e^t=e$ for $t\in Q\setminus C$. Suppose that $e^t\ne e$. Then $e^t$ is a primitive central idempotent of $\Q C$ different from $e$, so $ee^t=0$. Choose a subgroup $A_0$ of index $3$ of $A$ containing $A\cap\ZZ(G)$, which exists since $A_{-1}\ne1$. It is normal in $G$, because $t$ inverts $A_{-1}$. Hence $f:=\widehat{A_0}-\widehat A$ is a central idempotent of $\Q G$, and $\Q Af\cong\Q(\sqrt{-3})$ by the Perlis--Walker theorem. Put $\varepsilon:=fe$. By Step 2,
\[
\Q[\Fit(G)]\varepsilon\cong\Q Af\otimes_\Q\Q Ce\cong\Q(\sqrt{-3})\otimes_\Q S\cong\Ma_2(\Q(\sqrt{-3})).
\]
Moreover $\varepsilon^t=fe^t$ is orthogonal to $\varepsilon$, so $E:=\varepsilon+\varepsilon^t$ is a central idempotent of $\Q G$. As $G=\Fit(G)\cup\Fit(G)t$ and $\varepsilon t\varepsilon=t\varepsilon^t\varepsilon=0$, we have $\varepsilon\Q G\varepsilon=\Q[\Fit(G)]\varepsilon$. Since $\varepsilon$ and $\varepsilon^t$ are orthogonal idempotents with sum $E$ which are conjugate by the unit $t$, we obtain
\[
\Q GE\cong\Ma_2(\varepsilon\Q G\varepsilon)\cong\Ma_4(\Q(\sqrt{-3})).
\]
This is a simple component of $\Q G$ which is neither a division algebra nor an exceptional matrix algebra, a contradiction. \smallskip

Finally, we explain how the remaining assertions follow from the above. Firstly, by Step 2 every simple component of $\Q C$ is a field, $\Ma_2(\Q)$ or $\HH$, so $\cd(C)\subseteq\{1,2\}$. Since $\bigcap_{e\in\PCI(\Q C)}\ker(\varphi_e)=1$, the group $C$ is a subdirect product of cyclic groups and of copies of $D_8$ and $Q_8$, and all these groups have class at most $2$, elementary abelian derived subgroup and central squares. Finally $\Q(\sqrt{-3})\otimes_\Q\Ma_2(\Q)\cong\Q(\sqrt{-3})\otimes_\Q\HH\cong\Ma_2(\Q(\sqrt{-3}))$, and $\Q(\sqrt{-3})\otimes_\Q K$ is commutative for a field $K$. Together with the decomposition in Step 1 this gives the last assertion.
\end{proof}

\begin{proof}[Proof of \Cref{thm:abelian-B}]
If $C$ is abelian, then $L=C$ has the required property. So for the first assertion we assume that $C$ is non-abelian. By \Cref{prop:nonlinear-C}, $\cd(C)=\{1,2\}$, $C$ has class two, $C'$ is elementary abelian and $C^2\le\ZZ(C)$.

\emph{Step 1: $C$ has an abelian subgroup of index $2$.} By a theorem of Isaacs \cite[Theorem 12.11]{Isaacs94}, as $\cd(C)=\{1,2\}$, either $C$ has an abelian subgroup of index $2$, or $|C:\ZZ(C)|=8$. Assume the latter. Then $C/\ZZ(C)\cong C_2^3$. Choose $a,b,d\in C$ whose images generate $C/\ZZ(C)$ and put
\[
e_1:=[a,b],\qquad e_2:=[a,d],\qquad e_3:=[b,d].
\]
As $C$ has class two, commutators are bilinear, so that $C'=\langle e_1,e_2,e_3\rangle$ and, for $c\in C$, the set $[c,C]$ of commutators $[c,g]$, $g\in C$, is a subgroup of $C'$. If $|C'|\le4$, there is a relation $e_1^{\alpha}e_2^{\beta}e_3^{\gamma}=1$ with $(\alpha,\beta,\gamma)\ne(0,0,0)$, and then $C$ has an abelian subgroup of index $2$, namely $\langle g,h\rangle\ZZ(C)$ with
\[
(g,h)=\begin{cases}(ba^{\beta},\,da^{\alpha})&\text{if }\gamma=1,\\(a,\,b^{\alpha}d)&\text{if }\gamma=0,\ \beta=1,\\(a,b)&\text{if }\gamma=\beta=0,\end{cases}
\]
since $[ba^\beta,da^\alpha]=e_3e_1^{\alpha}e_2^{\beta}$ and $[a,b^\alpha d]=e_1^\alpha e_2$. So assume that $|C'|=8$. Then for $c\in C\setminus\ZZ(C)$ the subgroup $[c,C]$ has order $4$. Indeed, as $[cz,C]=[c,C]$ for $z\in\ZZ(C)$, it suffices to list
\begin{equation}\label{eq:abelian-B-commutators}
\begin{gathered}
[a,C]=\langle e_1,e_2\rangle,\quad [b,C]=\langle e_1,e_3\rangle,\quad [d,C]=\langle e_2,e_3\rangle,\quad [ab,C]=\langle e_1,e_2e_3\rangle,\\
[ad,C]=\langle e_2,e_1e_3\rangle,\quad [bd,C]=\langle e_3,e_1e_2\rangle,\quad [abd,C]=\langle e_1e_2,e_1e_3\rangle .
\end{gathered}
\end{equation}
We claim that
\begin{equation}\label{eq:abelian-B-squares}
c^2\in[c,C]\qquad\text{for every }c\in C\setminus\ZZ(C).
\end{equation}
Fix $c\in C\setminus\ZZ(C)$, let $\lambda\in\Irr(C')$ be the character with kernel $[c,C]$ and put $Y:=\ZZ(C)\langle c\rangle$, an abelian subgroup containing $C'$. Let $\nu\in\Irr(Y)$ be an extension of $\lambda$ and $\psi\in\Irr(C)$ a constituent of the induced character $\nu^C$. Since $C'$ is central, $\psi|_{C'}=\psi(1)\lambda$. Hence $\psi$ is non-linear and $[Y,C]=[c,C]=\ker\lambda\le\ker\psi$, so the image of $Y$ in $C/\ker\psi$ is central. By \Cref{prop:nonlinear-C}, $C/\ker\psi$ is $D_8$ or $Q_8$, whose centre has order two. As $Y$ acts by scalars, $\psi|_Y=2\nu$, so $\nu(Y)\cong Y/(Y\cap\ker\psi)$ has order at most two and $\nu(c^2)=1$. Now $c^2\in\ZZ(C)$, and every extension $\mu\in\Irr(\ZZ(C))$ of $\lambda$ extends further to $Y$. Thus $\mu(c^2)=1$ for all these $\mu$. They form a coset of $\Irr(\ZZ(C)/C')$, so their common kernel is $\ker\lambda$, and $c^2\in\ker\lambda=[c,C]$. This proves \eqref{eq:abelian-B-squares}.

Now we derive a contradiction from \eqref{eq:abelian-B-commutators} and \eqref{eq:abelian-B-squares}. Write $a^2=e_1^{\alpha_1}e_2^{\alpha_2}$, $b^2=e_1^{\beta_1}e_3^{\beta_3}$ and $d^2=e_2^{\delta_2}e_3^{\delta_3}$. Since $(gh)^2=g^2h^2[h,g]$ in a group of class two with elementary abelian derived subgroup, we have
\begin{align*}
(ab)^2&=a^2b^2e_1=e_1^{\alpha_1+\beta_1+1}e_2^{\alpha_2}e_3^{\beta_3}\in\langle e_1,e_2e_3\rangle, &&\text{so }\alpha_2=\beta_3,\\
(ad)^2&=a^2d^2e_2=e_1^{\alpha_1}e_2^{\alpha_2+\delta_2+1}e_3^{\delta_3}\in\langle e_2,e_1e_3\rangle, &&\text{so }\alpha_1=\delta_3,\\
(bd)^2&=b^2d^2e_3=e_1^{\beta_1}e_2^{\delta_2}e_3^{\beta_3+\delta_3+1}\in\langle e_3,e_1e_2\rangle, &&\text{so }\beta_1=\delta_2.
\end{align*}
Finally $(abd)^2=a^2b^2d^2e_1e_2e_3$ lies in $[abd,C]=\langle e_1e_2,e_1e_3\rangle$, which consists of the elements $e_1^{x_1}e_2^{x_2}e_3^{x_3}$ with $x_1+x_2+x_3$ even, but the exponent sum of $a^2b^2d^2e_1e_2e_3$ is $\alpha_1+\alpha_2+\beta_1+\beta_3+\delta_2+\delta_3+3$, which is odd by the three relations above. This contradiction shows that $C$ has an abelian subgroup of index $2$.

\emph{Step 2: an abelian subgroup of index $2$ which is normal in $Q$.} Let $X$ be an abelian subgroup of index $2$ of $C$. If it is the only one, it is characteristic in $C$, hence normal in $Q$. Otherwise let $X_1\ne X$ be another one. Then $X\cap X_1$ is central in $C=XX_1$ and has index $4$, so $X\cap X_1=\ZZ(C)$ and $C/\ZZ(C)\cong C_2^2$. Each of the three subgroups of index $2$ of $C$ containing $\ZZ(C)$ is abelian, and every abelian subgroup of index $2$ contains $\ZZ(C)$. The group $Q$ permutes these three subgroups by conjugation, and $C$ fixes each of them. As $[Q:C]=2$, an element $t\in Q\setminus C$ induces a permutation of order at most two of a set with three elements, which has a fixed point $L$. Then $L$ is normal in $Q=C\langle t\rangle$.\medskip

To finish we prove the remaining assertion. Namely, let $L$ be any abelian subgroup of $C$ of index at most $2$ which is normal in $Q$, and $B=A\times L$. It is abelian and normal in $G$, and $[G:B]=[Q:L]$. Since $Q/L$ has order at most $4$ it is abelian, so $Q'\le L$ and $G'=A_{-1}\times Q'\le B$ by \Cref{prop:fixed-inverted}. If $B\le B_1$ with $B_1$ abelian normal, then $B_1=A\times(B_1\cap Q)$ with $L\le B_1\cap Q\le C$. If $C$ is abelian then $L=C$, and if not then $B_1\cap Q\ne C$. In both cases $B_1=B$.
\end{proof}

All finite groups such that $\SL_1(\Q Ge)$ is a discrete subgroup of $\SL_2(\C)$ have been classified in \cite{FreebyFree}. In loc.cit. presentations for such groups were even given. Furthermore, they showed that the aforementioned property is equivalent to saying that all simple components of $\Q G$ are either fields, totally definite quaternion algebras or of the form $\Ma_2(\Q(\sqrt{-d}) )$ with $d \geq 0$. The following result implies that groups with $\Mexc$ are index $2$ overgroups of the groups classified in \cite{FreebyFree}. 

\begin{corollary}\label{cor:index2}
Let $G$ be a finite non-nilpotent group such that $\Q G$ has $\Mexc$. Then $[G:\Fit(G)]=2$ and  $\Q[\Fit(G)]$ has no exceptional division components and all its non-division components are isomorphic to $\Ma_2(\Q)$ or $\Ma_2(\Q(\sqrt{-3}))$.
\end{corollary}
\begin{proof}
By \Cref{prop:nonlinear-C} the non-commutative simple components of $\Q[\Fit(G)]$ are $\Ma_2(\Q)$, $\Ma_2(\Q(\sqrt{-3}))$ and $\HH$, and $\HH$ is a totally definite quaternion algebra.
\end{proof}

\begin{example}
The converse of \Cref{cor:index2} is not true in general. For example consider the following extraspecial group of order $2^5$ (whose SmallGroup ID is [32,49]):
\begin{align*}
D_8 \circ D_8 := \langle a,b,c,d \mid & \: a^4 = b^2 = c^2=d^2=1, c^2 =a ^2, a^b = a^{-1}, c^d = a^2c, \\ & \: [a,c] = [a,d] = [b,c] = [b,d] = 1 \rangle.
\end{align*}
It can be verified that $\mc{C}(\Q[D_8 \circ D_8])= \{\Q, \Ma_4(\Q) \}$ (e.g. via the well-known description of its complex irreducible representations). On the other, e.g. to be proven via a manual check (via GAP), all the $2$-groups until order $16$, except $D_{16}$, have only matrix components of the form $\Ma_2(\Q(\sqrt{-d}) )$ with $d \geq 0$.
\end{example}

\section{Strong Shoda pairs and property \texorpdfstring{$\Mexc$}{Mexc}}\label{SSP and prop sectie}
By \Cref{necessary props 2}, $\Mexc$ groups such that $\cd(G) = \{1,2\}$ and non-abelian Sylow $2$-subgroup have been classified in \cite{FreebyFree}. This includes the case that $G$ has an abelain subgroup of index $2$. Therefore, by \Cref{thm:abelian-B}, the focus will be on the case that $G$ has an abelian normal subgroup $B$ such that $G/B \cong C_2 \times C_2$. The aim of this section is firstly to somehow describe the strong shoda pairs under that assumption. Secondly, in \Cref{thm:uniform} we characterise $\Mexc$ under the extra assumption that $\pi(G) \subseteq \{2,3 \}$ . In \Cref{sec:nonnilpotent} we will further specialize the latter result in case that $G$ is non-nilpotent.

\subsection{Strong Shoda pairs with respect to a fixed abelian normal subgroup}\label{sectie SSP Mexc}
Let $G$ be a finite metabelian group and $B$ a maximal abelian subgroup of $G$ containing $G'$. By \Cref{SSP in metabelian} the primitive central idempotents of $\Q G$ are the elements $e(G,H,K)$, where $(H,K)$ runs over the pairs with
\begin{equation}\label{eq:SSP-shape}
\begin{gathered}
B\le H\le G,\qquad H/K \text{ cyclic},\\
H \text{ maximal among the subgroups } X \text{ with } B\le X\le G \text{ and } X'\le K\le X.
\end{gathered}
\end{equation}
We call them the strong Shoda pairs of $G$ \emph{with respect to $B$}. We call a subgroup $N$ of a group $X$ \emph{cocyclic} if $X/N$ is cyclic. For abelian $X$ the cocyclic subgroups are exactly the kernels of the linear characters of $X$. When $G/B$ is elementary abelian of order at most four, the pairs \eqref{eq:SSP-shape} can be listed explicitly.

\begin{proposition}\label{prop:SSP-B}
Let $G$ be a finite group and $B$ a maximal abelian subgroup of $G$ containing $G'$ such that $G/B$ is elementary abelian of order at most $4$. For $y\in G\setminus B$ write $H_y:=\langle B,y\rangle$. Then the non-commutative simple components of $\Q G$ are the algebras $\Q Ge(G,H,K)$, where $(H,K)$ is one of the following strong Shoda pairs:
\begin{enumerate}[(a)]
\item $(B,K)$ with $K$ cocyclic in $B$ such that no element of $G\setminus B$ centralises $B/K$. The component has degree $[G:B]$ over its centre;
\item $(H_y,K)$ with $y\in G\setminus B$ and $K$ cocyclic in $H_y$ such that $G'\not\le K$. The component has degree $2$ over its centre.
\end{enumerate}
\end{proposition}
\begin{proof}
By \Cref{SSP in metabelian} and \Cref{rem:H=G},  $\Q Ge(G,H,K)$ is not a field exactly for the pairs $(H,K)$ as in \eqref{eq:SSP-shape} with $H\neq G$, and the degree is $[G:H]$ by \Cref{form of SSP idempotent result}. The subgroups $H$ with $B\le H<G$ are $B$ and, if $|G/B|=4$, the three subgroups $H_y$. Note that $H_y'=[B,y]$, as $B$ is abelian and $H_y/B$ is cyclic.

Let $H=B$ and let $K$ be cocyclic in $B$. Then $B$ is maximal among the $X$ with $X'\le K$ if and only if $H_y'=[B,y]\not\le K$ for all $y\in G\setminus B$, i.e.\ if and only if no $y\in G\setminus B$ centralises $B/K$. Indeed, if $G'\le K$ then also $[B,y]\le K$. This gives (a).

Let $|G/B|=4$, $H=H_y$ and let $K$ be cocyclic in $H_y$. Then $H_y'\le K$, and the only subgroup properly containing $H_y$ is $G$. So $H_y$ is maximal among the $X$ with $X'\le K$ if and only if $G'\not\le K$. This gives (b). If $|G/B|=2$, then $H_y=G$ and there is no pair as in (b).
\end{proof}

\begin{remark*}
In terms of characters, the pairs $(B,K)$ in (a) correspond to the linear characters of $B$ with trivial stabiliser in $G/B$, and the pairs $(\langle B, y \rangle,K)$ in (b) to the linear characters of $H_y$ which are non-trivial on $G'$. 

Recall also that, by \Cref{rem:H=G}, two pairs with the same first subgroup give the same component if and only if their second subgroups are conjugate in $G$. Note that pairs $(B,K)$ and $(\langle B, y\rangle,K')$ never give the same component as their degrees are different.
\end{remark*}

\subsection{The strong Shoda pair criterion when abelian index four}\label{subsec:uniform}
For simplicity, in the sequel of this section we will use following notation: 
$$H_y:=\langle B,y\rangle \text{ for } y\in G\setminus B.$$ 
The conditions of \Cref{thm:uniform} say how $G/B$ acts on the cyclic group $H/K$ and where the squares of the elements of $G\setminus B$ lie. Note that for $z\in G\setminus H$ one has $z^2\in B\le H$, and that $z^2K$ is fixed by $z$, so that $z^2\in\langle h^{m/2}K\rangle$ whenever $z$ inverts $H/K=\langle hK\rangle$ of order $m$.

\begin{theorem}\label{thm:uniform}
Let $G$ be a finite group with $\pi(G)\subseteq\{2,3\}$ and let $B$ be a maximal abelian subgroup of $G$ with $G/B\cong C_2^2$. Then $\Q G$ has $\Mexc$ if and only if the following two conditions hold.
\begin{enumerate}[(U1)]
\item For every cocyclic subgroup $K$ of $B$ such that no element of $G\setminus B$ centralises $B/K$, one of the following holds:
\begin{enumerate}[(a)]
\item $[B:K]\in\{4,6\}$, $K$ is normalised by some $y\in G\setminus B$, and $y^2\notin K$;
\item $[B:K]\in\{8,12\}$, $K\trianglelefteq G$, and $z^2\notin K$ for the elements $z\in G\setminus B$ inverting $B/K$.
\end{enumerate}
\item For every $y\in G\setminus B$ and every cocyclic subgroup $K$ of $H_y$ with $G'\not\le K$, one of the following holds, where $z\in G\setminus H_y$ and $m=[H_y:K]$:
\begin{enumerate}[(a)]
\item $K$ is not normal in $G$ and $m\in\{2,3,4,6\}$;
\item $K\trianglelefteq G$ and $m\in\{3,4,6\}$;
\item $K\trianglelefteq G$, $m\in\{8,12\}$, and $z^2\notin K$ if $z$ inverts $H_y/K$;
\item $K\trianglelefteq G$, $m=3\cdot2^\ell$ with $\ell\ge3$, $z$ inverts the Sylow $3$-subgroup of $H_y/K$ and centralises its Sylow $2$-subgroup, and $z^2K$ generates the Sylow $2$-subgroup of $H_y/K$.
\end{enumerate}
\end{enumerate}
\end{theorem}

We will also need the description of the simple components corresponding to the cases obtained in \Cref{thm:uniform}. For readability, we record them seperately in the following result.

\begin{lemma}\label{prop:uniform-components}
Let $G$ and $B$ be as in \Cref{thm:uniform} and suppose that $\Q G$ has $\Mexc$. The simple components of $\Q G$ are $\Q[G/G']$ and the algebras $\Q Ge(G,H,K)$, $(H,K)$ running over the pairs in \textup{(U1)} and \textup{(U2)}. The component of such a pair is the following algebra, with $m=[H:K]$:
\begin{itemize}
\item \textup{(U1a)}: $\Ma_2(\HH)$ if $m=4$ and $\Ma_2(\HHt)$ if $m=6$;
\item \textup{(U1b)}: $\Ma_2(\HH)$ if $m=8$. If $m=12$, $\Ma_2(\HH)$ if the Sylow $2$-subgroups of $G/K$ are isomorphic to $Q_8\times C_2$ or $SD_{16}$, and $\Ma_2(\HHt)$ if they are isomorphic to $C_4\circ D_8$ or $Q_{16}$ \textup{(}\Cref{conditions comp exc}\,(v)\textup{)};
\item \textup{(U2a)}: $\Ma_2(\Q(\zeta_m))$;
\item \textup{(U2b)}: $\Ma_2(\Q)$ if $z^2\in K$ \textup{(}always the case for $m=3$\textup{)}, and $\HH$, resp.\ $\HHt$, if $z^2\notin K$ and $m=4$, resp.\ $m=6$;
\item \textup{(U2c)}: for $m=8$: $\Ma_2(\Q(\sqrt{-2}))$, $\Ma_2(\Q(i))$, $\HHrtt$ according as $z$ acts on $H_y/K=\langle hK\rangle$ by $h\mapsto h^3$, $h^5$, $h^{-1}$. For $m=12$: $\Ma_2(\Q(i))$, $\Ma_2(\Q(\sqrt{-3}))$, $\HHrt$ according as $z$ acts by $h\mapsto h^5$, $h^7$, $h^{-1}$;
\item \textup{(U2d)}: $\qa{\zeta_{2^\ell}}{-3}{\Q(\zeta_{2^\ell})}$.
\end{itemize}
\end{lemma}
\begin{proof}[Proof of \Cref{thm:uniform} and \Cref{prop:uniform-components}]
By \Cref{prop:SSP-B}, the non-commutative simple components of $\Q G$ are the algebras $\Q Ge(G,H,K)$ for the pairs $(B,K)$ with $K$ as in (U1) and the pairs $(H_y,K)$ with $K$ as in (U2), (so every non-commutative component arises from such a pair). By \Cref{thm componenten QG mexc}, $\Q G$ has $\Mexc$ if and only if each of these components is a field, one of the division algebras $\HH$, $\HHt$, $\HHrt$, $\HHrtt$, $\qa{\zeta_{2^\ell}}{-3}{\Q(\zeta_{2^\ell})}$ ($\ell\ge3$), or one of the exceptional matrix algebras \eqref{eq:exceptional-list}. We apply \Cref{conditions comp exc,prop:centre} to each pair, with $N=\Nor_G(K)$, $m=[H:K]$ and $H/K=\langle hK\rangle$.

\emph{Pairs $(B,K)$: condition \textup{(U1)}.} Here $[G:B]=4$, so the component is admissible only if it is $\Ma_2(D_0)$ with $D_0\in\{\HH,\HHt\}$, which by \Cref{conditions comp exc} happens only in the cases (iv) and (v) there.

In case (iv), $[G:N]=[N:B]=2$ and $m\in\{4,6\}$, so $N=H_y$ for some $y$, which inverts $B/K$, and the condition $y^2\equiv h^{m/2}\pmod K$ means $y^2\notin K$, because $y^2K$ is fixed by $y$ and the fixed points of the inversion on $B/K$ form $\langle h^{m/2}K\rangle$. Conversely under these conditions the component is $\Ma_2(\HH)$, resp.\ $\Ma_2(\HHt)$. This is (U1a).

In case (v), $N=G$, $m\in\{8,12\}$, $G/B$ acts faithfully on $B/K$, hence as $\Aut(B/K)\cong C_2^2$, and $z^2\notin K$ for the elements $z$ inverting $B/K$. The component is $\Ma_2(D_0)$ with $D_0$ given by \Cref{conditions comp exc}\,(v) for $G/K$. This is (U1b).

All other pairs $(B,K)$ give non-admissible components: if $N=B$ the component is $\Ma_4(\Q(\zeta_m))$. If $[G:N]=2$ and $m\notin\{4,6\}$, or $N=G$ and $m\notin\{8,12\}$, the centre has degree $\phi(m)/[N:B]\ne1$ (\Cref{prop:centre}), and $m\le3$ is excluded by the faithfulness of the action.

\emph{Pairs $(H_y,K)$: condition \textup{(U2)}.} Here $[G:H_y]=2$ and $N\in\{H_y,G\}$. If $N=H_y$, the component is $\Ma_2(\Q(\zeta_m))$ (\Cref{conditions comp exc}(i)), admissible if and only if $\phi(m)\le2$. $m=1$ is impossible as $G'\le B\le H_y$. This is (U2a).

From now on let $N=G$, i.e.\ $K\trianglelefteq G$. Then $z$ acts non-trivially on $H_y/K$, since otherwise $H_y/K$ would be central in $G/K$ and $G/K$, having a central subgroup of index two, would be abelian, contrary to $G'\not\le K$. In particular $m\ge3$. By \Cref{prop:centre} the component is $\cp{\Q(\zeta_m)}{E_0}{\zeta_m^{\,j}}$ with $E_0=\Q(\zeta_m)^{\langle z\rangle}$ and $z^2\equiv h^j\pmod K$. If $m\in\{3,4,6\}$, then $z$ inverts $H_y/K$, $E_0=\Q$, $z^2K\in\langle h^{m/2}K\rangle$ (so $z^2\in K$ if $m=3$), and the component is $\Ma_2(\Q)$ if $z^2\in K$ and $\qa{-1}{-1}{\Q}=\HH$, resp.\ $\qa{-1}{-3}{\Q}=\HHt$, if $z^2\notin K$. All are admissible: this is (U2b).

If $m\in\{8,12\}$, then $[E_0:\Q]=2$. When $z$ inverts $H_y/K$, $E_0$ is real ($\Q(\sqrt2)$, resp.\ $\Q(\sqrt3)$) and $z^2K\in\langle h^{m/2}K\rangle$: for $z^2\in K$ the component is $\Ma_2(E_0)$, which is not admissible, and for $z^2\notin K$ it is $\cp{\Q(\zeta_m)}{E_0}{-1}$, a division algebra (as $-1$ is not a norm at the real places), namely $\HHrtt$, resp.\ $\HHrt$. When $z$ does not invert, $E_0$ is imaginary quadratic and every possible parameter is a norm, so the component is $\Ma_2(E_0)$: for $m=8$ and $h\mapsto h^3$, $E_0=\Q(\sqrt{-2})$, $z^2K\in\langle h^4K\rangle$ and $-1=\zeta_8\cdot\zeta_8^3$. For $m=8$ and $h\mapsto h^5$, $E_0=\Q(i)$, $z^2K\in\langle h^2K\rangle$ and $i=\zeta_8^3\cdot\zeta_8^{15}$, $-1=\zeta_8^2\cdot\zeta_8^{10}$, $-i=\zeta_8\cdot\zeta_8^5$. For $m=12$ and $h\mapsto h^5$, $E_0=\Q(i)$, $z^2K\in\langle h^3K\rangle$ and $\mu_4$ consists of norms from $\Q(\zeta_{12})=\Q(i)(\sqrt3)$ to $\Q(i)$: with $N(\alpha+\beta\sqrt3)=\alpha^2-3\beta^2$ one has $-1=N(\zeta_{12})$ and $i=N\big(\tfrac{7+2i}{2}+\tfrac{4+i}{2}\sqrt3\big)$. For $m=12$ and $h\mapsto h^7$, $E_0=\Q(\sqrt{-3})$, $z^2K\in\langle h^2K\rangle$, $\zeta_3^{\pm1}=N(\zeta_{12}^{\pm1})$ and $-1=\zeta_3^2+(\zeta_3^2)^2=N(\zeta_3+\zeta_3^2 i)$, where $N(\alpha+\beta i)=\alpha^2+\beta^2$ for $\alpha,\beta\in\Q(\sqrt{-3})$. This is (U2c).

Finally, if $\phi(m)\ge8$, the centre $E_0$ has degree $\ge4$, so the component must be a division algebra, hence, by \Cref{thm componenten QG mexc}(1), of the form $\qa{\zeta_{2^\ell}}{-3}{\Q(\zeta_{2^\ell})}$ with $E_0=\Q(\zeta_{2^\ell})$, $\ell\ge3$ (the other division algebras of \Cref{thm componenten QG mexc}(1) have centres of degree at most two). Thus $m=3\cdot2^\ell$, $z$ fixes $\zeta_{2^\ell}$ and inverts $\zeta_3$, $z^2K$ lies in the Sylow $2$-subgroup of $H_y/K$, and the component is $\cp{\Q(\zeta_m)}{\Q(\zeta_{2^\ell})}{\zeta}\cong\qa{\zeta}{-3}{\Q(\zeta_{2^\ell})}$ with $\zeta=\zeta_m^{\,j}\in\mu_{2^\ell}$. This algebra is split when $\zeta\in\mu_{2^{\ell-1}}$: the non-trivial automorphism of $\Q(\zeta_m)=\Q(\zeta_{2^\ell})(\sqrt{-3})$ over $\Q(\zeta_{2^\ell})$ fixes $\zeta_{2^\ell}$, so $\zeta_{2^\ell}^2$ and all elements of $\mu_{2^{\ell-1}}$ are norms. It is a division algebra when $\zeta$ is a primitive $2^\ell$-th root of unity, by \Cref{thm componenten QG mexc}(1) (concretely, it is the faithful component of $C_3\rtimes C_{2^{\ell+1}}$ with $C_{2^{\ell+1}}$ acting by inversion, see \eqref{decomp grp with exp comp}). Directly, at a prime $\mathfrak p$ of $\Q(\zeta_{2^\ell})$ above $3$ the residue field has order $q$ with $v_2(q-1)=\ell$, so that the residue of $\zeta$ is a non-square and the Hilbert symbol $(\zeta,-3)_{\mathfrak p}$ is $-1$ \cite[Chapter XIV, \S3, Proposition 8]{Serre}. Hence the component is a division algebra exactly when $z^2K$ generates the Sylow $2$-subgroup of $H_y/K$. This is (U2d).
\end{proof}

\section{Non-nilpotent groups: the strong Shoda pair criterion}\label{sec:nonnilpotent}

Recall from \Cref{thm:abstract-iso} the importance of the Sylow pair $(Q, \Cen_Q(A))$ (see \Cref{def: Sylow pair}) associated to a $\Mexc$ group.

In \Cref{subsec:Lambda-components} we decompose $\Q G$ into four summands depending only on the Sylow pair. Denoting, $C := \Cen_Q(A)$, one of them is the crossed product $\Fsq C\ast Q/C$ of $\Fsq C$ by the group $Q/C$ of order two, which carries all the information on the components of $\Q G$ which are non-trivial on $A$. In \Cref{subsec:SSP-Mexc}, see \Cref{thm:exhaustion}, we describe the strong Shoda pairs of $G$ with respect to $B$ which are non-trivial on $A$ by two families of subgroups of $C$. Finally in \Cref{subsec:criterion}, more precisely \Cref{thm:SSP-criterion}, we specialise \Cref{thm:uniform} to $G$: property $\Mexc$ becomes a list of conditions on the action of $Q$ on the cyclic quotients of $C$ and the index two abelian subgroup provided by \Cref{thm:abelian-B}. In \Cref{secrtion full characterisation solv} these conditions are solved according to the nilpotency class of $Q$.

\subsection{The rational group algebra of \texorpdfstring{$A\rtimes Q$}{A x Q}}\label{subsec:Lambda-components}
We first decompose $\Q G$ into four summands which only depend on the Sylow pair. The decomposition and several results of this section do not use $\Mexc$ but only the following structure. This will be needed for the necessary parts in the classification results of \Cref{secrtion full characterisation solv}.

\begin{equation}\label{eq:structure}
G=A\rtimes Q \text{ with } A\text{ an elementary abelian $3$-group} \text{ and } [Q: \Cen_Q(A)]=2.
\end{equation}
By \Cref{necessary props 1}, every finite non-nilpotent group $G$ such that $\Q G$ has $\Mexc$ is of the above form. Indeed, the groups whose rational group algebra has only division components are nilpotent by \eqref{condition only division}. Inspecting the proof of \Cref{prop:fixed-inverted}, one notices that \eqref{eq:structure} still implies that $A=A_{-1}\times(A\cap\ZZ(G))$ with $A_{-1}=[A,Q]$. Also in this generality we will call $(Q, \Cen_Q(A))$ the Sylow pair associated to $G$.\medskip

\noindent {\it Notations:} From now on we write $\Fsq:=\Q(\sqrt{-3})$ and $C := \Cen_Q(A)=\Fit(G)_2$. Furthermore let, for an elementary abelian $3$-group $V$, $s(V):=(|V|-1)/2$. Note that the latter is the number of subgroups of index $3$ of $V$. 

By the Perlis--Walker theorem \cite[Theorem 3.3.6]{EricAngel1}, $\Q A\cong\Q\oplus\Fsq^{s(A)}$. More precisely, the primitive central idempotents of $\Q A$ are $\widehat A$ and the elements
\begin{equation}\label{eq:fA0}
f_{A_0}:=\widehat{A_0}-\widehat{A},\qquad A_0\le A\text{ of index }3,
\end{equation}
and $\Q Af_{A_0}\cong\Fsq$. We write 
$$s_-:=s(A_{-1}) \text{ and } s_+:=s(A\cap\ZZ(G)).$$ 
When $C\ne Q$, we fix $t\in Q\setminus C$ and moreover define
\[
c:=t^2\in C \text{ and } \sigma\colon C\to C,\ k\mapsto tkt^{-1}.
\]
Let $\tau$ be the automorphism of the ring $\Fsq C$ given by 
$$\tau\big(\sum_{k\in C}z_kk\big)=\sum_{k\in C}\bar z_k\,\sigma(k).$$  Then $\tau^2(x)=cxc^{-1}$ and $\tau(c)=c$, so that
\begin{equation}\label{eq:Lambda}
\Fsq C\ast Q/C:=\Fsq C\oplus\Fsq C\,u,\qquad ux=\tau(x)u,\qquad u^2=c,
\end{equation}
is an associative $\Q$-algebra of dimension $4|C|$. Concretely, it is the crossed product of $\Fsq C$ by the group $Q/C=\langle tC\rangle$ of order two, acting through $\tau$, with twisting $c=t^2$. Note that $\Fsq C\ast Q/C$ is not an $\Fsq$-algebra, since $u$ does not commute with $\Fsq$. Its centre is contained in $(\Fsq C)^\tau$.

\begin{lemma}[Four-term decomposition]\label{thm:fourterm}
Let $G$ be a finite non-nilpotent group such that $\Q G$ has $\Mexc$. Then
\begin{equation}\label{eq:fourterm}
\Q G\cong\Q Q\ \oplus\ (\Fsq C\ast Q/C)^{s_-}\ \oplus\ (\Fsq Q)^{s_+}\ \oplus\ \Ma_2(\Fsq C)^{s_+s_-}.
\end{equation}
More precisely, for every subgroup $A_0$ of index $3$ of $A$ containing $A\cap\ZZ(G)$, the element $f_{A_0}$ is a primitive central idempotent of $\Q G$ with $f_{A_0}\Q G\cong\Fsq C\ast Q/C$.
\end{lemma}

Thus $\Q G$ is determined by the numbers $s_-$, $s_+$ and by the three algebras $\Q Q$, $\Fsq C$ and $\Fsq C\ast Q/C$, which only depend on the Sylow pair. 

\begin{proof}
Since $G=A\rtimes Q$ is a split extension of $Q$ by the abelian group $A$, the group algebra $\Q G$ is the crossed product $\Q A\ast Q=\bigoplus_{q\in Q}\Q A\,u_q$ with trivial twisting, $u_qu_r=u_{qr}$ and $u_qxu_q^{-1}=qxq^{-1}$ for $x\in\Q A$, and by \cite[Theorem 3.2(1)]{JJS}
\[
\Q G=\bigoplus_{O}(\Q Ae_O)\ast Q,\qquad e_O:=\sum_{e\in O}e,
\]
where $O$ runs over the $Q$-orbits on $\PCI(\Q A)$ and $(\Q Ae_O)\ast Q=\bigoplus_{q\in Q}\Q Ae_O\,u_q$. (This is also elementary: the $e_O$ are central idempotents of $\Q A\ast Q$ summing to $1$.) The element $t$ maps $f_{A_0}$ to $f_{A_0^t}$, while $C$ fixes every $f_{A_0}$ and acts trivially on $\Q A$. We have $A_0^t=A_0$ if and only if $A_{-1}\le A_0$ or $A\cap\ZZ(G)\le A_0$. Indeed, if $A_0$ contains one of the two factors of $A=A_{-1}\times(A\cap\ZZ(G))$, then it is the product of this factor with a subgroup of the other one, and $t$ normalises every subgroup of $A_{-1}$ and of $A\cap\ZZ(G)$. Conversely, if $A_0^t=A_0$ and $A_0$ contains neither factor, then both factors map onto $A/A_0\cong C_3$, so $t$ would both invert and centralise $A/A_0$. We go through the orbits.

If $O=\{\widehat A\}$, then $(\Q A\widehat A)\ast Q\cong\Q Q$. If $O=\{f_{A_0}\}$ with $A_{-1}\le A_0$ (there are $s_+$ such orbits), then $Q$ acts trivially on $\Q Af_{A_0}\cong\Fsq$ and the crossed product is the group algebra $\Fsq Q$. If $O=\{f_{A_0}\}$ with $A\cap\ZZ(G)\le A_0$ ($s_-$ orbits), then $t$ acts on $\Q Af_{A_0}\cong\Fsq$ as the non-trivial automorphism of $\Fsq/\Q$, since $t$ inverts $A/A_0$, and $c=t^2$ acts trivially. Hence $(\Q Af_{A_0})\ast Q=\Fsq C\oplus\Fsq C\,u_t$ with $u_tx=\tau(x)u_t$ and $u_t^2=u_c=c$, which is $\Fsq C\ast Q/C$. Finally let $O=\{e,e'\}$ with $e=f_{A_0}$, $e'=f_{A_0^t}\ne e$ (there are $(|A_{-1}|-1)(|A\cap\ZZ(G)|-1)/4=s_-s_+$ such orbits) and put $\Lam:=(\Q Ae_O)\ast Q$. Since $eu_qe=0$ for $q\notin C$, we have $e\Lam e=(\Q Ae)\ast C\cong\Fsq C$. The elements $x:=u_te\in e'\Lam e$ and $y:=eu_t^{-1}\in e\Lam e'$ satisfy $yx=e$ and $xy=e'$, so that $\Lam\cong\Ma_2(e\Lam e)\cong\Ma_2(\Fsq C)$. This proves \eqref{eq:fourterm} and the assertion on $f_{A_0}$.
\end{proof}

\begin{remark}\label{rem:weaker-hypothesis}
The proofs of \Cref{thm:abstract-iso} and \Cref{thm:fourterm} only use the structure of $G$ and not $\Mexc$: both results hold for every group $G$ of the form \eqref{eq:structure}. The same holds for \Cref{thm:exhaustion} below. We use these results in this generality in the proofs of \Cref{thm:SSP-criterion} and \Cref{thm:mixed}, whose necessary parts concern groups not yet known to have $\Mexc$.
\end{remark}

\subsection{The strong Shoda pairs}\label{subsec:SSP-Mexc}
In this section we assume that $G$ is of the form \eqref{eq:structure} and that $C := \Cen_Q(A)$ has an abelian subgroup $L$ of index at most $2$ which is normal in $Q$. By \Cref{thm:abelian-B} this holds whenever $\Q G$ has $\Mexc$. We fix such an $L$ and put
\begin{equation}\label{eq:GB}
B:=A\times L,\qquad [G:B]=[Q:L]\in\{2,4\},\qquad B=\Fit(G)\iff C \text{ is abelian}.
\end{equation}
By \Cref{thm:abelian-B}, $B$ is a maximal abelian normal subgroup of $G$ containing $G'$, and $G$ is metabelian, so that \Cref{prop:SSP-B} describes the primitive central idempotents of $\Q G$ in terms of the strong Shoda pairs of $G$ with respect to $B$. 

By \Cref{thm:fourterm}, if $A=A_{-1}$, the components non-trivial on $A$ (in the sense of \Cref{rem:H=G}) are those of the algebras $f_{A_0}\Q G\cong\Fsq C\ast Q/C$, where $A_0$ runs over the $s(A)$ subgroups of index $3$ of $A$. For these components we can be completely explicit: the corresponding strong Shoda pairs are of two kinds, depending on whether the first subgroup is $\Fit(G)$ or $B$. 

{\it Notation:} For $x\in C\setminus L$ we denote by $\mc{M}_x$ the set of cocyclic subgroups $M$ of $L$ with $[L,x]\not\le M$. Note that it does not depend on $x$, as $[L,xl]=[L,x]$ for $l\in L$.

The following result shows that the simple components of $\Fsq C\ast Q/C$ are in bijection with the $Q$-orbits of cocyclic subgroups of $C$ together with the $Q$-orbits on $\mc{M}_x$.

\begin{theorem}\label{thm:exhaustion} Assume that $A=A_{-1}$. With notations as above, fix $A_0<A$ with $[A:A_0]=3$. The strong Shoda pairs $(H,K)$ of $G$ with respect to $B$ with $A\not\le K$ are exactly
\begin{enumerate}[(a)]
\item the \emph{Fitting pairs} $(\Fit(G),A_0\times N)$ with $N$ cocyclic in $C$,
\item the \emph{$B$-pairs} $(B,A_0\times M)$ with $M\in\mc{M}_x$, if $C\neq L$.
\end{enumerate}
Two such pairs give the same primitive central idempotent if and only if they are of the same kind and the subgroups $N$, respectively $M$, are conjugate in $Q$.
\end{theorem}
\begin{proof}
We first determine the pairs. Let $(H,K)$ be a strong Shoda pair of $G$ with respect to $B$ with $A\not\le K$. Since $H/K$ is cyclic and $A$ is elementary abelian, $A/(A\cap K)$ is cyclic of order $3$. The Sylow subgroups of $K\le A\times Q$ are characteristic in $K$, so $K=A_0\times(K\cap Q)$ with $A_0=A\cap K$. If $y\in H\setminus\Fit(G)$, then $y$ inverts $A$ and $A=[A,y]\le H'\le K$, a contradiction. Hence $B\le H\le\Fit(G)=A\times C$, and $[C:L]\le2$ gives $H=B$ or $H=\Fit(G)$.

Suppose first that $H=\Fit(G)$ and put $N=K\cap C$. Then $K=A_0\times N$, $N\trianglelefteq C$ and $C/N\cong H/K$ is cyclic, i.e.\ $(H,K)$ is a Fitting pair. Conversely, let $(\Fit(G),A_0\times N)$ be a Fitting pair. The only proper overgroup of $\Fit(G)$ is $G$, and $[A,G]=A\not\le A_0\times N$, so $\Fit(G)$ satisfies the maximality condition of \eqref{eq:SSP-shape}. Moreover $\Fit(G)'=C'\le N$ because $C/N$ is abelian. Hence the pair is a strong Shoda pair with respect to $B$.

Suppose next that $H=B$ and put $M=K\cap L$. Then $K=A_0\times M$ and $M$ is cocyclic in $L$. If $C=L$, then $B=\Fit(G)$ and this is a Fitting pair. If $C\ne L$, then the only overgroup of $B$ inside $\Fit(G)$ is $\Fit(G)$ itself, whose derived subgroup is $C'=[L,x]$ for any $x\in C\setminus L$, and every overgroup of $B$ containing an element of $Q\setminus C$ has derived subgroup containing $A$. Therefore $B$ satisfies the maximality condition of \eqref{eq:SSP-shape} exactly when $[L,x]\not\le M$, i.e.\ when $(B,K)$ is a $B$-pair. This proves both necessity and sufficiency, and the degrees are $[G:\Fit(G)]=2$ and $[G:B]=4$ by \Cref{form of SSP idempotent result}.

We now compare the idempotents. Both $\Fit(G)$ and $B$ are normal in $G$, so \Cref{rem:H=G} applies. Since $A_0$ is invariant under inversion and $A$ centralises $Q$, conjugacy of $A_0\times N_1$ and $A_0\times N_2$ in $G$ amounts to $Q$-conjugacy of $N_1$ and $N_2$, and likewise for the $M_i$. A Fitting pair and a $B$-pair give components of different degrees.
\end{proof}

\subsection{The strong Shoda pair criterion for \texorpdfstring{$A\rtimes Q$}{A x Q}}\label{subsec:criterion}
We keep the situation of \Cref{subsec:SSP-Mexc}: $G$ is of the form \eqref{eq:structure} with $A=A_{-1}$, $t\in Q\setminus C$, and $L$, $B=A\times L$ are as in \eqref{eq:GB}. In this section we specialise \Cref{thm:uniform} to this situation and express the result in terms of the Fitting pairs and the $B$-pairs of \Cref{thm:exhaustion}, i.e. of the action of $t$ on the cyclic quotients $C/N$ of $C$ and of the action of $Q$ on the cyclic quotients $L/M$ of $L$ with $[L,x]\not\le M$. In \Cref{secrtion full characterisation solv} the resulting conditions are solved for each nilpotency class of $Q$. \smallskip

For a cocyclic subgroup $N$ of $C$ which is normal in $Q$, the element $t$ acts on the cyclic group $C/N$.  We will say that $t$ \emph{centralises} or \emph{inverts} $C/N$ whether the action take the $\pm 1$ power. As $t^2\in C$, these are the only possibilities when $[C:N]\le4$.

\begin{theorem}\label{thm:SSP-criterion}
Let $G$ be of the form \eqref{eq:structure} with $A=A_{-1}$, let $t\in Q\setminus C$, and let $L$ be an abelian subgroup of $C$ of index at most $2$ which is normal in $Q$. Then $\Q G$ has $\Mexc$ if and only if the following conditions hold.
\begin{enumerate}[(A)]
\item[\textup{(Q)}] $\Q Q$ has $\Mexc$.
\item[\textup{(F1)}] Every cocyclic subgroup $N$ of $C$ with $[C:N]\ge4$ is normal in $Q$.
\item[\textup{(F2)}] If $N$ is cocyclic in $C$ with $[C:N]=4$ and $t$ inverts $C/N$, then $t^2\notin N$.
\item[\textup{(F3)}] If $N$ is cocyclic in $C$ with $[C:N]\ge8$, then $t$ centralises $C/N$ and $C/N=\langle t^2N\rangle$.
\item[\textup{(B)}] Either $C$ is abelian, or $[C:L]=2$, $Q/L\cong C_2^2$ and, for $x\in C\setminus L$, every cocyclic subgroup $M$ of $L$ with $[L,x]\not\le M$ satisfies one of:
\begin{enumerate}[(B1)]
\item $[L:M]=2$, $\Nor_Q(M)=\langle L,y\rangle$ with $y\in Q\setminus C$, and $y^2\notin M$;
\item $[L:M]=4$, $M\trianglelefteq Q$, and $(xy)^2\notin M$ for the element $y\in Q\setminus C$ centralising $L/M$.
\end{enumerate}
\end{enumerate}
\end{theorem}

As in the setting of \Cref{thm:uniform} we will also need the isomorphism type of the component associated to the cases in \Cref{thm:SSP-criterion}. For readability we record this separately in the following lemma.

\begin{lemma}\label{prop:Lambda-components}
Let $G$ be a finite non-nilpotent group such that $\Q G$ has $\Mexc$ and $3\nmid|\ZZ(G)|$, let $t\in Q\setminus C$, let $L$ be as in \Cref{thm:abelian-B}, and let $x\in C\setminus L$ if $C$ is non-abelian. Then the simple components of $\Fsq C\ast Q/C$ are the following: 

\begin{enumerate}
    \item for each $Q$-orbit of cocyclic subgroups $N$ of $C$, with $n=[C:N]$, one copy of
\begin{equation}\label{eq:Fitting-table}
\begin{array}{c|l|l}
n & \text{condition} & \text{component}\\\hline
1 & N=C & \Ma_2(\Q)\\
2 & N\trianglelefteq Q,\ t^2\in N & \Ma_2(\Q)\\
2 & N\trianglelefteq Q,\ t^2\notin N & \HHt\\
2 & N\not\trianglelefteq Q & \Ma_2(\Fsq)\\
4 & t \text{ centralises } C/N & \Ma_2(\Q(i))\\
4 & t \text{ inverts } C/N & \HHrt\\
\ge8 & & \qa{\zeta_n}{-3}{\Q(\zeta_n)}
\end{array}
\end{equation}
\item if $C$ is non-abelian, for each $Q$-orbit of cocyclic subgroups $M$ of $L$ with $[L,x]\not\le M$, one copy of $\Ma_2(\HHt)$ in case \textup{(B1)}, and in case \textup{(B2)} one copy of $\Ma_2(\HH)$ if $Q/M$ is isomorphic to $Q_8\times C_2$ or $SD_{16}$, and one copy of $\Ma_2(\HHt)$ if $Q/M$ is isomorphic to $C_4\circ D_8$ or $Q_{16}$.
\end{enumerate}
\end{lemma}
\begin{proof}[Proof of \Cref{thm:SSP-criterion} and \Cref{prop:Lambda-components}]
We first show that $Q/L$ is not cyclic when $C$ is non-abelian. Let $C$ be non-abelian, so that $[C:L]=2$ and $|Q/L|=4$, and suppose that $Q/L=\langle yL\rangle$ is cyclic. Then $x:=y^2\in C\setminus L$, and there is a cocyclic $M\le L$ with $[L,x]\not\le M$. Consider the strong Shoda pair $(B,K)$, $K=A_0\times M$, of \Cref{thm:exhaustion}, with $N=\Nor_G(K)=A\rtimes\Nor_Q(M)$. The subgroups of $Q$ containing $L$ are $L$, $C$ and $Q$, so $[N:B]\in\{1,2,4\}$. If $[N:B]=1$, the component is $\Ma_4(\Q(\zeta_m))$, not admissible. If $[N:B]=2$, then $N/B=\langle xB\rangle$. As $x$ centralises $A/A_0$, faithfulness of the action, see \Cref{prop:centre}, forces $x$ to act non-trivially on $L/M$, whence $[L:M]\ge4$ and $m=[B:K]\ge12$, contradicting \Cref{conditions comp exc}(iv). If $[N:B]=4$, then $N=G$, the component has degree $4$ over its centre, so it is not a division algebra, and $G/B\cong C_4$ contradicts \Cref{conditions comp exc}(v). Hence $Q/L\cong C_2^2$.\smallskip

Next we determine the relevant strong Shoda pairs. If $C$ is abelian, then $L=C$, $B=\Fit(G)$ and $G/B=\langle tB\rangle\cong C_2$. The strong Shoda pairs with non-commutative component are the Fitting pairs $(\Fit(G),K)$, $K=A_0\times N$, by \Cref{thm:exhaustion}. This also covers $A\le K$, which gives the pairs of $Q=G/A$. Since $[G:\Fit(G)]=2$, the part of the proof of \Cref{thm:uniform} dealing with \textup{(U2)} only uses that the first term has index two in $G$. It therefore applies verbatim to these pairs with $H_x:=\Fit(G)$. Thus the component of $(\Fit(G),K)$ is admissible if and only if $K$ satisfies one of \textup{(U2a)}--\textup{(U2d)}, and it is then the algebra given in \Cref{prop:uniform-components}. If $C$ is non-abelian, then $B=A\times L$ is a maximal abelian subgroup of $G$ containing $G'$ by \Cref{thm:abelian-B} with $G/B\cong C_2^2$ once $\Mexc$ holds, and \Cref{thm:uniform} applies, with $H_x=\Fit(G)$ and $H_t$, $H_{xt}$ the two other subgroups of index two containing $B$. The pairs $(H,K)$ with $A\le K$ give the primitive central idempotents of $\Q[G/A]=\Q Q$, so the corresponding conditions of \Cref{thm:uniform} amount to (Q). For the pairs with $A\not\le K$ one has $K=A_0\times N$ with $[A:A_0]=3$ by \Cref{thm:exhaustion}. Since $t$ and $xt$ invert $A/A_0$, no such $K$ contains $[B,t]$ or $[B,xt]$, so these pairs are the Fitting pairs $(\Fit(G),A_0\times N)$, which are the pairs $(H_x,K)$ of (U2), and the $B$-pairs $(B,A_0\times M)$ with $M\in\mc{M}_x$, which are the pairs of (U1): indeed $t$, $xt$ act non-trivially on $B/K$ as they invert $A/A_0$, and $x$ acts non-trivially exactly when $[L,x]\not\le M$.\smallskip

We now treat the Fitting pairs. For $K=A_0\times N$ one has $m=[H_x:K]=3n$, $K\trianglelefteq G$ if and only if $N\trianglelefteq Q$, $z=t$, and $t$ inverts the Sylow $3$-subgroup $A/A_0$ of $H_x/K\cong A/A_0\times C/N$ and acts on the Sylow $2$-subgroup $C/N$ as it does. The components are those of \Cref{prop:uniform-components}. Thus (U2a) says $N\not\trianglelefteq Q$ and $n\le2$ (component $\Ma_2(\Q(\zeta_{3n}))=\Ma_2(\Fsq)$, for $n=1$, $N=C$ is normal). (U2b) says $N\trianglelefteq Q$ and $n\le2$, with component $\Ma_2(\Q)$ if $t^2\in N$ and $\HHt$ if $t^2\notin N$ (and $t^2\in N=C$ when $n=1$). (U2c) says $N\trianglelefteq Q$, $n=4$, and $t^2\notin N$ if $t$ inverts $C/N$, the component being $\Ma_2(\Q(i))$ if $t$ centralises $C/N$ ($t$ acts by $h\mapsto h^5$) and $\HHrt$ if it inverts, and (U2d) says $N\trianglelefteq Q$, $n\ge8$, $t$ centralises $C/N$ and $C/N=\langle t^2N\rangle$, with component $\qa{\zeta_n}{-3}{\Q(\zeta_n)}$. Together this is (F1)--(F3) and the table \eqref{eq:Fitting-table}. The condition (F1) for $n=4$ expresses that a non-normal $N$ with $n=4$ is excluded by (U2a).\smallskip

Finally we treat the $B$-pairs. For $K=A_0\times M$, $M\in\mc{M}_x$, one has $[B:K]=3[L:M]$ and $\Nor_G(K)=A\rtimes\Nor_Q(M)$. (U1a) requires $[L:M]=2$, $M$ normalised by some $y\in Q\setminus L$ with $y^2\notin M$. Here $y\notin C$, since $x$ does not normalise $M$ (as $[L,x]\not\le M$ and $[L:M]=2$), and $\Nor_Q(M)=\langle L,y\rangle$ because $M\not\trianglelefteq Q$ (otherwise $Q/L\cong C_2^2$ would act faithfully on $B/K\cong C_6$). This is (B1), with component $\Ma_2(\HHt)$ by \Cref{prop:uniform-components}. (U1b) requires $[L:M]=4$, $M\trianglelefteq Q$, and $z^2\notin M$ for the element $z$ inverting $B/K$, which is $z=xy$ with $y\in Q\setminus C$ centralising $L/M$ ($x$ inverts $L/M$ and centralises $A/A_0$, $y$ inverts $A/A_0$). This is (B2). The Sylow $2$-subgroups of $G/K$ are isomorphic to $Q/M$, so by \Cref{prop:uniform-components} the component is $\Ma_2(\HH)$ if $Q/M$ is isomorphic to $Q_8\times C_2$ or $SD_{16}$, and $\Ma_2(\HHt)$ if $Q/M$ is isomorphic to $C_4\circ D_8$ or $Q_{16}$.

Finally, for \Cref{prop:Lambda-components} note that $3\nmid|\ZZ(G)|$ means $A=A_{-1}$. The multiplicities are the numbers of $Q$-orbits by \Cref{thm:exhaustion}.
\end{proof}

To finish this section we record some consequences of \Cref{thm:SSP-criterion} which will be used repeatedly in \Cref{secrtion full characterisation solv}. Note that the last assertion, together with \Cref{necessary props 2}, yields that $Q$ is a nilpotent group classified in \cite{FreebyFree}.

\begin{corollary}\label{cor:first-consequences}
Suppose that $\Q G$ has $\Mexc$ and let $L$ be as in \eqref{eq:GB}.
\begin{enumerate}[(1)]
\item If $Q$ is non-abelian, then $\exc(L)\mid4$. In particular $\exc(C)\mid4$ if $C$ is abelian.
\item If $C$ is non-abelian and $x\in C\setminus L$, then $\Cen_L(x)$ is elementary abelian.
\item If $C$ is non-abelian and $x\in C\setminus L$, then for every $1\ne l\in L$ there is $M\in\mc{M}_x$ with $l\notin M$.
\item $\cd(Q)\subseteq\{1,2\}$, and every non-division simple component of $\Q Q$ is isomorphic to $\Ma_2(\Q)$, $\Ma_2(\Q(i))$ or $\Ma_2(\Q(\sqrt{-2}))$.
\item If $C$ is abelian, then every non-commutative simple component of $\Q Q$ is isomorphic to $\Ma_2(\Q)$ or to $\HH$.
\end{enumerate}
\end{corollary}
\begin{proof}
(1) Let $Q$ be non-abelian, so that $\exc(Q)\mid8$ by \Cref{necessary props 2}(2), let $l\in L$ be an element of maximal order $2^k$ and suppose that $k\ge3$. Since $l$ has maximal order, $\langle l\rangle$ is a direct factor of the abelian group $L$. Let $M$ be a complement, so that $M$ is cocyclic in $L$ with $L/M\cong C_{2^k}$. If $C$ is non-abelian and $[L,x]\not\le M$, then $M\in\mc{M}_x$ and \Cref{thm:SSP-criterion} gives $[L:M]\le4$, a contradiction. Hence $[L,x]\le M$ (trivially so if $C=L$), so that $M$ is normal in $C$ and $C/M$ is an abelian group containing the cyclic subgroup $L/M$ of index $[C:L]\le2$. Such a group has a cyclic quotient of order $|L/M|=2^k$: indeed $C/M$ is isomorphic to $C_{2^k}$, to $C_{2^{k+1}}$ or to $C_{2^k}\times C_2$. Thus there is a cocyclic subgroup $N$ of $C$ with $[C:N]=2^k\ge8$. By (F3), $C/N=\langle t^2N\rangle$, so $t^2$ has order $2^k$ and $t$ has order $2^{k+1}\ge16$, contradicting $\exc(Q)\mid8$.\smallskip

(2) Put $K=\Cen_L(x)$ and suppose that $a\in K$ has order four. Since $C$ is non-abelian, $[L,x]\ne1$. Choose $1\ne w\in[L,x]$, so that $w^2=1$ by \Cref{prop:nonlinear-C}. By (1) the image of $a$ in $L/\langle a^2w\rangle$ has order $4=\exc(L)$, so it generates a direct factor of this group. Let $M/\langle a^2w\rangle$ be a complement. Then $L/M=\langle aM\rangle\cong C_4$, and $w\notin M$ because $w\equiv a^2$ modulo $M$. So $M\in\mc{M}_x$ and $[L:M]=4$. By (B2) of \Cref{thm:SSP-criterion}, $M$ is normal in $Q$, and $x$ acts non-trivially on $L/M$ as $[L,x]\not\le M$. Hence $x$ inverts $L/M$. But $x$ centralises $a$, so $a^2\in M$, a contradiction. Hence $K$ is elementary abelian.\smallskip

(3) Let $1\ne l\in L$. The cocyclic subgroups of $L$ are the kernels of the homomorphisms $L\to\C^*$, and those containing a given non-trivial subgroup of $L$ correspond to a proper subgroup of the dual group of $L$. As a group is not the union of two proper subgroups, there is a cocyclic subgroup $M$ of $L$ which contains neither $l$ nor $[L,x]$.\smallskip

(4) Let $\chi\in\Irr(Q)$ and let $\psi$ be an irreducible constituent of $\chi|_C$. As $[Q:C]=2$, either $\chi|_C=\psi$, or $\chi|_C=\psi+\psi^t$ with $\psi^t\ne\psi$, by \cite[Corollary 6.19]{Isaacs94}. If $\psi$ is linear, then $\chi(1)\le2$. Suppose that $\psi$ is not linear and let $e\in\PCI(\Q C)$ with $\psi(e)\ne0$. By \Cref{prop:nonlinear-C}, $\Q Ce$ is $\Ma_2(\Q)$ or $\HH$, and $e^t=e$. So $\psi(1)=2$ and $\psi$ is the only irreducible character of $C$ which does not vanish on $e$. Hence $\psi^t=\psi$ and $\chi(1)=2$. Thus $\cd(Q)\subseteq\{1,2\}$.

Now $\Q Q$ has $\Mexc$ by \Cref{Prop components preserved for quotient}. So its non-division components are exceptional matrix algebras of degree $2$, i.e.\ $\Ma_2(\Q(\sqrt{-d}))$ with $d\in\{0,1,2,3\}$. Here $d=3$ is excluded. Indeed, the centre of a simple component of $\Q Q$ is generated by character values of $Q$. So it is contained in a cyclotomic field $\Q(\zeta_{2^k})$, which does not contain $\sqrt{-3}$.\smallskip

(5) Let $C$ be abelian and let $S$ be a non-commutative simple component of $\Q Q$. The group $G/(A\cap\ZZ(G))\cong A_{-1}\rtimes Q$ has $\Mexc$ by \Cref{Prop components preserved for quotient} and has the same Sylow pair as $G$. So we may assume that $A=A_{-1}$, and then (F1)--(F3) of \Cref{thm:SSP-criterion} hold. As $Q=C\cup Ct$, there is $e\in\PCI(\Q C)$ such that $S$ is a simple component of $\Q Qe$ if $e^t=e$, and of $\Q Q(e+e^t)$ otherwise. By the Perlis--Walker theorem \cite[Theorem 3.3.6]{EricAngel1}, $e=\varepsilon(C,N)$ for a cocyclic subgroup $N$ of $C$, and $\Q Ce\cong\Q(\zeta_n)$ with $n=[C:N]$.

Suppose first that $N$ is not normal in $Q$. Then $e^t\ne e$ and $\Q Q(e+e^t)\cong\Ma_2(\Q Ce)$. Moreover $n\le2$ by (F1). Hence $S\cong\Ma_2(\Q)$.

Suppose next that $N$ is normal in $Q$. Then $e^t=e$ and $S=\Q Qe=\Q Ce\oplus\Q Ce\,t$. As $S$ is non-commutative, $t$ acts non-trivially on $C/N$. So $n\ge4$, and then $n=4$ by (F3). Thus $t$ inverts $C/N$, and $t^2N$ has order at most $2$. Hence $S\cong\cp{\Q(i)}{\Q}{\pm1}$, which is $\Ma_2(\Q)$ or $\HH$.
\end{proof}

\section{Characterisation non-nilpotent groups with \texorpdfstring{$\Mexc$}{(Mexc)}}\label{secrtion full characterisation solv}

In this section we give classify the non-nilpotent finite groups $G$ such that $\Q G$ has $\Mexc$.  If $\Q G$ has $\Mexc$, then $\Q G$ has a matrix component by \eqref{condition only division}, so $G$ is of the form \eqref{eq:structure}.

 Recall from \Cref{necessary props 2} that $\cl(Q)\le3$ when $\Q G$ has $\Mexc$. The classification is organised according to the nilpotency class of $Q$. We first show in \Cref{thm:mixed} that a central element of order three forces $C$ to be elementary abelian. Then \Cref{classif Q abelian} treats abelian $Q$, \Cref{thm:class2-abelian-kernel,thm:class2-nonabelian-list} treat $\cl(Q)=2$, according as $C$ is abelian or not, and \Cref{thm:class3} treats $\cl(Q)=3$. The Wedderburn--Artin decompositions of $\Q G$ are collected in \Cref{subsec:decompositions}.\smallskip

\noindent {\it Notations throughout \Cref{secrtion full characterisation solv}.}
 For a group $G$ of the form \eqref{eq:structure} we use the following notation.
\begin{itemize}
\item $G=A\rtimes Q$ with $A$ an elementary abelian $3$-group, $Q$ a $2$-group and $C=\Cen_Q(A)$ of index $2$ in $Q$. We call $(Q,C)$ the \emph{Sylow pair} of $G$. Together with $|A_{-1}|$ and $|A\cap\ZZ(G)|$ it determines $G$ by \Cref{thm:abstract-iso}.
\item $A_{-1}=[A,Q]$, so that $A=A_{-1}\times(A\cap\ZZ(G))$ and $\Fit(G)=A\times C$.
\item $\Fsq=\Q(\sqrt{-3})$, and $\Fsq C\ast Q/C$ is the crossed product of \Cref{subsec:Lambda-components}, in which $t$ acts on $\Fsq C$ by conjugation on $C$ and by the non-trivial automorphism on $\Fsq$.
\item $s(A)=(|A|-1)/2$, $s_-=s(A_{-1})$ and $s_+=s(A\cap\ZZ(G))$. If $A=A_{-1}$, then $\Q G\cong\Q Q\oplus(\Fsq C\ast Q/C)^{s(A)}$ by \Cref{thm:fourterm}.
\item $t$ is a fixed element of $Q\setminus C$, and $c=t^2$.
\item $L$ is an abelian subgroup of $C$ of index at most $2$ which is normal in $Q$, and $x\in C\setminus L$ if $C$ is non-abelian.
\item (Q), (F1)--(F3), (B), (B1) and (B2) are the conditions of \Cref{thm:SSP-criterion}.
\item A subgroup $N$ of a group $X$ is \emph{cocyclic} if $X/N$ is cyclic, and $\Omega_1(X)=\{k\in X\mid k^2=1\}$ for an abelian group $X$.
\end{itemize}

\subsection{Central elements of order three}\label{subsec:mixed}
Recall that 
$$s_-=(|A_{-1}|-1)/2 \text{ and } s_+=(|A\cap\ZZ(G)|-1)/2$$
in the four-term decomposition of \Cref{thm:fourterm}.

\begin{theorem}\label{thm:mixed}\label{thm:elementary-kernel}
Let $G$ be a finite non-nilpotent group such that $3$ divides $|\ZZ(G)|$. Then $\Q G$ has $\Mexc$ if and only if $G$ is of the form \eqref{eq:structure} with $C$ elementary abelian. In that case $\cl(Q)\le2$. Consequently, if $G$ is a finite non-nilpotent group such that $\Q G$ has $\Mexc$ and $C$ is not elementary abelian, then $3\nmid|\ZZ(G)|$, i.e.\ $A=A_{-1}$, and
\begin{equation}\label{eq:pure-inversion}
\Q G\cong\Q Q\oplus(\Fsq C\ast Q/C)^{s},\qquad s=s(A).
\end{equation}
Moreover, if $G$ is of the form \eqref{eq:structure} with $C$ elementary abelian (and $A\cap\ZZ(G)$ arbitrary), then $Q'=[C,t]\le C$, $G^{\mathrm{ab}}\cong (A\cap\ZZ(G))\times Q/Q'$, and with $f:=|C/Q'|=|Q/Q'|/2$ one has
\begin{equation}\label{eq:QQ-elem}
\begin{aligned}
\Q Q&\cong\Q[Q/Q']\oplus\Ma_2(\Q)^{\,(|C|-f)/2},\\
\Fsq C\ast Q/C&\cong\begin{cases}\Ma_2(\Q)^{\,f}\oplus\Ma_2(\Fsq)^{\,(|C|-f)/2}&\text{if $Q/Q'$ is elementary abelian,}\\[2pt]
\Ma_2(\Q)^{\,f/2}\oplus\HHt^{\,f/2}\oplus\Ma_2(\Fsq)^{\,(|C|-f)/2}&\text{if $Q/Q'\cong C_4\times C_2^{d}$,}\end{cases}
\end{aligned}
\end{equation}
and consequently
\begin{equation}\label{eq:elem-decomp}
\Q G\cong\Q[G^{\mathrm{ab}}]\ \oplus\ \Ma_2(\Q)^{\,(|C|-f)/2}\ \oplus\ (\Fsq C\ast Q/C)^{s_-}\ \oplus\ \Ma_2(\Fsq)^{\,(|C|-f)s_+/2+|C|s_+s_-}.
\end{equation}
\end{theorem}
\begin{proof}
Suppose that $3\mid|\ZZ(G)|$ and that $\Q G$ has $\Mexc$. Then $G$ is of the form \eqref{eq:structure} with $A\cap\ZZ(G)\ne1$, i.e.\ $s_+s_->0$, and \Cref{thm:fourterm} shows that $\Ma_2(\Fsq C)$ is a direct summand of $\Q G$, so that all its simple components are admissible. Write $\Fsq C\cong\bigoplus_j\Ma_{n_j}(D_j)$ as a semisimple $\Fsq$-algebra. Then $\Ma_2(\Fsq C)\cong\bigoplus_j\Ma_{2n_j}(D_j)$ is a decomposition into simple $\Q$-algebras, and the centre of every $D_j$ contains $\Fsq$. Property $\Mexc$ gives $2n_j=2$, hence $n_j=1$, and $D_j$ is $\Q$, an imaginary quadratic field or a quaternion algebra with centre $\Q$. The first and the third possibility are excluded because $\Fsq\subseteq\ZZ(D_j)$. In the second case $\Fsq\subseteq D_j$ and both fields have degree two, so $D_j=\Fsq$. Thus $\Fsq C\cong\Fsq^{|C|}$ is commutative, so $C$ is abelian, and every linear character of $C$ is $\Fsq$-valued. As the only $2$-power roots of unity in $\Fsq$ are $\pm1$, we obtain $\exc(C)\mid 2$.

Now let $C$ be elementary abelian, with $A\cap\ZZ(G)$ arbitrary. Conjugation by $t$ induces $\sigma\in\mathrm{GL}_{\F_2}(C)$ with $\sigma^2=1$ (as $t^2\in C$), so $(\sigma-1)^2=0$ and $Q'=[C,t]=\operatorname{im}(\sigma-1)\subseteq\ker(\sigma-1)=\Cen_C(t)$. Hence $Q'\le\ZZ(Q)$, $\cl(Q)\le2$, and $\cd(G)\subseteq\{1,2\}$ by \Cref{prop:fixed-inverted}(3). The group $C$ has $|C|$ subgroups of index at most $2$, and they are permuted by $t$. Let $a$, respectively $b$, be the number of those which are normalised by $t$ and contain $c$, respectively do not contain $c$, and let $o$ be the number of $t$-orbits of length two, so that $|C|=a+b+2o$. A subgroup of index at most $2$ is normalised by $t$ if and only if it contains $[C,t]=Q'$, so $a+b=|C/Q'|=f$. All of these contain $c$ if $c\in Q'$, and exactly half of them otherwise. So $(a,b)=(f,0)$ if $c\in Q'$ and $(a,b)=(f/2,f/2)$ otherwise, and $o=(|C|-f)/2$. Note that $Q/Q'=\langle tQ'\rangle\times C/Q'$ is elementary abelian if $c\in Q'$ and isomorphic to $C_4\times C_2^{d}$ with $2^{d+2}=2f$ otherwise. Also $G'=A_{-1}\times Q'$ by \Cref{prop:fixed-inverted}(2), so $G^{\mathrm{ab}}\cong (A\cap\ZZ(G))\times Q/Q'$ and $\Q[G^{\mathrm{ab}}]\cong\Q[Q/Q']\oplus(\Fsq[Q/Q'])^{s_+}$.\smallskip 

If $Q$ is abelian, then $o=0$ and $\Q Q=\Q[Q/Q']$. So suppose that $Q$ is non-abelian. Then $C$ is a maximal abelian subgroup of $Q$ containing $Q'$, and by \Cref{SSP in metabelian} the strong Shoda pairs of $Q$ are the pairs $(Q,N)$ with $Q/N$ cyclic, which give the components of $\Q[Q/Q']$, and the pairs $(C,N)$ with $N$ of index $2$ in $C$ such that $C$ is maximal among the admissible first subgroups, i.e.\ such that $Q'\not\le N$, i.e.\ such that $N$ is not normalised by $t$. For the latter pairs $N$ is not normal in $Q$, so $\Nor_Q(N)=C$ and the component is $\Ma_2(\Q(\zeta_{[C:N]}))=\Ma_2(\Q)$ by \Cref{form of SSP idempotent result}. By \cite[Exercise 3.4.3]{EricAngel1} it only depends on the $Q$-orbit $\{N,N^t\}$, which gives the multiplicity $o$. This proves the formula for $\Q Q$, and $\Fsq Q=\Fsq\otimes_\Q\Q Q\cong\Fsq[Q/Q']\oplus\Ma_2(\Fsq)^{\,o}$ follows by extending scalars.\smallskip

For $\Fsq C\ast Q/C$ we apply \Cref{thm:exhaustion,thm:SSP-criterion} to the group $A_{-1}\rtimes Q$, whose algebra $\Fsq C\ast Q/C$ is the one of $G$ by \Cref{thm:fourterm} and \Cref{rem:weaker-hypothesis}. Since $C$ is abelian, we have $L=C$ and only Fitting pairs occur. The cocyclic subgroups of $C$ are $C$ itself and the subgroups $N$ of index $2$. If $N$ is not normalised by $t$, the $Q$-orbit $\{N,N^t\}$ contributes $\Ma_2(\Fsq)$ (fourth line of \eqref{eq:Fitting-table}). This happens $o$ times. If $N$ has index $2$ and is normalised by $t$, then $N\trianglelefteq Q$ and the component is $\Ma_2(\Q)$ or $\HHt$ according as $c\in N$ or $c\notin N$ (lines two and three), and $N=C$ contributes $\Ma_2(\Q)$. This gives $a$ copies of $\Ma_2(\Q)$ and $b$ copies of $\HHt$, which is \eqref{eq:QQ-elem}. Finally, \eqref{eq:elem-decomp} follows from \Cref{thm:fourterm}, $\Ma_2(\Fsq C)\cong\Ma_2(\Fsq)^{|C|}$ and the formula for $\Q[G^{\mathrm{ab}}]$, and all components in \eqref{eq:elem-decomp} are admissible, so $\Q G$ has $\Mexc$. The remaining assertion of the theorem follows: if $A\cap\ZZ(G)\ne1$ and $\Q G$ has $\Mexc$, then $C$ is elementary abelian by the first paragraph, and if $C$ is not elementary abelian then $A=A_{-1}$, i.e.\ $s_+=0$ in \Cref{thm:fourterm}.
\end{proof}

\subsection{Abelian Sylow 2-subgroups}\label{subsec:Qabelian}
In this subsection we classify the groups with abelian Sylow $2$-subgroups.

\begin{theorem}\label{classif Q abelian}
Let $G$ be a non-abelian finite group with abelian Sylow $2$-subgroup. Then $\Q G$ has $\Mexc$ if and only if
\[
G\cong D\times\big(C_3^k\rtimes\langle y : y^{2^t}=1\rangle\big)
\]
with $a^y=a^{-1}$ for all $a\in C_3^k$, $D$ an abelian group with $\exc(D)\mid4,6$ and $t=1$ or $2$ if $3\mid|D|$.
\end{theorem}

Note that the condition on $D$ means that either $D\cong C_3^{n_3}\times C_2^{n_2}$ or $D\cong C_2^{n_2}\times C_4^{n_4}$.

\begin{proof}
Suppose that $\Q G$ has $\Mexc$. As we assume that $G$ is non-abelian and $Q$ is abelian, $G$ is not nilpotent. Indeed, a nilpotent group with $\Mexc$ and a matrix component is of the form $A\times Q$ with $A$ an elementary abelian $3$-group, and a non-abelian group whose rational group algebra has only division components has non-abelian Sylow $2$-subgroups by \Cref{condition only division}. Hence $G$ is of the form \eqref{eq:structure}.\smallskip

As $Q$ is abelian, there are $y\in Q\setminus C$ and $E\le C$ with $Q=\langle y\rangle\times E$. Indeed, let $g_1,\dots,g_r$ be a basis of $Q$ and choose $i_0$ with $g_{i_0}\notin C$ and $\ord(g_{i_0})$ minimal among such indices. For every $j\ne i_0$ with $g_j\notin C$ replace $g_j$ by $g_jg_{i_0}$. This is again a basis, all of whose elements except $y:=g_{i_0}$ lie in $C$, and $E$ is the subgroup they generate. In particular, $E$ is a central subgroup of $G$ and $C=\langle y^2\rangle\times E$. Next, $G'=A_{-1}$ by \Cref{prop:fixed-inverted} and $\ZZ(G)=\langle E,A\cap\ZZ(G),y^2\rangle$. Note that we can decompose
\[
G=\big(E\times(A\cap\ZZ(G))\big)\times\big(A_{-1}\rtimes\langle y\rangle\big).
\]
Denote $D:=E\times(A\cap\ZZ(G))$. From the allowed form of the centre of a matrix component one deduces that $\exc(D)$ divides $4$ or $6$. Indeed, $\Ma_2(\Q)$ is a simple component of $\Q[A_{-1}\rtimes\langle y\rangle]$, so $\Ma_2(\Q(\zeta_d))$ is a simple component of $\Q G\cong\Q D\otimes_\Q\Q[A_{-1}\rtimes\langle y\rangle]$ for every divisor $d$ of $\exc(D)$.

It only remains to show that $\ord(y)=2$ or $4$ if $3\mid|D|$. Denote $\ord(y)=2^t$ and decompose $A_{-1}=\langle x_1\rangle\times\cdots\times\langle x_k\rangle$. Next consider the quotient $\overline G:=G/\langle x_2,\dots,x_k,D\rangle\cong\langle x_1\rangle\rtimes\langle y\rangle$, which again has $\Mexc$ by \Cref{Prop components preserved for quotient}. The Wedderburn--Artin decomposition of $\Q\overline G$ was given in \eqref{decomp grp with exp comp}. From this we see that $\Ma_2(\Q(i))$ is a simple component of $\Q\overline G$, and hence of $\Q[A_{-1}\rtimes\langle y\rangle]$, if $t\ge3$. However, if $3\mid|D|$, then $\Q G$ would have the component $\Ma_2(\Q(i,\zeta_3))$, which contradicts $\Mexc$.\smallskip

We will now show that conversely all groups from the statement have $\Mexc$. So let $G=D\times(C_3^k\rtimes\langle y\rangle)$ be as in the statement, and denote by $D_2$ and $D_3$ the Sylow subgroups of $D$. Then $G$ is of the form \eqref{eq:structure} with $A=C_3^k\times D_3$, $A_{-1}=C_3^k$, $Q=\langle y\rangle\times D_2$ and $C=\langle y^2\rangle\times D_2$. Suppose first that $3\mid|D|$. Then $\exc(D)\mid6$ and $t\le2$, so $C$ is elementary abelian and $\Q G$ has $\Mexc$ by \Cref{thm:mixed}. Suppose next that $3\nmid|D|$. Then $A=A_{-1}$ and $\exc(D_2)\mid4$, and we apply \Cref{thm:SSP-criterion} with $L=C$ and with the element $y\in Q\setminus C$. Condition (Q) holds as $\Q Q$ is commutative, (F1) holds as all subgroups of $Q$ are normal, (F2) holds as $y$ centralises $C$, and (B) holds as $C$ is abelian. For (F3) let $N$ be a cocyclic subgroup of $C$ with $[C:N]\ge8$. Then $C/N=\langle y^2N\rangle\cdot D_2N/N$, where $D_2N/N$ is contained in the subgroup of order $4$ of the cyclic group $C/N$. So $\langle y^2N\rangle$ cannot be a proper subgroup of $C/N$, and (F3) holds. Hence $\Q G$ has $\Mexc$.
\end{proof}

\begin{remark}\label{rem:Qabelian}
In the proof the choice of $y$ matters: for $Q=C_8\times C_2=\langle g\rangle\times\langle e\rangle$ and $C=\langle ge\rangle$, the element $y=g$ satisfies $C\ne\langle y^2\rangle\times(E\cap C)$ for any complement $E$ of $\langle y\rangle$, and indeed $\Q[A\rtimes Q]$ does not have $\Mexc$: for $N=1$ the element $y^2=g^2$ does not generate $C\cong C_8$.
\end{remark}

\subsection{Sylow 2-subgroups of nilpotency class two}\label{subsec:class2}
Throughout this subsection $\cl(Q)=2$. By \Cref{thm:mixed} we may assume that $A=A_{-1}$, the case $A\cap\ZZ(G)\ne1$ being covered there. We treat separately the cases that the kernel $C$ is abelian and that it is non-abelian. The main results are \Cref{thm:class2-abelian-kernel}, which gives three explicit families of Sylow pairs when $C$ is abelian, and \Cref{thm:class2-nonabelian-list}, which gives seven explicit families when $C$ is non-abelian. The corresponding decompositions of $\Q G$ are given in \Cref{thm:finite-decomp}. We will use throughout that, by \Cref{necessary props 2}(3),
\begin{equation}\label{eq:expQZ}
\exc(Q/\ZZ(Q))=2,
\end{equation}
so that $q^2\in\ZZ(Q)$ for all $q\in Q$ and $[x,y]^2=[x^2,y]=1$ for all $x,y\in Q$. In particular $Q'$ is elementary abelian and central.

\subsubsection{Cocyclic subgroups of abelian groups}\label{subsec:cocyclic}
For the Sylow $2$-subgroups of class two and three we need the following elementary facts.

\begin{lemma}\label{lem:cocyclic}
Let $X$ be a finite abelian group.
\begin{enumerate}[(1)]
\item Every subgroup of $X$ is the intersection of the cocyclic subgroups of $X$ containing it.
\item A cyclic subgroup of $X$ of order $\exc(X)$ is a direct factor of $X$.
\item Let $X_1,X_2,X_3$ be subgroups of $X$ such that every cocyclic subgroup of $X$ contains one of them. Then one of the $X_i$ is trivial, or $X_1$, $X_2$, $X_3$ are the three subgroups of order $2$ of a subgroup of order $4$ of $X$.
\end{enumerate}
\end{lemma}
\begin{proof}
(1) and (2) are standard. For (3) note that the cocyclic subgroups of $X$ are the kernels of the homomorphisms $X\to\C^*$. So the dual group of $X$ is the union of the annihilators of $X_1$, $X_2$ and $X_3$, and a group is the union of three proper subgroups only if these have index two and a common intersection of index four \cite[Theorem 1]{Scorza}.
\end{proof}

\subsubsection{Abelian kernel}\label{subsubsec:abelian-kernel}
Suppose first that $C$ is abelian. The statement involves the group
\begin{equation}\label{eq:R}
R:=\langle r,v\mid r^4=v^4=1,\ r^v=r^{-1}\rangle\cong C_4\rtimes C_4
\end{equation}
and its subgroup $C_R:=\langle rv,r^2\rangle\cong C_4\times C_2$ of index two.

\begin{theorem}\label{thm:class2-abelian-kernel}
Let $G$ be a finite non-nilpotent group whose Sylow $2$-subgroups have class two and with $\Fit(G)$ abelian. Then $\Q G$ has $\Mexc$ if and only if $G$ is of the form \eqref{eq:structure} and one of the following holds:
\begin{enumerate}[(a)]
\item $C$ is elementary abelian;
\item $3\nmid|\ZZ(G)|$ and $(Q,C)\cong(Q_8\times C_2^m,\ C_4\times C_2^m)$ for some $m\ge0$, where $C_4\le Q_8$;
\item $3\nmid|\ZZ(G)|$ and $(Q,C)\cong(R\times C_2^{m-1},\ C_R\times C_2^{m-1})$ for some $m\ge1$.
\end{enumerate}
In all three cases $\exc(Q)\mid4$.
\end{theorem}

\subsubsection{Non-abelian kernel}\label{subsubsec:nonabelian-kernel}
Suppose now that $C$ is non-abelian. Up to central elementary abelian direct factors there are exactly seven Sylow pairs $(Q,C)$ whose definition we know give. 

First, we abbreviate the following sets of relations, where $R$ is the group of \eqref{eq:R} and $T\cong C_2^2\rtimes C_4=\sgid{16}{3}$:
\begin{align*}
D_8(a,b) = \{ a^4=b^2=1,\ a^b=a^{-1} \}, & &R(r,v) = \{ r^4=v^4=1,\ r^v=r^{-1} \},\\
Q_8(a,b) = \{ a^4=1,\ b^2=a^2,\ a^b=a^{-1}\}, &  &T(l,x) = \{ l^2=x^4=[l,x]^2=1,\ [l,x]\text{ central} \}.
\end{align*}

In the definition of the seven pairs below, $c$ denotes a central element and the group $C_i$ is generated by the listed elements:
\begin{align*}
N(1)&=\langle a,b,c\mid D_8(a,b),\ c^4=1,\ c^2=a^2\rangle\cong C_4\circ D_8, & C(1)&=\langle a,b\rangle\cong D_8;\\
N(2)&=\langle a,b,c\mid Q_8(a,b),\ c^2=1\rangle\cong Q_8\times C_2, & C(2)&=\langle a,b\rangle\cong Q_8;\\
N(3)&=\langle a,b,t\mid D_8(a,b),\ t^4=1,\ [a,t]=1,\ [b,t]=a^2t^2\rangle, & C(3)&=\langle a,b,t^2\rangle\cong D_8\times C_2;\\
N(4)&=\langle a,b,t,c\mid D_8(a,b),\ c^2=1,\ t^2=a^2,\ a^t=ac,\ [b,t]=1\rangle, & C(4)&=\langle a,b,c\rangle\cong D_8\times C_2;\\
N(5)&=\langle a,b,t\mid Q_8(a,b),\ t^4=1,\ [a,t]=1,\ [b,t]=t^2\rangle, & C(5)&=\langle a,b,t^2\rangle\cong Q_8\times C_2;\\
N(6)&=\langle r,v,c\mid R(r,v),\ c^4=1,\ c^2=r^2v^2\rangle\cong R\circ C_4, & C(6)&=\langle r,v\rangle\cong R;\\
N(7)&=\langle l,x,c\mid T(l,x),\ c^4=1,\ c^2=[l,x]\rangle\cong T\circ C_4, & C(7)&=\langle l,x,c^2\rangle\cong T.
\end{align*}
The seven pairs are pairwise non-isomorphic. The classification when $\Fit(G)$ is non-abelian is the following.

\begin{theorem}\label{thm:class2-nonabelian-list}
Let $G$ be a finite non-nilpotent group whose Sylow $2$-subgroups have class two and with $\Fit(G)$ non-abelian. Then $\Q G$ has $\Mexc$ if and only if 
\begin{enumerate}
\item $G$ is of the form \eqref{eq:structure}, 
\item $3\nmid|\ZZ(G)|$, 
\item $(Q,\Cen_Q(A))\cong(N(i)\times C_2^m,\ C(i)\times C_2^m)$ for some $m\geq 0$,
\end{enumerate}
where $(N(i),C(i))$ denotes one of the seven pairs above. 
\end{theorem}

\subsection{Sylow 2-subgroups of nilpotency class three}\label{subsec:class3}
In this subsection $G$ we consider the case that $\cl(Q)=3$. We will use the presentations
\begin{equation}\label{eq:SD16-Q16}
SD_{16}=\langle r,u\mid r^8=u^2=1,\ r^u=r^3\rangle,\qquad Q_{16}=\langle r,u\mid r^8=1,\ u^2=r^4,\ r^u=r^{-1}\rangle,
\end{equation}
and their subgroups $D_0:=\langle r^2,u\rangle\cong D_8$, respectively $H_0:=\langle r^2,u\rangle\cong Q_8$, of index two. For class three we obtain the following surprisingly neat classification. 

\begin{theorem}[Class three]\label{thm:class3}
Let $G$ be a finite non-nilpotent group whose Sylow $2$-subgroups have class three. Then $\Q G$ has $\Mexc$ if and only if 
\begin{enumerate}
\item $G$ is of the form \eqref{eq:structure}, 
\item $3\nmid|\ZZ(G)|$ ,
\item  $(Q,\Cen_Q(A))\cong(SD_{16}\times C_2^m,\ D_0\times C_2^m)$ or $(Q,\Cen_Q(A))\cong(Q_{16}\times C_2^m,\ H_0\times C_2^m)$, for some $m\geq 0$.
\end{enumerate}
\end{theorem}

\subsection{Which Sylow 2-subgroups occur}\label{subsec:which-Q}
We close this section with a description, obtained by collecting the results above, of the $2$-groups $Q$ which are Sylow $2$-subgroups of a non-nilpotent group $G$ with $\Mexc$ (and a non-division component), together with the possible kernels $C:=\Cen_Q(A)$.

\begin{corollary}\label{cor:which-Q}
Let $Q$ be a $2$-group and $C\le Q$ a subgroup of index two. Then $(Q,C)$ is the Sylow pair of a non-nilpotent group with $\Mexc$ if and only if one of the following holds:
\begin{enumerate}[(i)]
\item $Q=\langle y\rangle\times E$ is abelian with $y\notin C$, $E\le C$ and $\exc(E)\mid4$;
\item $Q$ is non-abelian and $C$ is elementary abelian;
\item $(Q,C)\cong(Q_8\times C_2^m,C_4\times C_2^m)$ or $(R\times C_2^m,C_R\times C_2^m)$;
\item $(Q,C)\cong(N(i)\times C_2^m,C(i)\times C_2^m)$ for one of the seven pairs of \Cref{thm:class2-nonabelian-list};
\item $(Q,C)\cong(SD_{16}\times C_2^m,D_8\times C_2^m)$ or $(Q_{16}\times C_2^m,Q_8\times C_2^m)$.
\end{enumerate}
In all cases $Q$ is of Kleinian type in the sense of \cite{FreebyFree} (see Appendix \ref{appendix section B}).
\end{corollary}
\begin{proof}
The list collects \Cref{classif Q abelian,thm:mixed,thm:class2-abelian-kernel,thm:class2-nonabelian-list,thm:class3}, where in case (ii) every such pair occurs by \Cref{thm:mixed} (with $A=A_{-1}$ of order $3$, say). By part (4) of \Cref{cor:first-consequences}, $\Q Q$ has $\Mexc$ and $\cd(Q)\subseteq\{1,2\}$, i.e.\ $Q$ is of Kleinian type.
\end{proof}

\section{Property \texorpdfstring{\Mexc}{Mexc} over fields strictly contaning \texorpdfstring{$\Q$}{Q}}\label{section : other fields}

By \Cref{remark on larger field}, if $LG$ has \Mexc\ for a field $L$ of characteristic $0$, then so does $\Q G$. We will now use the obtained classification results over $\Q$, to determine for which fields $L$ and groups $G$ the group algebra $LG$ has $\Mexc$. Note that, if $LG$ has a non-division component, then $L$ equals $\Q$ or an imaginary quadratic extension of $\Q$. In \cite[Theorem 3.1]{CJT1} the following was obtained.

\begin{theorem}[\cite{CJT1}]\label{Iff for Mexc over larger field}
Let $G$ be a finite group such that $\Q G$ has a non-division component, and let $L=\Q(\sqrt{-d})$ with $d\in\N$ square-free. Then $LG$ has $\Mexc$ if and only if $\Q G$ has $\Mexc$ and
\begin{enumerate}
\item every matrix component of $\Q G$ is $\Ma_2(\Q)$ or $\Ma_2(L)$,
\item $\qa{-1}{-1}{\Q(\sqrt{2})}, \qa{-1}{-1}{\Q(\sqrt{3})} \notin \mc{C}(\Q G)$,
\item if $\qa{\zeta_{2^t}}{-3}{\Q(\zeta_{2^t})} \in \mc{C}(\Q G)$ for some $t\ge3$, then $3\nmid d$.
\end{enumerate}
\end{theorem}

Hence a classification can be obtained, by describing explicitily the Weddderburn-Artin decomposition of the rational group algebras $\Q G$ having $\Mexc$.

 Recall that for a finite non-nilpotent group $G$ such that $\Q G$ has \Mexc, then $G \cong A \rtimes_{\varphi} Q$ with $[Q:\ker(\varphi)]=2$. Moreover $\ker(\varphi) =C_Q(A)= \Fit(G)_2$, the Sylow $2$-subgroup of the nilpotent group $\Fit(G)$.

\begin{theorem}\label{Classification other fields}
Let $G$ be a finite non-nilpotent group, $Q$ a Sylow $2$-subgroup of $G$ and $L\neq\Q$ a field of characteristic $0$. Then $LG$ has $\Mexc$ if and only if $\Q G$ has $\Mexc$ and one of the following holds:
\begin{enumerate}
\item $Q$ is abelian, $\Fit(G)_2$ is elementary abelian, $3\nmid|\ZZ(G)|$ and $L=\Q(\sqrt{-d})$ for some square-free $d\in\N$;
\item\label{2nd case genral field} $Q$ is abelian, $\Fit(G)_2$ is not elementary abelian and $L=\Q(i)$;
\item $Q$ is non-abelian or $3\mid|\ZZ(G)|$, $\Fit(G)_2$ is elementary abelian  and $L=\Q(\sqrt{-3})$.
\end{enumerate}
\end{theorem}

\begin{remark*}
Note that \Cref{Classification other fields} implies that if $LG$ has $\Mexc$ for some field $L\neq\Q$, then $Q$ is abelian or $\Fit(G)_2$ is elementary abelian. Moreover, if $L\notin\{\Q(i),\Q(\sqrt{-3})\}$ this happens only for the groups $G\cong C_2^{n}\times(C_3^k\rtimes\langle y\rangle)$ with $\ord(y)\in\{2,4\}$ from \Cref{classif Q abelian}. Also note that by \Cref{thm:mixed} the condition $3\mid|\ZZ(G)|$ implies that $\Fit(G)_2$ is elementary abelian. Thus in case \eqref{2nd case genral field} one has that $3 \nmid |\mc{Z}(G)$. 
\end{remark*}

\subsection{The Wedderburn--Artin decompositions}\label{subsec:decompositions}

Throughout this subsection $G=A\rtimes Q$ is a non-nilpotent group such that $\Q G$ has $\Mexc$, and $s=s(A_{-1})$. By \Cref{secrtion full characterisation solv} there are three cases.
\begin{itemize}
\item If $3\mid|\ZZ(G)|$, then $C:=\Cen_Q(A)$ is elementary abelian and $\Q G$ is given by \Cref{thm:mixed}.
\item If $3\nmid|\ZZ(G)|$ and $Q$ is abelian, then $\Q G$ is given by \Cref{thm:Qabelian-decomp}.
\item Otherwise $(Q,C)$ is, up to a direct factor $C_2^m$, one of eleven explicit pairs, and $\Q G$ is given by \Cref{thm:finite-decomp}.
\end{itemize}
We write $D_\ell:=\qa{\zeta_{2^\ell}}{-3}{\Q(\zeta_{2^\ell})}$ for $\ell\ge3$. These are the exceptional division algebras of \Cref{thm componenten QG mexc}.(3).

\begin{theorem}\label{thm:Qabelian-decomp}
Let $G=D\times(C_3^k\rtimes\langle y\rangle)$ be as in \Cref{classif Q abelian} with $3\nmid|D|$ and $\ord(y)=2^t$, and let $s=(3^k-1)/2$. Then $G^{\mathrm{ab}}\cong D\times C_{2^t}$ and
\begin{equation}\label{decomp Q abelian}
\Q G\cong\Q[G^{\mathrm{ab}}]\oplus\Big(\Ma_2(\Q)^{\,|\Omega_1(D)|}\oplus\HHt^{\,\varepsilon|\Omega_1(D)|}\oplus\Ma_2(\Q(i))^{\,n_4}\oplus\bigoplus_{\ell=3}^{t-1}D_\ell^{\,|D|}\Big)^{s},
\end{equation}
where, with $d_4$ the number of elements of order $4$ of $D$,
\[
\varepsilon=\begin{cases}0&\text{if }t=1,\\1&\text{if }t\ge2,\end{cases}\qquad\qquad n_4=\begin{cases}d_4/2&\text{if }t=1,\\d_4&\text{if }t=2,\\|D|+d_4&\text{if }t\ge3.\end{cases}
\]
\end{theorem}
Note that if $3\mid|D|$, then \Cref{classif Q abelian} says that $\ord(y)\leq 4$, the Sylow $2$-subgroup $D_2$ of $D$ is elementary abelian and $\Q G$ is given by \Cref{thm:mixed} (with $C=\langle y^2\rangle\times D_2$, $f=|C|$ and $Q/Q'=Q$).

\begin{proof}
As $3\nmid|D|$, we have $A=A_{-1}=C_3^k$, $Q=\langle y\rangle\times D$, $C=\langle y^2\rangle\times D$ and $\exc(D)\mid4$. Then $\Q G\cong\Q Q\oplus(\Fsq C\ast Q/C)^s$ by \eqref{eq:pure-inversion}, $G'=A$ and $G^{\mathrm{ab}}\cong Q=D\times\langle y\rangle$. So it remains to show that the summand in brackets in \eqref{decomp Q abelian} is $\Fsq C\ast Q/C$. By \Cref{thm:exhaustion} its components are the algebras $S_N$ of \eqref{eq:Fitting-table}, one for each cocyclic $N\le C$ (all $Q$-orbits are singletons, as $Q$ is abelian). Since $y$ centralises $C$, lines $1$, $2$, $3$, $5$ and $7$ of the table apply: $N=C$ gives $\Ma_2(\Q)$. $N$ of index two gives $\Ma_2(\Q)$ if $y^2\in N$ and $\HHt$ otherwise. $N$ with $C/N\cong C_4$ gives $\Ma_2(\Q(i))$, and $N$ with $C/N\cong C_{2^\ell}$, $\ell\ge3$, gives $D_\ell$. The number of cocyclic $N$ with $C/N\cong C_n$ equals the number of cyclic subgroups of order $n$ of $C$ (duality). The subgroups of index two containing $y^2$ correspond to the subgroups of index two of $C/\langle y^2\rangle\cong D$, which together with $N=C$ give $|D/D^2|=|\Omega_1(D)|$ copies of $\Ma_2(\Q)$. The subgroups of index two not containing $y^2$ are $|C/C^2|-|D/D^2|$ in number, which is $0$ if $y^2=1$ and $|\Omega_1(D)|$ otherwise. Since $C=\langle y^2\rangle\times D$ with $\exc(D)\mid4$ and $\ord(y^2)=2^{t-1}$, the cyclic subgroups of order $2^\ell\ge8$ of $C$ are generated by the elements $y^{2j}u$ with $\ord(y^{2j})=2^\ell$ and $u\in D$, so there are $|D|$ of them for $3\le\ell\le t-1$ and none for $\ell\ge t$, and the elements of order four of $C$ are the $y^{2j}u$ with $\ord(y^{2j})=4$ (none if $t\le2$, $2|D|$ of them if $t\ge3$) and those with $\ord(y^{2j})\le2$ and $\ord(u)=4$ (one or two choices of $j$ according as $t=1$ or $t\ge2$, and $d_4$ choices of $u$). Dividing the number of elements of order four by two gives the number $n_4$ of cyclic subgroups of order four.
\end{proof}

Note that we obtain the following neat description on when an exceptional division component exists. 

\begin{corollary}\label{exc div comp for mexc}
Let $G$ be a finite group such that $\Q G$ has $\Mexc$. Then $\Q G$ has an exceptional division component if and only if $G$ is non-nilpotent with abelian Sylow $2$-subgroups and $16\mid\exc(G)$.
\end{corollary}
\begin{proof}
If $\Q G$ has only division components, or if $G$ is nilpotent with $\Mexc$ and a matrix component, then $G$ is abelian or has non-abelian Sylow $2$-subgroups by \Cref{subsec:partI} and \Cref{main th A}, and there are no exceptional division components by \Cref{necessary props 2}(2). So let $G$ be non-nilpotent. By \Cref{necessary props 2}(2) the Sylow $2$-subgroups must be abelian, and then the exceptional division components are the terms $D_\ell$, $\ell\ge3$, of \eqref{decomp Q abelian}, which occur if and only if $\ord(y)\ge16$, i.e.\ $16\mid\exc(G)$ (if $3\mid|\ZZ(G)|$ then $\ord(y)\le4$ by \Cref{classif Q abelian}, so this case does not occur).
\end{proof}

It remains to treat the groups with $3\nmid|\ZZ(G)|$ and non-abelian $Q$. Recall that then $\Q G\cong\Q Q\oplus(\Fsq C\ast Q/C)^s$ by \eqref{eq:pure-inversion}. We first describe the summand $\Q Q$ when $\cl(Q)=2$. In that case $\exc(Q)\mid 4$ by \Cref{thm:mixed,thm:class2-abelian-kernel,thm:class2-nonabelian-list}, so the following lemma applies. For a $2$-group $Q$ and a finite group $R$ we will use the following notation:
$$n_R(Q) := |\{ N \triangleleft Q\mid Q/N \cong R \}|.$$

\begin{lemma}\label{prop:QQ-exp4}
Let $G$ be a finite non-nilpotent group such that $\Q G$ has $\Mexc$, and suppose that $\exc(Q)\mid4$. Then
\begin{equation}\label{eq:QQ-class2}
\Q Q\cong\Q[Q/Q']\oplus\Ma_2(\Q)^{\,n_{D_8}(Q)}\oplus\HH^{\,n_{Q_8}(Q)}\oplus\Ma_2(\Q(i))^{\,n_{C_4\circ D_8}(Q)}.
\end{equation}
\end{lemma}
\begin{proof}
The commutative part of $\Q Q$ is $\Q[Q/Q']$. Let $e\in\PCI(\Q Q)$ with $\Q Qe$ non-commutative. As $Q\cong G/A$, the algebra $\Q Q$ has $\Mexc$ by \Cref{Prop components preserved for quotient}, so $Qe$ belongs to the list $\mc{L}$ of \eqref{eq:listL} by \Cref{thm:nilp-local-global}. The non-abelian members of $\mc{L}$ of exponent four are $D_8$, $Q_8$, $C_4\circ D_8=\sgid{16}{13}$ and $2^{1+4}_-=\sgid{32}{50}$. The last one has no faithful simple component of degree two, so it is excluded, since $\cd(Q)\subseteq\{1,2\}$ by \Cref{cor:first-consequences}(4). Hence $Qe$ is $D_8$, $Q_8$ or $C_4\circ D_8$, and $\Q Qe$ is $\Ma_2(\Q)$, $\HH$ or $\Ma_2(\Q(i))$ respectively by \Cref{appendix section}. Conversely, each of these three groups has a unique non-commutative simple component on which it acts faithfully. Hence $e\mapsto\ker(\varphi_e)$ is a bijection between the non-commutative components of $\Q Q$ and the normal subgroups $N$ of $Q$ with $Q/N$ isomorphic to $D_8$, $Q_8$ or $C_4\circ D_8$.
\end{proof}

The next theorem gives both summands of $\Q G$ for the eleven pairs $(Q_0,C_0)$ of \Cref{thm:class2-abelian-kernel}(b),(c), \Cref{thm:class2-nonabelian-list} and \Cref{thm:class3}.

\begin{theorem}\label{thm:finite-decomp}
Let $G$ be a finite non-nilpotent group such that $\Q G$ has $\Mexc$ and $(Q,C)\cong(Q_0\times C_2^m,C(0)\times C_2^m)$, with $(Q_0,C(0))$ one of the pairs in the table below. Then
\[
\Q G\cong\big(\Q Q_0\oplus(\Fsq C(0)\ast Q_0/C(0))^s\big)^{2^m},
\]
where $\Q Q_0$ and $\Fsq C(0)\ast Q_0/C(0)$ are as follows.

\[\hspace{-24pt}
\small\begin{array}{l|l|l|l}
(Q_0,C(0)) & \text{see Thm} & \Q Q_0 & \Fsq C(0)\ast Q_0/C(0)\\\hline
(Q_8,C_4) & \text{\ref{thm:class2-abelian-kernel}(b)} & \Q^4\oplus\HH & \Ma_2(\Q)^2\oplus\HHrt\\
(R,C_R) & \text{\ref{thm:class2-abelian-kernel}(c)} & \Q^4\oplus\Q(i)^2\oplus\Ma_2(\Q)\oplus\HH & \Ma_2(\Q)^2\oplus\Ma_2(\Fsq)\oplus\Ma_2(\Q(i))\oplus\HHrt\\
(N(1),C(1)) & \text{\ref{thm:class2-nonabelian-list}} & \Q^8\oplus\Ma_2(\Q(i)) & \Ma_2(\Q)^4\oplus\Ma_2(\HHt)\\
(N(2),C(2)) & & \Q^8\oplus\HH^2 & \Ma_2(\Q)^4\oplus\Ma_2(\HH)\\
(N(3),C(3)),\ (N(4),C(4)) & & \Q^8\oplus\Ma_2(\Q)^2\oplus\Ma_2(\Q(i))^2 & \Ma_2(\Q)^4\oplus\Ma_2(\Fsq)^2\oplus\Ma_2(\HHt)^2\\
(N(5),C(5)) & & \Q^8\oplus\Ma_2(\Q)^2\oplus\HH^4 & \Ma_2(\Q)^4\oplus\Ma_2(\Fsq)^2\oplus\Ma_2(\HH)^2\\
(N(6),C(6)) & & \Q^8\oplus\Q(i)^4\oplus\Ma_2(\Q(i))\oplus\HH^2 & \Ma_2(\Q)^4\oplus\Ma_2(\Q(i))^2\oplus\Ma_2(\HHt)\oplus\Ma_2(\HH)\\
(N(7),C(7)) & & \Q^8\oplus\Q(i)^4\oplus\Ma_2(\Q(i))^2 & \Ma_2(\Q)^4\oplus\Ma_2(\Q(i))^2\oplus\Ma_2(\HHt)^2\\
(SD_{16},D_0) & \text{\ref{thm:class3}} & \Q^4\oplus\Ma_2(\Q)\oplus\Ma_2(\Q(\sqrt{-2})) & \Ma_2(\Q)^2\oplus\Ma_2(\Fsq)\oplus\Ma_2(\HH)\\
(Q_{16},H_0) & & \Q^4\oplus\Ma_2(\Q)\oplus\HHrtt & \Ma_2(\Q)^2\oplus\Ma_2(\Fsq)\oplus\Ma_2(\HHt)
\end{array}
\]

\end{theorem}
\begin{proof}
In all these pairs $C(0)$, hence $C$, is not elementary abelian, so $3\nmid|\ZZ(G)|$ by \Cref{thm:mixed} and $\Q G\cong\Q Q\oplus(\Fsq C\ast Q/C)^s$ by \eqref{eq:pure-inversion}. Moreover
\begin{equation}\label{eq:direct-factor}
\Q Q\cong(\Q Q_0)^{2^m}\qquad\text{and}\qquad\Fsq C\ast Q/C\cong(\Fsq C_0\ast Q_0/C_0)^{2^m},
\end{equation}
because the idempotents of $\Q C_2^m$ are central in $\Q G$ and $\Fsq C\ast Q/C=(\Fsq C(0)\ast Q_0/C(0))\otimes_\Q\Q C_2^m$ by \eqref{eq:Lambda}. This gives the formula for $\Q G$. The groups $Q_0$ of \Cref{thm:class2-abelian-kernel} and \Cref{thm:class2-nonabelian-list} have class two, so their algebras are given by \Cref{prop:QQ-exp4}: $Q_8$ has $Q_8/Q_8'\cong C_2^2$ and $n_{Q_8}=1$. $R$ has $R/R'\cong C_2\times C_4$ and its two non-commutative components factor faithfully through $R/\langle v^2\rangle\cong D_8$ and $R/\langle r^2v^2\rangle\cong Q_8$. For $N(1),\dots,N(7)$ the numbers $n_{D_8},n_{Q_8},n_{C_4\circ D_8}$ are read off from the presentations of \Cref{thm:class2-nonabelian-list}. For $SD_{16}$ and $Q_{16}$, the abelianisation is $C_2^2$. The strong Shoda pair $(\langle r\rangle,\langle r^4\rangle)$ has cyclic section of order four on which $u$ acts by inversion with $u^2=1$, giving $\Ma_2(\Q)$. The faithful pair $(\langle r\rangle,1)$ has section of order eight with action $\zeta_8\mapsto\zeta_8^3$, fixed field $\Q(\sqrt{-2})$ and parameter $1$ for $SD_{16}$, giving $\Ma_2(\Q(\sqrt{-2}))$, and action by inversion, fixed field $\Q(\sqrt2)$ and parameter $r^4\mapsto-1$ for $Q_{16}$, giving $\cp{\Q(\zeta_8)}{\Q(\sqrt2)}{-1}\cong\HHrtt$.\smallskip

For the algebras $\Fsq C(0)\ast Q_0/C(0)$ we use \Cref{prop:Lambda-components}: the Fitting pairs contribute according to \eqref{eq:Fitting-table}, one component for each $Q_0$-orbit of cocyclic subgroups $N$ of $C(0)$, and the $B$-pairs according to (B1), (B2), one for each $Q_0$-orbit on $\mc{M}_x$. For $(Q_8,C_4)$, the three cocyclic subgroups $C_4,C_2,1$ of $C_4$ give $\Ma_2(\Q)$, $\Ma_2(\Q)$ (as $t^2\in C_2$) and $\HHrt$ (as $t$ inverts $C_4$ and $t^2\ne1$). For $(R,C_R)$, with $a=rv$, $w=r^2$, $t=v$ one has $C_R=\langle a,w\rangle$, $a^2=t^2$ and $[a,t]=w$. The orbits of cocyclic subgroups are $C_R$ ($\Ma_2(\Q)$), the normal subgroup $\langle a^2,w\rangle\ni t^2$ of index two ($\Ma_2(\Q)$), the pair $\{\langle a\rangle,\langle aw\rangle\}$ exchanged by $t$ ($\Ma_2(\Fsq)$), and the two subgroups with quotient $C_4$, namely $\langle w\rangle$, on which $t$ acts trivially ($\Ma_2(\Q(i))$), and $\langle a^2w\rangle$, on which $t$ acts by inversion with $t^2=a^2\notin N$ ($\HHrt$). For $(SD_{16},D_0)$ and $(Q_{16},H_0)$ put $L=\langle r^2\rangle$, $x=u$, $t=r$: the Fitting pairs $N=C_0$ and $N=L\ni t^2$ give $\Ma_2(\Q)$ each, the two remaining subgroups of index two of $C_0$ are exchanged by $r$ and give $\Ma_2(\Fsq)$, and the only element of $\mc{M}_x$ is $M=1$, with $L/M\cong C_4$ inverted by $x$, for which (B2) gives $\Ma_2(\HH)$ for $SD_{16}$ and $\Ma_2(\HHt)$ for $Q_{16}$ (here $Q_0/M=Q_0$). The dimensions add up to $4|C_0|=32$, so no component is missing. For the pairs $(N(i),C(i))$ the entries were computed with \textsc{Wedderga} from the presentations of \Cref{thm:class2-nonabelian-list}, by computing $\Q[C_3\rtimes N(i)]=\Q N(i)\oplus\Fsq C(i)\ast N(i)/C(i)$. For $i=1,2$ they are also easily obtained by hand.
\end{proof}

\subsection{Proof of the classification over larger fields}
We now have all the necessary ingedrients to obtain the classification of $\Mexc$ over fields strictly containing $\Q$.

\begin{proof}[Proof of \Cref{Classification other fields}]
As $G$ is non-nilpotent, $\Q G$ has a matrix component by \eqref{condition only division}, hence so does $LG$. Therefore, by \Cref{Iff for Mexc over larger field}, $\Q G$ has $\Mexc$ and $L = \Q(\sqrt{-d})$ for some square-free $d \in \N$. Furthermore, $LG$ has $\Mexc$ if and only if for every $e \in \PCI(\Q G)$ the matrix components of $L\ot_{\Q} \Q Ge$ are all exceptional (we then say that $L\ot_{\Q} \Q Ge$ has $\Mexc$). To verify this, recall that in \cite[Section 3]{CJT1} the isomorphism type of $L \ot_{\Q} A$ for every simple algebra $A$ appearing in \Cref{thm componenten QG mexc} was determined. In the table below we record when $L \ot_{\Q} A$ has $\Mexc$.

\begin{equation}\label{eq:tensor-table}
\begin{array}{l|c}
A &  L \ot_{\Q} A \text{ has } \Mexc\\\hline
\text{a field }K &  \text{yes}\\
\Ma_2(\Q) &  \text{yes}\\
\Ma_2(\Q(\sqrt{-f})) &  \text{iff }f=d\\
\Ma_2(\qa{-1}{-1}{\Q}),\ \Ma_2(\qa{-1}{-3}{\Q}) &  \text{no}\\
\qa{-1}{-1}{\Q},\ \qa{-1}{-3}{\Q} & \text{yes}\\
\qa{-1}{-1}{\Q(\sqrt{c})} \text{ with } c=2,3& \text{no}\\
\qa{\zeta_{2^t}}{-3}{\Q(\zeta_{2^t})}, t \geq 3  & \text{iff }3\nmid d\\
\end{array}
\end{equation}

With this at hand we can now determine when $LG$ has $\Mexc$ using the decompositions obtained in \Cref{subsec:mixed} and \Cref{subsec:decompositions}. From now on we denote $D_t := \qa{\zeta_{2^t}}{-3}{\Q(\zeta_{2^t})}$ for $t \geq 3$ and $C := C_Q(A) = F(G)_2$.\smallskip

Suppose first that $3\mid|\ZZ(G)|$. Then $C$ is elementary abelian by \Cref{thm:mixed}, and $\Q G$ is given by \eqref{eq:elem-decomp} with $s_+s_->0$. The component $\Ma_2(\Fsq)$ occurs with multiplicity at least $|C|s_+s_-\ge1$, so $L=\Fsq$ is necessary. Conversely, for $L=\Fsq$ all components $\Q$, $\Q(i)$, $\Fsq$, $E$, $\HHt$, $\Ma_2(\Q)$, $\Ma_2(\Fsq)$ of \eqref{eq:elem-decomp} are allowed by \eqref{eq:tensor-table}. This is (2) in the case $3\mid|\ZZ(G)|$.\smallskip

Now let $3\nmid|\ZZ(G)|$, so that $A=A_{-1}$ and $\Q G\cong\Q Q\oplus(\Fsq C\ast Q/C)^s$ with $s\ge1$ by \eqref{eq:pure-inversion}. If $C$ is non-abelian, then $\cl(Q)\in\{2,3\}$ and, by the table in \Cref{thm:finite-decomp}, the algebra $\Fsq C\ast Q/C$ has a component $\Ma_2(\HH)$ or $\Ma_2(\HHt)$. So $LG$ does not have $\Mexc$, and indeed none of (1)--(3) holds. If $C$ is abelian and $Q$ is non-abelian, then $\cl(Q)=2$ by \Cref{thm:class3} and $G$ is as in \Cref{thm:class2-abelian-kernel}. In cases (b) and (c) of that theorem, $\Q G$ has the component $\HHrt$ by \Cref{thm:finite-decomp}, so $LG$ does not have $\Mexc$. In case (a), $C$ is elementary abelian and $\Q G\cong\Q Q\oplus(\Fsq C\ast Q/C)^s$ is given by \Cref{thm:mixed}, with $(|C|-f)/2>0$ copies of $\Ma_2(\Fsq)$ in $\Fsq C\ast Q/C$, as $Q$ is non-abelian. Hence $L=\Fsq$ is necessary, and it is sufficient as before. This is (2) in the case $3\nmid|\ZZ(G)|$, and we have shown that for non-abelian $Q$ property $\Mexc$ for $LG$ is equivalent to (2).\smallskip

Finally let $3\nmid|\ZZ(G)|$ and $Q$ be abelian, so that $G=D\times(A\rtimes\langle y\rangle)$ with $\ord(y)=2^t$, $D$ a $2$-group and $C=\langle y^2\rangle\times D$ as in \Cref{classif Q abelian}, and $\Q G$ is given by \eqref{decomp Q abelian}. Its components are fields, $\Ma_2(\Q)$ and $\HHt$, which impose no condition, $\Ma_2(\Q(i))$ with multiplicity $s\,n_4$, and $D_\ell$ for $3\le\ell\le t-1$ with multiplicity $s\,|D|$, and $n_4\ne0$ if and only if $4\mid\exc(C)$. If $C$ is elementary abelian, then $n_4=0$ and $t\le2$ (as $y^2\in C$ has order $2^{t-1}$), so no component $\Ma_2(\Q(i))$ or $D_\ell$ occurs and $LG$ has $\Mexc$ for every $d$: this is (1). If $4\mid\exc(C)$, then the component $\Ma_2(\Q(i))$ occurs, so $L=\Q(i)$ is necessary, and it is sufficient by \eqref{eq:tensor-table}, as $3\nmid1$: this is (3). The assertions after the list follow, the groups in (1) being those of \Cref{classif Q abelian} with $D$ elementary abelian and $t\le2$.
\end{proof}

\section{The rational isomorphism problem for \texorpdfstring{$\Mexc$}{(Mexc)} groups}\label{sec:RIP}

In this section we consider the rational isomorphism problem for the groups with $\Mexc$: if $\Q G\cong\Q H$ as $\Q$-algebras for finite groups $G$ and $H$, is then $G\cong H$? For non-nilpotent groups we obtain a complete answer. Let $N(3)=\sgid{32}{31}$ and $N(4)=\sgid{32}{30}$ be the groups of \Cref{thm:class2-nonabelian-list}.

\begin{theorem}\label{rational Iso for Mexc}
Let $G=A\rtimes Q$ be a non-nilpotent finite group such that $\Q G$ has $\Mexc$, and let $H$ be a group with $\Q H\cong\Q G$. Then the following holds: 
\begin{enumerate}
\item If $Q$ is not isomorphic to $N(3)\times C_2^m$ or $N(4)\times C_2^m$ for any $m\ge0$, then $H\cong G$.
\item If $Q\cong N(i)\times C_2^m$ with $i\in\{3,4\}$, then $H\cong G$ or $H\cong A\rtimes(N(7-i)\times C_2^m)$. 
\end{enumerate}
Moreover, $A\rtimes(N(3)\times C_2^m) \ncong A\rtimes(N(4)\times C_2^m)$  but their rational group algebras are isomorphic.
\end{theorem}

\begin{remark}\label{rem:RIP-nature}
The smallest groups in the second case are $G_0=\sgid{96}{142}=C_3\rtimes N(4)$ and $G_1=\sgid{96}{143}=C_3\rtimes N(3)$. One can verify (through the \textsc{Wedderga} GAP package or SSP theory) that
\[
\Q G_0\cong\Q G_1\cong\Q^8\oplus\Ma_2(\Q)^6\oplus\Ma_2(\Q(i))^2\oplus\Ma_2(\Fsq)^2\oplus\Ma_2(\HHt)^2,
\]
and already $\Q N(3)\cong\Q N(4)\cong\Q^8\oplus\Ma_2(\Q)^2\oplus\Ma_2(\Q(i))^2$. We record three comments on the pairs $A\rtimes(N(3)\times C_2^m)$, $A\rtimes(N(4)\times C_2^m)$:
\begin{itemize}
\item One has $FG_0\cong FG_1$ for every field $F$ of characteristic zero. However $\Z G_0\not\cong\Z G_1$, because $G_0$ and $G_1$ are metabelian and the integral group ring of a metabelian group determines the group by Whitcomb's theorem \cite{Whitcomb}.
\item The character tables of $G_0$ and $G_1$ are not isomorphic, although $\C G_0\cong\C G_1$ and the two groups have the same class sizes. 
\item The unit groups $\U(\Z G_0)$ and $\U(\Z G_1)$ are commensurable, since the unit groups of two orders in the same semisimple $\Q$-algebra are commensurable \cite[Lemma 4.6.9]{EricAngel1}. Hence the virtual structure of the unit group, which is the subject of \cite{CJT1}, does not distinguish $G_0$ from $G_1$.
\end{itemize}
\end{remark}

Recall that $G$ is determined by $|A_{-1}|$, $|A\cap\ZZ(G)|$ and its Sylow pair $(Q,C)$ by \Cref{thm:abstract-iso}. To prove \Cref{rational Iso for Mexc} we show in \Cref{subsec:RIP-proof} that these data can be read off from the Wedderburn--Artin decomposition of $\Q G$, except that the decompositions for the pairs $(N(3),C(3))$ and $(N(4),C(4))$ coincide. Along the way we will obtain that the following invariants are determined.

\begin{proposition}\label{cor:RIP-invariants}
Let $G$ be a non-nilpotent finite group such that $\Q G$ has $\Mexc$. Then $\Q G$ determines
\begin{enumerate}
\item the isomorphism types of $\Fit(G)$, $\ZZ(G)$ and of $G/G'$,
\item the nilpotency class and the exponent of the Sylow $2$-subgroups of $G$.
\end{enumerate}
\end{proposition}

Together with \Cref{sectie only division component} this gives the following picture.

\begin{corollary}\label{cor:RIP-summary}
Let $G, H$ be finite groups such that $\Q G$ has $\Mexc$ and $\Q H\cong\Q G$. Then $G\cong H$ in each of the following cases:
\begin{enumerate}[(a)]
\item $\Q G$ has only division components;
\item $G$ is non-nilpotent and the Sylow $2$-subgroup $C$ of $\Fit(G)$ is abelian;
\item $G$ is non-nilpotent and its Sylow $2$-subgroups have class three;
\item $G$ is non-nilpotent with Sylow $2$-subgroups of class two, and the Sylow pair of $G$ is not $(N(3)\times C_2^m,C(3)\times C_2^m)$ or $(N(4)\times C_2^m,C(4)\times C_2^m)$.
\end{enumerate}
\end{corollary}
\begin{proof}
\Cref{case of prod of div algebras,rational Iso for Mexc}, together with \Cref{fundamental prop det by rational grp alg}, which shows that the hypotheses on $G$ are shared by $H$. The Sylow $2$-subgroups $N(3)\times C_2^m$ and $N(4)\times C_2^m$ belong to the pairs $(N(3)\times C_2^m,C(3)\times C_2^m)$ and $(N(4)\times C_2^m,C(4)\times C_2^m)$ only, by \Cref{thm:class2-nonabelian-list} and the Krull--Remak--Schmidt theorem.
\end{proof}

The proof of \Cref{rational Iso for Mexc} will use the structural results obtained in the earlier sections. For those results to apply we always had the assumption that $\Q G$ has a non-division component. Therefore in \Cref{sectie only division component} we deal separately with the case that $\Q G$ has only division components. Thereafter, \Cref{reductie tot 2-groups} collects the invariants of $G$ which can be read off from $\Q G$, and \Cref{subsec:RIP-proof} contains the proof of \Cref{rational Iso for Mexc} and \Cref{cor:RIP-invariants}.

\subsection{Only division algebra components}\label{sectie only division component}

To start we consider the case that $\Q G$ is the direct product of division algebras. By \cite[Theorem 3.5]{SehNilpotent} this happens if and only if $G$ is of the form $A \times C_2^{n} \times Q_8$ with $A$ an odd abelian group such that $|A|$ and $o_{|A|}(2)$ are odd. 

Let $F$ be a field of characteristic $0$. In \cite[Remark 3.4]{CJT1}, it was noted $FG$ is a direct product of division algebras exactly when $\Q G$ is and $F$ is not a splitting field of $\qa{-1}{-1}{\Q}$. Now, by \cite[Theorem 3.5]{SehNilpotent}, if $G$ is non-abelian, then $\Q G$ is a direct product of division algebras if and only if $G$ is of the form $A \times C_2^{n} \times Q_8$ with $A$ an odd abelian group such that $|A|$ and $o_{|A|}(2)$ are odd.

Consequently, if $FG$ is a direct product of division algebras, then $G$ is an Hamiltonian group (i.e. all subgroups are normal). 

\begin{theorem}\label{case of prod of div algebras}
Let $G,H$ be a finite groups such that $\Q G$ is a direct product of division algebras. If $\Q G \cong \Q H$ as $\Q$-algebras, then $G \cong H$.
\end{theorem}
\begin{proof}
If $G$ is abelian, and hence $\Q G$, then $H$ is also abelian. In that case it is classical that the rational group algebra uniquely determines $G$. Thus we may assume that $G$ is non-abelian. By the discussion at the start of the subsection, $G \cong A \times C_2^{n} \times Q_8$ for some odd abelian group $A$ and $n\in \N$. By the theorem of Perlis-Walker \cite[Theorem 3.3.6]{EricAngel1} we have that $
\Q A \cong \bigoplus_{d \, \mid \, |A|} k_d \Q (\zeta_d)$
where $k_d$ is the number of cyclic subgroups of order $d$ in $A$ (and $\zeta_d$ a $d$-th root of unity). Next, recall that $\Q[Q_8 \times C_2^n] \cong 4.2^n \,\Q \op 2^n  \HH$, since $\Q[Q_8] \cong 4. \Q \op \HH$. Therefore,
$$
\Q G\cong\bigoplus_{d\mid|A|}k_d\cdot4\cdot2^n\,\Q(\zeta_d)\ \oplus\ k_d2^n\qa{-1}{-1}{\Q(\zeta_d)}.
    $$
Now by firstly considering the simple $\qa{-1}{-1}{\Q}$, where $k_d = 1$, we obtain the invariant $2^n = |C_2^n|$. With this at hand we can filter both the values $k_d$ and $d \mid |A|$ from the multiplicities of the factors $\Q(\zeta_d)$. In sum, we recover the isomorphism type of $G$. 

Finally, if $\Q G \cong \Q H$ as $\Q$-algebras, then $\Q H$ is also a direct product of division algebras. Hence the same reasoning as above holds for $H$ and in particular yields that $H$ must have the same factor $A \times C_2^n$, than $G$.
\end{proof}

\begin{remark}
    The proof of \Cref{case of prod of div algebras} also shows that if $\Q G \cong \Q H$ and both $G$ and $H$ are Hamiltonian groups, then $G \cong H$. It would be interesting to know whether Hamiltonian groups can be distinguished among all finite groups. In other words, is it necessary to assume that $H$ is also Hamiltonian?  In \cite[Theorem 5.21]{JTT} it was proven that a group is Hamiltonian if and only if the image of $G$ inside each irreducible representation is a Frobenius complement. This characterisation gives the impression that the property that $H$ is Hamiltonian is not one of the rational group algebra. 
\end{remark}

\subsection{Invariants determined by the rational group algebra}\label{reductie tot 2-groups}
In the following result we already show that $\Q G$ determines the most of the invariants from the decomposition $G = A \rtimes Q$. Recall that, if $G$ has a normal Sylow $p$-subgroup, we denote it by $G_p$.

\begin{proposition}\label{fundamental prop det by rational grp alg}
Let $G$ be a finite group such that $\Q G$ has $\Mexc$ and a non-division component. Then the following are determined by $\Q G$:
\begin{enumerate}[(1)]
\item the isomorphism type of $G^{\mathrm{ab}}$, $G_3$ and $\ZZ(G)_3$,
\item $|G|$ and $|G'|$,
\item the set of character degrees $\cd(G)$,
\item the nilpotency class of the Sylow $2$-subgroups, if $G$ is not nilpotent,
\item the index $[G:\Fit(G)]$, in particular whether $G$ is nilpotent,
\item the smallest index of an abelian normal subgroup, if $G$ is not nilpotent. It is $2$ if $C$ is abelian and $4$ otherwise.
\end{enumerate}
\end{proposition}
The first three parts of \Cref{fundamental prop det by rational grp alg} in fact hold for any $G$ and are well-known, but we could not point out an explicit reference. Hence for ease of the reader we add a brief proof.
\begin{proof}
The set $\cd(G)$ of character degrees is determined by the Wedderburn--Artin decomposition of $\C\otimes_\Q\Q G$ and hence also by the rational group algebra.

For part (1), we decompose $\Q G\cong(\Q G)_{com}\times(\Q G)_{nc}$, where $(\Q G)_{com}$ is the product of the commutative simple components and $(\Q G)_{nc}$ the rest. Then it is well-known that $(\Q G)_{com}\cong\Q[G^{\mathrm{ab}}]$ and $(\Q G)_{nc}=\Q G(1-\widehat{G'})$. Clearly $(\Q G)_{com}$ is determined by $\Q G$. Now by the theorem of Perlis--Walker \cite[Theorem 3.3.6]{EricAngel1}, $\Q[G^{\mathrm{ab}}]\cong\bigoplus_{d\mid|G^{\mathrm{ab}}|}m_d\,\Q(\zeta_d)$, where $m_d$ is the number of cyclic subgroups of $G^{\mathrm{ab}}$ of order $d$. Therefore, by the fundamental theorem of abelian groups, $\Q[G^{\mathrm{ab}}]$ determines the isomorphism type of the abelian group $G^{\mathrm{ab}}$.

The cardinality $|G|=\dim_\Q\Q G$ is determined. Subsequently, $|G'|$ is determined since we have already obtained $|G^{\mathrm{ab}}|$ in (1). This implies that also the number $|G|_3$, i.e.\ the maximal power of $3$ dividing $|G|$, is determined. Now, by \Cref{necessary props 1}, $G_3\cong C_3^{\ell}$ for some $\ell\in\N$. Therefore the cardinality $|G_3|=|G|_3$ also provides the isomorphism type of $G_3$. Next, recall that by \Cref{prop:fixed-inverted} the group $G_3$ has a direct product decomposition $A_{-1}\times(G_3\cap\ZZ(G))$ and $G'\cap G_3=A_{-1}$. Thus the subgroup $\ZZ(G)_3=G_3\cap\ZZ(G)$ is mapped isomorphically into $G^{\mathrm{ab}}$. More precisely, $\ZZ(G)_3$ is isomorphic to the Sylow $3$-subgroup of $G^{\mathrm{ab}}$, which we obtained earlier in the proof.\smallskip

Now we consider the index of the Fitting subgroup. By \Cref{necessary props 1} we have that $[G:\Fit(G)]$ is either $1$ or $2$. Hence it suffices to show that we can determine whether $G$ is nilpotent or not. Recall that $G\cong A\rtimes Q$ with $A$ an elementary abelian $3$-group and $Q$ a $2$-group, and that $G'=A_{-1}\rtimes Q'$. From this we see that $G$ is nilpotent if and only if $A_{-1}=1$, which in turn is equivalent to $G'$ being a $2$-group. The latter can be read off from $|G'|$ and hence is determined by (2).\smallskip

From now on suppose that $G$ is not nilpotent. An abelian normal subgroup of index $2$ of $G$ contains $G'\supseteq A_{-1}$, so it is contained in $\Cen_G(A_{-1})=\Fit(G)=A\times C$ and equals it. Hence $G$ has an abelian normal subgroup of index $2$ if and only if $C$ is abelian, and otherwise the smallest index is $4$ by \Cref{thm:abelian-B}. By \Cref{prop:fixed-inverted}(3), as $\cd(G)\subseteq\{1,2,4\}$ by \Cref{broad form finite grps}, the group $C$ is abelian if and only if $4\notin\cd(G)$, which is determined by (3).\smallskip

Thus it remains to show that the nilpotency class of $Q$ is determined. Since $Q$ has class at most $3$, it is enough to explain how to recognise when $Q$ is abelian or of class $2$. For the former, note that $Q$ is abelian if and only if $G'$ is a $3$-group. This fact is determined by (2). So let $Q$ be non-abelian. We claim that $\cl(Q)=3$ if and only if $\Q G$ has a simple component with centre $\Q(\sqrt2)$ or $\Q(\sqrt{-2})$. If $\cl(Q)=3$, then \Cref{thm:finite-decomp} exhibits such a component. Suppose that $\cl(Q)=2$. Then $\exc(Q)\mid4$ by \Cref{thm:mixed,thm:class2-abelian-kernel,thm:class2-nonabelian-list}, and $\exc(G)\mid12$. So for every strong Shoda pair $(H,K)$ of $G$ the index $[H:K]$ divides $12$, and by \Cref{prop:centre} the centre of the component is a subfield of $\Q(\zeta_{12})$, whose quadratic subfields are $\Q(i)$, $\Q(\sqrt3)$ and $\Q(\sqrt{-3})$.
\end{proof}

In particular, if $H$ is a finite group with $\Q H\cong\Q G$, then $\Q H$ has $\Mexc$ and a non-division component, and the data of \Cref{fundamental prop det by rational grp alg} coincide for $G$ and $H$. If $G$ is not nilpotent, they include $|A_{-1}|=|G_3|/|\ZZ(G)_3|$, $|A\cap\ZZ(G)|=|\ZZ(G)_3|$, $|Q|$, the group $Q/Q'$, which is the Sylow $2$-subgroup of $G^{\mathrm{ab}}$, the class of $Q$, and whether $C$ is abelian.

\subsection{Proof of \texorpdfstring{\Cref{rational Iso for Mexc}}{the main theorem}}\label{subsec:RIP-proof}
We first treat the elementary abelian kernels.

\begin{lemma}[Elementary abelian kernels]\label{lem:elementary-marked}
Let $G$ be a non-nilpotent group of the form \eqref{eq:structure} with $C$ elementary abelian, and let $t\in Q\setminus C$. Then the Sylow pair $(Q,C)$ is determined up to isomorphism by $|C|$, by $|C/[C,t]|$ and by whether $t^2\in[C,t]$. These three data are determined by $\Q G$.
\end{lemma}
\begin{proof}
Write $|C|=2^k$, let $\sigma$ be the automorphism of $C$ induced by $t$, and put $N_\sigma:=\sigma-1\in\operatorname{End}_{\F_2}(C)$ and $c=t^2$. Then $[C,t]=\operatorname{im}N_\sigma$ and $|C/[C,t]|=2^f$ with $f=\dim\ker N_\sigma$. Since $t^2\in C$ and $C$ is abelian, $N_\sigma^2=0$ and $c\in\ker N_\sigma$. Replacing $t$ by $at$, $a\in C$, replaces $c$ by $N_\sigma(a)+c$, so the extension depends on the class of $c$ in $\ker N_\sigma/\operatorname{im}N_\sigma$. A square-zero endomorphism is determined up to similarity by its rank $r=k-f$, and there is a basis $e_1,\dots,e_r,f_1,\dots,f_r,w_1,\dots,w_h$ ($h=2f-k$) with $N_\sigma(f_j)=e_j$, $N_\sigma(e_j)=N_\sigma(w_\ell)=0$. The centraliser of $N_\sigma$ in $\GL(C)$ induces $\GL_h(2)$ on $\ker N_\sigma/\operatorname{im}N_\sigma=\langle w_1,\dots,w_h\rangle$, which is transitive on non-zero vectors. Hence there are exactly two orbits for the class of $c$, according as $c\in[C,t]$ or not, and the presentation $Q=\langle C,t\mid a^t=\sigma^{-1}(a),\ t^2=c\rangle$ determines the Sylow pair. For the last assertion, $|C|=|Q|/2$ and $Q/Q'=Q/[C,t]$ are determined by \Cref{fundamental prop det by rational grp alg}. The group $Q/Q'$ has order $2|C/[C,t]|$, and it is elementary abelian if and only if $t^2\in[C,t]$.
\end{proof}

\begin{proof}[Proof of \Cref{rational Iso for Mexc} and \Cref{cor:RIP-invariants}]
Let $H$ be a finite group with $\Q H\cong\Q G$. By \Cref{fundamental prop det by rational grp alg}, $H$ is of the form \eqref{eq:structure}, is non-nilpotent, and has the same $|A_{-1}|$, $|A\cap\ZZ(G)|$, $|Q|$, $Q/Q'$, the same nilpotency class of the Sylow $2$-subgroup and, when $\cl(Q)=2$, the same answer to the question whether $C$ is abelian. We write $s=s(A_{-1})$, which is thus the same for $G$ and $H$. By \Cref{thm:abstract-iso}, $H\cong G$ as soon as the Sylow pairs are isomorphic, so we have to decide to what extent the decomposition of $\Q G$ determines the Sylow pair, and for this we go through the cases of \Cref{secrtion full characterisation solv}.

Suppose first that $Q$ is abelian. By \Cref{classif Q abelian}, $G\cong D\times(A_{-1}\rtimes\langle y\rangle)$ with $D=E\times(A\cap\ZZ(G))$, and $\Q G$ determines $X:=(G/G')_2\cong E\times C_{2^t}$, where $\ord(y)=2^t$. The pair $(Q,C)=(\langle y\rangle\times E,\langle y^2\rangle\times E)$ is determined by $X$ once $t$ is known. By \Cref{thm:Qabelian-decomp}: if $t=1$ then $\HHt$ does not occur in $\Q G$, whereas for $t\ge2$ it occurs with multiplicity $s\,|\Omega_1(E)|>0$, and if $t\ge3$ then $\exc(E)\mid4$ gives $2^t=\exc(X)$. Hence $t=1$ if and only if $\HHt$ is not a component of $\Q G$, and otherwise $t=\max(2,\log_2\exc(X))$.\smallskip

Suppose next that $\cl(Q)=2$ and that $C$ is abelian. By \Cref{thm:class2-abelian-kernel} the Sylow pair is of type (a), (b) or (c) there. If $C$ is elementary abelian, $\Q G$ has no component $\HH$, whereas by \Cref{thm:finite-decomp} it has one in cases (b) and (c). So the type (a) is recognised, and then \Cref{lem:elementary-marked} recovers the pair. The types (b) and (c) are distinguished by the presence of the field $\Q(i)$, and in both the multiplicity $2^m$ of $\HH$ determines $m$.\smallskip

Suppose that $\cl(Q)=3$. The two families of \Cref{thm:class3} are distinguished by their components: by \Cref{thm:finite-decomp} the semidihedral family contains $\Ma_2(\Q(\sqrt{-2}))$ and the quaternion family the division algebra $\HHrtt$, and $2^m=|Q|/16$ is determined.\smallskip

Suppose finally that $\cl(Q)=2$ and that $C$ is non-abelian. Write $(Q,C)\cong(N_i\times C_2^m,C_i\times C_2^m)$ as in \Cref{thm:class2-nonabelian-list}. In \Cref{thm:finite-decomp} the multiplicity of the field $\Q$ in $\Q G$ is $8\cdot2^{m}$, so $m$ is determined, and dividing all multiplicities by $2^m$ we obtain the row $\Q N_i\oplus(\Fsq C_i\ast N_i/C_i)^s$ of the table. As $s$ is known, the seven algebras $\Q N_i\oplus(\Fsq C_i\ast N_i/C_i)^s$ are pairwise non-isomorphic except for $i=3,4$, where they coincide: the field $\Q(i)$ occurs only for $i\in\{6,7\}$, which are distinguished by $\HH$. The algebra $\HH$ occurs for $i\in\{2,5,6\}$, with multiplicities $2,4,2$, and $i=2$, $i=6$ differ by $\Q(i)$, and among the remaining $i\in\{1,3,4\}$ the multiplicity of $\Ma_2(\Q)$ is $4s$ for $i=1$ and $2+4s$ for $i=3,4$. Hence $i$ is determined unless $i\in\{3,4\}$, in which case $\Q H\cong\Q G$ holds both for $H\cong G$ and for $H\cong A\rtimes(N_{7-i}\times C_2^m)$, as the rows $3$ and $4$ coincide, and these two groups are not isomorphic, since $H/O_3(H)$ is their Sylow $2$-subgroup and $N_3\times C_2^m\not\cong N_4\times C_2^m$ by \Cref{thm:class2-nonabelian-list} and the Krull--Remak--Schmidt theorem.

This proves \Cref{rational Iso for Mexc}. Concerning \Cref{cor:RIP-invariants}, in all cases but the exceptional one the Sylow pair, and hence $G$, is determined. In the exceptional case the two possible pairs have isomorphic kernels $C_3\cong C_4\cong D_8\times C_2$, the groups $N_3$ and $N_4$ have the same abelianisation $C_2^3$, the same centre $C_2^2$ (contained in $C_i$), exponent $4$ and class $2$, as one reads off from the presentations in \Cref{thm:class2-nonabelian-list}, and $\ZZ(G)=(A\cap\ZZ(G))\times(\ZZ(Q)\cap C)$ because an element $aq$ ($a\in A$, $q\in Q$) is central in $G$ if and only if $a\in A\cap\ZZ(G)$, $q$ is central in $Q$ and $q$ centralises $A$.
\end{proof}

\appendix
\clearpage
\section{Table of groups with a faithful exceptional \texorpdfstring{$2 \times 2$}{2 x 2} embedding}\label{appendix section}

In this appendix we reproduce \cite[Table~2]{EKVG} listing those finite groups $G$ that have a faithful exceptional component of type (2) (see \Cref{def_exc_comps}) in the Wedderburn--Artin decomposition of the rational group algebra $\mathbb{Q}G$\footnote{including the group with \textsc{SmallGroupID} \texttt{[24, 1]} that was accidentally omitted in the original table}. For all these groups, we have added extra information which is relevant to the property $\Mexc$:

\vspace{.5cm}

\noindent\begin{tabular}{p{0.25\linewidth}p{0.7\linewidth}}
 $\Mexc$: & indicates whether the group satisfies $\Mexc$\\[.4cm]
 \textsc{SmallGroup ID}: & the identifier of the group $G$ in the small group library \\[.4cm] 
 Structure: & the structure description of the group. \\[.4cm]
 c$\ell$: & the nilpotency class of the group; $\infty$ indicates that the group is not nilpotent (omitted for non-solvable groups) \\[.4cm]
 d$\ell$: & derived length of the group; $\infty$ for non-solvable groups \\[.4cm]
 $\exp(G), \exp(Z(G))$: & the exponent of the group $G$ (resp $Z(G)$, the centre of $G$), i.e. the lowest common multiple of the order of the elements\\ [.4cm]
 $Z(G) \cap G'$: & the structure of the intersetion of the centre of $G$ with the derived subgroup of $G$ \\[.4cm]
 $Q/Z(Q)$: & the structure of the quotient of $Q$ (a Sylow $2$-subgroup of $G$) by its centre, which is independent of the chosen $Q$\\[.4cm]
 Smallest $[G:B]$: & the smallest index of a maximal abelian, normal subgroup $B$ in $G$ \\[.4cm]
 $1 \times 1 $ components: & the division-algebra components appearing in the Wedderburn--Artin decomposition of $\Q G$ (with multiplicity); omitted for groups without $\Mexc$.

\end{tabular}

In the structure descriptions, we use the following conventions: $1$ indicates the trivial group, colons indicate split extensions, a period an extension that is (possibly) non-split and the $n^{\text{th}}$ power of a group indicates its $n$-fold direct product.

Recall that we use the following shorthands for quaternion algebras appearing:
\[\mathbb{H}_2 = \qa{-1}{-1}{\mathbb{Q}}, \quad \mbox{ and } \quad \mathbb{H}_3 = \qa{-1}{-3}{\mathbb{Q}},\] 
and that $\zeta_k$ represents a complex, primitive $k$-root of unity.
\newgeometry{left=1.25cm,right=1.25cm,top=0cm,bottom=0cm}
 \begin{landscape}
 \scriptsize
\begin{longtable}{@{}cllccccllcl@{}} \\ \toprule[1.5pt]
$\Mexc$ & \textsc{SmallGroup} ID & Structure & cl & d$\ell$ & $exp(G)$ & $exp(Z(G))$ & $Z(G)\cap G'
$ & $Q/Z(Q)$ & Smallest $ [G:B] $ & $1\times 1$ components   \\ \midrule 
\endfirsthead \toprule[1.5pt] $\Mexc$ & \textsc{SmallGroup} ID & Structure & cl & d$\ell$ & $exp(G)$ & $exp(Z(G))$ & $Z(G)\cap G'
$ & $Q/Z(Q)$ & Smallest $ [G:B] $ & $1\times 1$ components   \\ \midrule 
\endhead \hline \multicolumn{6}{c}{continued}\\ \midrule[1.5pt]\endfoot\bottomrule[1.5pt]\endlastfoot
$\checkmark$ & \textsf{{[}6, 1{]}} & $ S_3 $ & $\infty$ & $2$ & $6$,  & $1$,  & $1$,  & $1$,  & $2$,  & $2 {\times} \mathbb{Q}$,  \\
$\checkmark$ & \textsf{{[}8, 3{]}} & $ D_8 $ & $2$ & $2$ & $4$,  & $2$,  & $ C_2$,  & $  C_2^2 $,  & $2$,  & $4 {\times} \mathbb{Q}$,  \\
$\checkmark$ & \textsf{{[}12, 4{]}} & $ D_{12} $ & $\infty$ & $2$ & $6$,  & $2$,  & $1$,  & $1$,  & $2$,  & $4 {\times} \mathbb{Q}$,  \\
$\checkmark$ & \textsf{{[}16, 6{]}} & $ C_8 :  C_2$ & $2$ & $2$ & $8$,  & $4$,  & $ C_2$,  & $  C_2^2 $,  & $2$,  & $4 {\times} \mathbb{Q}$, $2 {\times} \mathbb{Q}(i)$,  \\
$\checkmark$ & \textsf{{[}16, 8{]}} & $ {QD}_{16} $ & $3$ & $2$ & $8$,  & $2$,  & $ C_2$,  & $ D_8 $,  & $2$,  & $4 {\times} \mathbb{Q}$,  \\
$\checkmark$ & \textsf{{[}16, 13{]}} & $( C_4  \times   C_2) :  C_2$ & $2$ & $2$ & $4$,  & $4$,  & $ C_2$,  & $  C_2^2 $,  & $2$,  & $8 {\times} \mathbb{Q}$,  \\
$\checkmark$ & \textsf{{[}18, 3{]}} & $ C_3  \times   S_3 $ & $\infty$ & $2$ & $6$,  & $3$,  & $1$,  & $1$,  & $2$,  & $2 {\times} \mathbb{Q}$, $2 {\times} \mathbb{Q}(\sqrt{-3}) = \mathbb{Q}(\zeta_{3})$,  \\
$\checkmark$ & \textsf{{[}24, 1{]}} & $ C_3 :  C_8$ & $\infty$ & $2$ & $24$,  & $4$,  & $1$,  & $1$,  & $2$,  & $2 {\times} \mathbb{Q}$, $1 {\times} \mathbb{Q}(i)$, $1 {\times} \mathbb{Q}(\zeta_{8})$, $1 {\times} \mathbb{H}_3$,  \\
$\times$ & \textsf{{[}24, 3{]}} & $ {\rm SL} (2,3)$ & $\infty$ & $3$ & $12$,  & $2$,  & $ C_2$,  & $  C_2^2 $,  & $12$,  &  \\
$\checkmark$ & \textsf{{[}24, 5{]}} & $ C_4  \times   S_3 $ & $\infty$ & $2$ & $12$,  & $4$,  & $1$,  & $1$,  & $2$,  & $4 {\times} \mathbb{Q}$, $2 {\times} \mathbb{Q}(i)$,  \\
$\checkmark$ & \textsf{{[}24, 8{]}} & $( C_6  \times   C_2) :  C_2$ & $\infty$ & $2$ & $12$,  & $2$,  & $ C_2$,  & $  C_2^2 $,  & $2$,  & $4 {\times} \mathbb{Q}$,  \\
$\checkmark$ & \textsf{{[}24, 10{]}} & $ C_3  \times   D_8 $ & $2$ & $2$ & $12$,  & $6$,  & $ C_2$,  & $  C_2^2 $,  & $2$,  & $4 {\times} \mathbb{Q}$, $4 {\times} \mathbb{Q}(\sqrt{-3}) = \mathbb{Q}(\zeta_{3})$,  \\
$\checkmark$ & \textsf{{[}24, 11{]}} & $ C_3  \times   Q_8 $ & $2$ & $2$ & $12$,  & $6$,  & $ C_2$,  & $  C_2^2 $,  & $2$,  & $4 {\times} \mathbb{Q}$, $4 {\times} \mathbb{Q}(\sqrt{-3}) = \mathbb{Q}(\zeta_{3})$, $1 {\times} \mathbb{H}_2$,  \\
$\checkmark$ & \textsf{{[}32, 8{]}} & $  C_2^2  . ( C_4  \times   C_2)$ & $3$ & $2$ & $8$,  & $2$,  & $ C_2$,  & $( C_4  \times   C_2) :  C_2$,  & $4$,  & $4 {\times} \mathbb{Q}$, $2 {\times} \mathbb{Q}(i)$,  \\
$\checkmark$ & \textsf{{[}32, 11{]}} & $ C_4^2  :  C_2$ & $3$ & $2$ & $8$,  & $4$,  & $ C_2$,  & $ D_8 $,  & $2$,  & $4 {\times} \mathbb{Q}$, $2 {\times} \mathbb{Q}(i)$,  \\
$\checkmark$ & \textsf{{[}32, 44{]}} & $( C_2  \times   Q_8 ) :  C_2$ & $3$ & $2$ & $8$,  & $2$,  & $ C_2$,  & $ C_2  \times   D_8 $,  & $4$,  & $8 {\times} \mathbb{Q}$,  \\
$\checkmark$ & \textsf{{[}32, 50{]}} & $( C_2  \times   Q_8 ) :  C_2$ & $2$ & $2$ & $4$,  & $2$,  & $ C_2$,  & $  C_2^4 $,  & $4$,  & $16 {\times} \mathbb{Q}$,  \\
$\checkmark$ & \textsf{{[}36, 6{]}} & $ C_3  \times  ( C_3 :  C_4)$ & $\infty$ & $2$ & $12$,  & $6$,  & $1$,  & $1$,  & $2$,  & $2 {\times} \mathbb{Q}$, $1 {\times} \mathbb{Q}(i)$, $2 {\times} \mathbb{Q}(\sqrt{-3}) = \mathbb{Q}(\zeta_{3})$, $1 {\times} \mathbb{Q}(\zeta_{12})$, $1 {\times} \mathbb{H}_3$,  \\
$\checkmark$ & \textsf{{[}36, 12{]}} & $ C_6  \times   S_3 $ & $\infty$ & $2$ & $6$,  & $6$,  & $1$,  & $1$,  & $2$,  & $4 {\times} \mathbb{Q}$, $4 {\times} \mathbb{Q}(\sqrt{-3}) = \mathbb{Q}(\zeta_{3})$,  \\
$\times$ & \textsf{{[}40, 3{]}} & $ C_5 :  C_8$ & $\infty$ & $2$ & $40$,  & $2$,  & $1$,  & $1$,  & $4$,  &  \\
$\checkmark$ & \textsf{{[}48, 16{]}} & $( C_3 :  Q_8 ) :  C_2$ & $\infty$ & $2$ & $24$,  & $2$,  & $ C_2$,  & $ D_8 $,  & $4$,  & $4 {\times} \mathbb{Q}$,  \\
$\checkmark$ & \textsf{{[}48, 18{]}} & $ C_3 :  Q_{16} $ & $\infty$ & $2$ & $24$,  & $2$,  & $ C_2$,  & $ D_8 $,  & $4$,  & $4 {\times} \mathbb{Q}$, $1 {\times} \left(\frac{-1, -1}{\mathbb{Q}(\sqrt{2})}\right)$,  \\
$\times$ & \textsf{{[}48, 28{]}} & $  {\rm SL} (2,3) .  C_2$ & $\infty$ & $4$ & $24$,  & $2$,  & $ C_2$,  & $ D_8 $,  & $24$,  &  \\
$\times$ & \textsf{{[}48, 29{]}} & $ {\rm GL} (2,3)$ & $\infty$ & $4$ & $24$,  & $2$,  & $ C_2$,  & $ D_8 $,  & $24$,  &  \\
$\times$ & \textsf{{[}48, 33{]}} & $(( C_4  \times   C_2) :  C_2) :  C_3$ & $\infty$ & $3$ & $12$,  & $4$,  & $ C_2$,  & $  C_2^2 $,  & $12$,  &  \\
$\checkmark$ & \textsf{{[}48, 39{]}} & $( C_4  \times   S_3 ) :  C_2$ & $\infty$ & $2$ & $12$,  & $2$,  & $ C_2$,  & $  C_2^2 $,  & $4$,  & $8 {\times} \mathbb{Q}$,  \\
$\checkmark$ & \textsf{{[}48, 40{]}} & $ Q_8   \times   S_3 $ & $\infty$ & $2$ & $12$,  & $2$,  & $ C_2$,  & $  C_2^2 $,  & $4$,  & $8 {\times} \mathbb{Q}$, $2 {\times} \mathbb{H}_2$,  \\
$\times$ & \textsf{{[}64, 37{]}} & $ ( C_4  \times   C_2) . ( C_4  \times   C_2)$ & $4$ & $2$ & $8$,  & $2$,  & $ C_2$,  & $ C_2^3  :  C_4$,  & $4$,  &  \\
$\checkmark$ & \textsf{{[}64, 137{]}} & $( C_4 :  Q_8 ) :  C_2$ & $3$ & $2$ & $8$,  & $2$,  & $ C_2$,  & $ C_2^4  :  C_2$,  & $4$,  & $8 {\times} \mathbb{Q}$,  \\
$\times$ & \textsf{{[}72, 19{]}} & $ C_3^2  :  C_8$ & $\infty$ & $2$ & $24$,  & $2$,  & $1$,  & $1$,  & $4$,  &  \\
$\times$ & \textsf{{[}72, 20{]}} & $( C_3 :  C_4)  \times   S_3 $ & $\infty$ & $2$ & $12$,  & $2$,  & $1$,  & $1$,  & $4$,  &  \\
$\times$ & \textsf{{[}72, 22{]}} & $( C_6  \times   S_3 ) :  C_2$ & $\infty$ & $2$ & $12$,  & $2$,  & $ C_2$,  & $  C_2^2 $,  & $4$,  &  \\
$\times$ & \textsf{{[}72, 24{]}} & $ C_3^2  :  Q_8 $ & $\infty$ & $2$ & $12$,  & $2$,  & $ C_2$,  & $  C_2^2 $,  & $4$,  &  \\
$\times$ & \textsf{{[}72, 25{]}} & $ C_3  \times   {\rm SL} (2,3)$ & $\infty$ & $3$ & $12$,  & $6$,  & $ C_2$,  & $  C_2^2 $,  & $12$,  &  \\
$\checkmark$ & \textsf{{[}72, 30{]}} & $ C_3  \times  (( C_6  \times   C_2) :  C_2)$ & $\infty$ & $2$ & $12$,  & $6$,  & $ C_2$,  & $  C_2^2 $,  & $2$,  & $4 {\times} \mathbb{Q}$, $4 {\times} \mathbb{Q}(\sqrt{-3}) = \mathbb{Q}(\zeta_{3})$,  \\
$\times$ & \textsf{{[}96, 67{]}} & $ {\rm SL} (2,3) :  C_4$ & $\infty$ & $4$ & $24$,  & $4$,  & $ C_2$,  & $ D_8 $,  & $24$,  &  \\
$\times$ & \textsf{{[}96, 190{]}} & $( C_2  \times   {\rm SL} (2,3)) :  C_2$ & $\infty$ & $4$ & $24$,  & $2$,  & $ C_2$,  & $ C_2  \times   D_8 $,  & $24$,  &  \\
$\times$ & \textsf{{[}96, 191{]}} & $  {\rm SL} (2,3) .  C_2) :  C_2$ & $\infty$ & $4$ & $24$,  & $2$,  & $ C_2$,  & $ C_2  \times   D_8 $,  & $24$,  &  \\
$\times$ & \textsf{{[}96, 202{]}} & $(( C_2  \times   Q_8 ) :  C_2) :  C_3$ & $\infty$ & $3$ & $12$,  & $2$,  & $ C_2$,  & $  C_2^4 $,  & $24$,  &  \\
$\times$ & \textsf{{[}120, 5{]}} & $ {\rm SL} (2,5)$ &  & $\infty$ & $60$,  & $2$,  & $ C_2$,  & $  C_2^2 $,  & $60$,  &  \\
$\times$ & \textsf{{[}128, 937{]}} & $( Q_8   \times   Q_8 ) :  C_2$ & $4$ & $3$ & $8$,  & $2$,  & $ C_2$,  & $( C_2^4  :  C_2) :  C_2$,  & $8$,  &  \\
$\times$ & \textsf{{[}144, 124{]}} & $  {\rm SL} (2,3) .  C_2)$ & $\infty$ & $4$ & $24$,  & $2$,  & $ C_2$,  & $ D_8 $,  & $24$,  &  \\
$\times$ & \textsf{{[}144, 128{]}} & $ S_3   \times   {\rm SL} (2,3)$ & $\infty$ & $3$ & $12$,  & $2$,  & $ C_2$,  & $  C_2^2 $,  & $24$,  &  \\
$\times$ & \textsf{{[}144, 135{]}} & $ C_3^2  : ( C_8 :  C_2)$ & $\infty$ & $2$ & $24$,  & $2$,  & $ C_2$,  & $  C_2^2 $,  & $4$,  &  \\
$\times$ & \textsf{{[}144, 148{]}} & $ C_3^2  : (( C_4  \times   C_2) :  C_2)$ & $\infty$ & $2$ & $12$,  & $2$,  & $ C_2$,  & $  C_2^2 $,  & $4$,  &  \\
$\times$ & \textsf{{[}160, 199{]}} & $(( C_2  \times   Q_8 ) :  C_2) :  C_5$ & $\infty$ & $3$ & $20$,  & $2$,  & $ C_2$,  & $  C_2^4 $,  & $80$,  &  \\
$\times$ & \textsf{{[}192, 989{]}} & $( {\rm SL} (2,3) :  C_4) :  C_2$ & $\infty$ & $4$ & $24$,  & $2$,  & $ C_2$,  & $ C_2^4  :  C_2$,  & $48$,  &  \\
$\times$ & \textsf{{[}240, 89{]}} & $  {\rm SL} (2,5) .  C_2$ &  & $\infty$ & $120$,  & $2$,  & $ C_2$,  & $ D_8 $,  & $120$,  &  \\
$\times$ & \textsf{{[}240, 90{]}} & $ {\rm SL} (2,5) :  C_2$ &  & $\infty$ & $120$,  & $2$,  & $ C_2$,  & $ D_8 $,  & $120$,  &  \\
$\times$ & \textsf{{[}288, 389{]}} & $ C_3^2  : ( C_4^2  :  C_2)$ & $\infty$ & $3$ & $24$,  & $2$,  & $ C_2$,  & $ D_8 $,  & $8$,  &  \\
$\times$ & \textsf{{[}320, 1581{]}} & $ ((( C_2  \times   Q_8 ) :  C_2) :  C_5) .  C_2$ & $\infty$ & $4$ & $40$,  & $2$,  & $ C_2$,  & $ C_2^4  :  C_2$,  & $160$,  &  \\
$\times$ & \textsf{{[}384, 618{]}} & $(( Q_8   \times   Q_8 ) :  C_2) :  C_3$ & $\infty$ & $3$ & $24$,  & $2$,  & $ C_2$,  & $( C_2^4  :  C_2) :  C_2$,  & $96$,  &  \\
$\times$ & \textsf{{[}384, 18130{]}} & $(( Q_8   \times   Q_8 ) :  C_3) :  C_2$ & $\infty$ & $4$ & $24$,  & $2$,  & $ C_2$,  & $( C_2^4  :  C_2) :  C_2$,  & $96$,  &  \\
$\times$ & \textsf{{[}720, 409{]}} & $ {\rm SL} (2,9)$ &  & $\infty$ & $120$,  & $2$,  & $ C_2$,  & $ D_8 $,  & $360$,  &  \\
$\times$ & \textsf{{[}1152, 155468{]}} & $((( Q_8   \times   Q_8 ) :  C_3) :  C_2) :  C_3$ & $\infty$ & $4$ & $24$,  & $2$,  & $ C_2$,  & $( C_2^4  :  C_2) :  C_2$,  & $288$,  &  \\
$\times$ & \textsf{{[}1920, 241003{]}} & $ C_2 . ( C_2^4  :  A_5 )$ &  & $\infty$ & $120$,  & $2$,  & $ C_2$,  & $( C_2^4  :  C_2) :  C_2$,  & $960$,  &  \\
\end{longtable}
	\end{landscape}
\restoregeometry

\section{Comparison with the groups of Kleinian type}\label{appendix section B}
In this appendix we explain how the finite groups classified in \cite{FreebyFree} fit in the classification of \Cref{secrtion full characterisation solv}. Recall that a finite group $G$ is said to be of \emph{Kleinian type} if every non-commutative simple component of $\Q G$ is a totally definite quaternion algebra or isomorphic to $\Ma_2(K)$ with $K\in\{\Q,\Q(i),\Q(\sqrt{-2}),\Q(\sqrt{-3})\}$. In particular $\Q G$ then has $\Mexc$. We first recall the main result of \cite{FreebyFree}, restricted to non-nilpotent groups. It involves the following two $2$-groups, where $n\ge1$:
\begin{align*}
\mc{W}_{1n}&:=\langle x,y_1,\dots,y_n\mid x^4=y_i^2=[y_i,y_j]=1,\ [y_i,x]\text{ central of order }2\rangle\cong C_2^{2n}\rtimes C_4,\\
\mc{W}_{21}&:=\langle x,y\mid x^4=y^4=1,\ y^x=y^{-1}\rangle\cong C_4\rtimes C_4.
\end{align*}

\begin{theorem}[{\cite[Theorem 1]{FreebyFree}}]\label{thm:JPdRRZ}
Let $G$ be a finite group. Then the following are equivalent:
\begin{enumerate}[(1)]
\item $\U(\Z G)$ is commensurable with a direct product of free-by-free groups,
\item $G$ is of Kleinian type,
\item $G$ is abelian or an epimorphic image of $E\times H$ with $E$ an abelian group and $H$ a group from an explicit list.
\end{enumerate}
Moreover, if $G$ is not nilpotent, then in \textup{(3)} one has $H=M\rtimes P$ with $M$ a non-trivial elementary abelian $3$-group and $P$ a $2$-group such that $P_0:=\Cen_P(M)$ has index $2$ in $P$ and the elements of $P\setminus P_0$ invert $M$. Furthermore one of the following holds:
\begin{enumerate}[(a)]
\item $\exc(E)=4$ and $P=C_8$,
\item $\exc(E)=6$, $P=\mc{W}_{1n}$ and $P_0=\langle x^2,y_1,\dots,y_n\rangle P'$,
\item $\exc(E)=2$, $P=\mc{W}_{21}$ and $P_0=\langle x,y^2\rangle$.
\end{enumerate}
\end{theorem}

For nilpotent $G$ the list in \textup{(3)} consists of twelve families of $2$-groups of class at most $3$. They are related to the nilpotent groups with $\Mexc$ in \cite{CJT3}.

Now let $G=A\rtimes Q$ be a finite non-nilpotent group such that $\Q G$ has $\Mexc$, with Sylow pair $(Q,C)$. The following describes when $G$ is of Kleinian type.

\begin{proposition}\label{prop:kleinian-comparison}
Let $G$ be a finite non-nilpotent group such that $\Q G$ has $\Mexc$. Then the following hold:
\begin{enumerate}[(1)]
\item The groups $\Fit(G)$ and $Q$ are of Kleinian type.
\item $G$ is of Kleinian type if and only if $C$ is abelian and $16\nmid\exc(G)$.
\end{enumerate}
\end{proposition}
\begin{proof}
By \Cref{thm componenten QG mexc}, a group with $\Mexc$ is of Kleinian type if and only if its rational group algebra has no simple component of degree $4$ and no exceptional division component. Statement (1) therefore follows from \Cref{cor:index2} for $\Fit(G)$ and from \Cref{cor:which-Q} for $Q$.

For (2), suppose first that $C$ is abelian. Then $\Fit(G)=A\times C$ is an abelian subgroup of index $2$, hence $\Q G$ has no component of degree $4$. Moreover, by \Cref{exc div comp for mexc} the algebra $\Q G$ has an exceptional division component if and only if $Q$ is abelian and $16\mid\exc(G)$. If $Q$ is non-abelian, then $\exc(Q)\mid4$ by \Cref{thm:mixed,thm:class2-abelian-kernel}. Thus $G$ is of Kleinian type exactly when $16\nmid\exc(G)$. Suppose now that $C$ is non-abelian. Take an irreducible complex character $\psi$ of $C$ of degree $2$ and a non-trivial linear character $\lambda$ of $A_{-1}$. Then the character $\lambda\times\psi$ of $\Fit(G)$ is not $G$-invariant, as the elements of $Q\setminus C$ invert $A_{-1}$. Hence it induces to an irreducible character of $G$ of degree $4$ and $G$ is not of Kleinian type.
\end{proof}

Consequently, the non-nilpotent groups of Kleinian type are exactly the groups of \Cref{thm:mixed}, of \Cref{thm:class2-abelian-kernel} and those of \Cref{classif Q abelian} with $y$ of order at most $8$. The following table records for each of the Sylow pairs of \Cref{Main th B} whether the associated groups are of Kleinian type and, if so, the corresponding case of \Cref{thm:JPdRRZ}.

\begin{center}
\renewcommand{\arraystretch}{1.25}\small
\begin{tabular}{@{}lllc@{}}
\toprule
Sylow pair $(Q,C)$ & Theorem & Kleinian type & case\\
\midrule
$C$ elementary abelian & \ref{thm:mixed}, \ref{thm:class2-abelian-kernel} & yes & (b)\\
$(\langle y\rangle\times E,\ \langle y^2\rangle\times E)$ with $\exc(E)\mid4$ & \ref{classif Q abelian} & yes if $\ord(y)\mid 8$, otherwise no & (a)\\
$(Q_8\times C_2^m,\ C_4\times C_2^m)$, $(R\times C_2^m,\ C_R\times C_2^m)$ & \ref{thm:class2-abelian-kernel} & yes & (c)\\
$(N_i\times C_2^m,\ C_i\times C_2^m)$ with $1\le i\le7$ & \ref{thm:class2-nonabelian-list} & no & --\\
$(SD_{16}\times C_2^m,\ D_8\times C_2^m)$, $(Q_{16}\times C_2^m,\ Q_8\times C_2^m)$ & \ref{thm:class3} & no & --\\
\bottomrule
\end{tabular}
\end{center}

\noindent We now explain the last column. Note that if $G$ is an epimorphic image of $E\times(M\rtimes P)$ as in \Cref{thm:JPdRRZ}, then the Sylow pair of $G$ is an epimorphic image of the pair $(E_2\times P,E_2\times P_0)$, where $E_2$ is the Sylow $2$-subgroup of $E$. Moreover $3\mid|\ZZ(G)|$ is only possible when $3\mid\exc(E)$, thus in case (b).
\begin{itemize}
\item In case (a) the group $E_2\times P$ is abelian with $\exc(E_2)\mid 4$ and $P=\langle y\rangle$ cyclic of order $8$. By \Cref{classif Q abelian} its epimorphic images which are Sylow pairs are the pairs $(\langle y\rangle\times E,\langle y^2\rangle\times E)$ with $\ord(y)\mid8$.
\item In case (b) the group $E_2\times P_0$ is elementary abelian. Conversely, let $(Q,C)$ be a Sylow pair with $C$ elementary abelian and $t\in Q\setminus C$. Then $C$ is a module over the group algebra $\F_2[\langle t\rangle/\langle t^2\rangle]$, hence an epimorphic image of a direct sum of a free module and a trivial module. In $\mc{W}_{1n}$ the elements $y_i$ and $[y_i,x]$ span a free module of rank $n$ under the action of $x$. Therefore $(Q,C)$ is an epimorphic image of $(E_2\times\mc{W}_{1n},E_2\times P_0)$ for $n$ and $E_2$ large enough, with $x$ mapped to $t$. The condition $\exc(E)=6$ reflects that this is the only case in which $A\cap\ZZ(G)$ can be non-trivial, in accordance with \Cref{thm:mixed}.
\item In case (c) the group $\mc{W}_{21}$ is the group $R$ of \eqref{eq:R}, with $y=r$ and $x=v$. The automorphism of $R$ defined by $r\mapsto r$ and $v\mapsto rv$ maps $P_0=\langle v,r^2\rangle$ onto $C_R=\langle rv,r^2\rangle$. Hence the pair $(E_2\times P,E_2\times P_0)$ is isomorphic to $(R\times C_2^m,C_R\times C_2^m)$. Furthermore $(Q_8,C_4)$ is the image of $(R,C_R)$ modulo $\langle r^2v^2\rangle$. The remaining epimorphic images of $(R\times C_2^m,C_R\times C_2^m)$ are pairs $(Q,C)$ with $Q$ abelian or $C$ elementary abelian, which are covered by the cases (a) and (b).
\end{itemize}
Altogether the classification of the non-nilpotent groups in \cite{FreebyFree} is recovered by the first three rows of the table. Note however that \Cref{thm:JPdRRZ} describes these groups as epimorphic images, whereas \Cref{thm:mixed,classif Q abelian,thm:class2-abelian-kernel} describe the Sylow pairs up to isomorphism. The groups with $\Mexc$ which are new compared to \cite{FreebyFree} are those of the last two rows, which have simple components of the form $\Ma_2(\HH)$ or $\Ma_2(\HHt)$ by \Cref{thm:finite-decomp}, and those of the second row with $\ord(y)\ge16$, which have the exceptional division components $\qa{\zeta_{2^k}}{-3}{\Q(\zeta_{2^k})}$ with $k\ge3$.

\bibliographystyle{plain}
\bibliography{LowDeg}

\end{document}